\documentclass[11pt]{article}

\usepackage{amsfonts,amssymb}
\usepackage[all, knot]{xy}
\usepackage{diagram}

\newcommand{\lt}{\triangleleft}

\newcommand{\C}{\mathbb{C}}
\newcommand{\Z}{\mathbb{Z}}

\def\g{{\mathfrak g}}

\newcommand{\Hom}{\mbox{Hom}}

\newcommand{\Set}{{\rm Set}}
\newcommand{\Vect}{{\rm Vect}}
\newcommand{\Coalg}{{\rm Coalg}}

\newcommand{\CoComCoalg}{{\rm CoComCoalg}}

\newtheorem{theorem}{Theorem}[section]
\newtheorem{corollary}[theorem]{Corollary}
\newtheorem{lemma}[theorem]{Lemma}
\newtheorem{proposition}[theorem]{Proposition}

\newtheorem{definition}[theorem]{Definition}

\newtheorem{example}[theorem]{Example}

\newtheorem{remark}[theorem]{Remark}

\input{epsf.sty}

\def\g{{\mathfrak g}}

\begin{document}

\title{Cohomology of Categorical Self-Distributivity   }

\author{
J. Scott Carter\footnote{Supported in part by NSF Grant DMS \#0301095.}
\\ University of South Alabama
\and
Alissa S. Crans \\
Loyola Marymount University \and
Mohamed
Elhamdadi
\\ University of South Florida
\and
Masahico Saito\footnote{Supported in part by NSF Grant DMS \#0301089.}
\\ University of South Florida
}

\maketitle

\begin{abstract}
We define self-distributive structures in the categories of
coalgebras and cocommutative coalgebras.  We obtain examples from
vector spaces whose bases are the elements of finite quandles, the
direct sum of a Lie algebra with its ground field, and Hopf
algebras.  The self-distributive operations of these structures
provide solutions of the Yang--Baxter equation, and, conversely,
solutions of the Yang--Baxter equation can be used to construct
self-distributive operations in certain categories.

Moreover, we present a cohomology theory that encompasses both Lie
algebra and quandle cohomologies, is analogous to Hochschild
cohomology, and can be used to study deformations of these
self-distributive structures.  All of the work here is informed
via diagrammatic computations.

\end{abstract}

\section{Introduction}
In the past several decades, operations satisfying
self-distributivity [$(a\lt b)\lt c= (a\lt c) \lt (b \lt c)$] have
secured an important role in knot theory. Such operations not only
provide solutions of the Yang--Baxter equation and satisfy a law
that is an algebraic distillation of the type (III) Reidemeister
move, but they also capture one of the essential properties of
group conjugation. Sets possessing such a binary operation are
called shelves. Adding an axiom corresponding to the type (II)
Reidemeister move amounts to the property that the set acts on
itself (on the right) bijectively and thus gives the structure of
a rack.  Further introducing a condition corresponding to the type
(I) Reidemeister move has the effect of making each element
idempotent and gives the structure of a quandle.  Keis, or
involutory quandles, satisfy an extra involutory condition.  Such
structures were discussed as early as the 1940s \cite{Takasaki}.

The primordial example of a self-distributive operation comes from
group conjugation: $x \lt y = y^{-1} x y$.
This operation satisfies the additional quandle axioms which are
stated in the sequel. Quandle cohomology has been studied
extensively in connection with applications to knots and knotted
surfaces~\cite{CJKLS,CKamS}.
Analogues of self-distributivity in a variety of categorical
settings have been discussed as adjoint maps in Lie algebras
~\cite{Alissa}
 and
quantum group theories (see for example
\cite{MajidPink,MajidGreen}). In particular, the adjoint map $x
\otimes y \mapsto S(y_{(1)})x y_{(2)}$ of Hopf algebras is a
direct analogue of group conjugation. Thus, analogues of
self-distributive operations are found in a variety of algebraic
structures where cohomology theories are also defined.

In this paper, we study how quandles and racks and their
cohomology theories are related to these other algebraic systems
and their cohomology theories. Specifically, we treat
self-distributive maps in a unified manner via a categorical
technique called internalization ~\cite{E}. Then we develop a
cohomology theory
and provide explicit relations to rack and Lie algebra cohomology
theories.
Furthermore, this
cohomology theory can be seen as a theory of obstructions to
deformations of self-distributive structures.

The organization of this paper is as follows:
Section~\ref{qcatsec} consists of a review of the fundamentals of
quandle theory, internalization in a category, and the definition
of a coalgebra. Section~\ref{coalsec} contains a collection of
examples  that possess a self-distributive binary operation. In
particular, a motivating example built from a Lie algebra is
presented. In Section~\ref{ybesec} we relate the ideas of
self-distributivity to solutions of the Yang-Baxter equation, and
demonstrate connections of these ideas to Hopf algebras.
Section~\ref{hochsec} contains  a review of Hochschild cohomology
from the diagrammatic point of view and in relation to
deformations of algebras. These ideas are imitated in
Section~\ref{cohsec} where the most original and substantial ideas
are presented. Herein a cohomology theory for shelves in the
coalgebra category
 is defined in low dimensions.
The theory is informed by the diagrammatic representation
of the self-distributive operation, the comultiplication, their axioms,
and their relationships.
Section~\ref{relsec} contains the main results of the paper. Theorems~\ref{q3prop}
through
\ref{nontrivLie2prop}
state that the cohomology theory is non-trivial, and that  non-trivial quandle cocycles and Lie algebra cocycles give non-trivial shelf cocycles in dimension $2$ and $3$.

\subsection*{Acknowledgements} In addition to the National Science Foundation who supported this research financially, we also wish to thank our colleagues whom
we engaged in a number of crucial discussions. The topology
seminar at South Alabama listened to a series of talks as the work
was being developed. Joerg Feldvoss gave two of us a wonderful
lecture on deformation theory of algebras and helped provide a key
example. John Baez was extremely helpful with some fundamentals of
categorical constructions. We also thank N. Apostolakis, L.H.
Kauffman, and D. Radford for valuable conversations.

\section{Internalized Shelves}\label{qcatsec}

\subsection{Review of Quandles}

A {\it quandle}, $X$, is a set with a binary operation $(a, b)
\mapsto a \lt b$ such that

(I) For any $a \in X$,
$a\lt a =a$.

(II) For any $a,b \in X$, there is a unique $c \in X$ such that
$a= c\lt b$.

(III)
For any $a,b,c \in X$, we have
$ (a \lt b) \lt c=(a\lt c)\lt (b\lt c). $

\noindent A {\it rack} is a set with a binary operation that
satisfies (II) and (III). Racks and quandles have been studied
extensively in, for example, \cite{Br88,FR,Joyce,Matveev}.

The following are typical examples of quandles: A group $G$ with
conjugation as the quandle operation: $a \lt b = b^{-1} a b$,
denoted by $X=$ Conj$(G)$, is a quandle. Any subset of $G$ that is
closed under such conjugation is also a quandle. More generally if
$G$ is a group, $H$ is a subgroup, and $s$ is an automorphism that
fixes the elements of $H$  ({\it i.e.} $s(h)=h \ \forall h \in
H$), then $G/H$ is a quandle with $\lt $ defined by $Ha\lt
Hb=Hs(ab^{-1})b.$ Any $\Lambda (={\Z }[t, t^{-1}])$-module $M$ is
a quandle with $a\lt b=ta+(1-t)b$, for $a,b \in M$, and is called
an {\it  Alexander  quandle}. Let $n$ be a positive integer, and
for elements  $i, j \in \{ 0, 1, \ldots , n-1 \}$, define $i\lt j
\equiv 2j-i \pmod{n}$. Then $\lt$ defines a quandle structure
called the {\it dihedral quandle},
  $R_n$, that
 coincides with the set  of reflections in the dihedral group
with composition given by conjugation.

The third quandle axiom $(a\lt b)\lt c=(a\lt c)\lt (b\lt c)$,
which corresponds to the type (III) Reidemeister move, can be
reformulated to make sense in a more general setting.  In fact, we
do not need the full-fledged structure of a quandle; we simply
need a structure having a binary operation satisfying the
self-distributive law.  We call a set together with a binary
operation satisfying the self-distributive axiom (III) a {\it
shelf}.

We reformulate the self-distributive operation of a shelf as
follows: Let $X$ be a shelf with the shelf operation denoted by a
map $q: X \times X \rightarrow X $. Define $\Delta: X \rightarrow
X \times X $ by $\Delta(x)=(x, x)$ for any $x \in X$, and $\tau: X
\times X \rightarrow X \times X $ by a transposition $\tau
(x,y)=(y,x)$ for $x, y \in X$. Then axiom (III) above can be
written as:
$$ q (q \times 1 ) = q (q \times q) ( 1 \times \tau \times 1)
 (1 \times 1 \times \Delta)
: X^3 \rightarrow X . $$

It is natural and useful to  formulate this axiom for morphisms in
certain categories. This approach was explored in \cite{Alissa}
(see also \cite{BC}) and involves a technique known as
internalization.

\subsection{Internalization}

All familiar mathematical concepts were defined in the category of
sets, but most of these can live in other categories as well. This
idea, known as internalization,
 is actually very familiar.  For
example, the notion of a group can be enhanced by looking at
groups in categories other than $\Set$, the category of sets and
functions between them. We have the notions of topological groups,
which are groups in the category of topological spaces, Lie
groups, groups in the category of smooth manifolds, and so on.
Internalizing a concept consists of first expressing it completely
in terms of commutative diagrams and then interpreting those
diagrams in some sufficiently nice ambient category, $K$. In this
paper, we consider the notion of a shelf in the categories of
coalgebras and cocommutative coalgebras. Thus, we define the
notion of an internalized shelf, or shelf in $K$. This concept is
also known as a shelf object in $K$ or internal shelf.

Given two objects $X$ and $Y$ in an arbitrary category, we define
their {\it product} to be any object $X \times Y$ equipped with
morphisms $\pi_1 : X \times Y \rightarrow X$ and $\pi_2: X \times
Y \rightarrow Y$ called projections, such that the following
universal property is satisfied: for any object $Z$ and morphisms
$f : Z \rightarrow X$ and $g: Z \rightarrow Y,$ there is a unique
morphism $h: Z \rightarrow X \times Y$ such that $f = \pi_1 h$ and
$g = \pi_2 h$. Note that this product does not necessarily exist,
nor is it unique. However, it is unique up to canonical
isomorphism, which is why we refer to {\it the} product when it
exists.  We say a category has {\it binary products} when every
pair of objects has a product. {\it Trinary products} $(X \times
Y) \times Z$ and $X \times (Y \times Z)$ are defined similarly,
are canonically isomorphic, and denoted  by $X \times Y \times Z$
if the isomorphism is the identity. Inductively, {\it $n$-ary}
products are defined. We say a category has {\it finite products}
if it has $n$-ary products for all $n \geq 0$. Note that whenever
$X$ is an object in some category for which the product $X \times
X$ exists, there is a unique morphism called the {\it diagonal}
$D: X \rightarrow X \times X$ such that $\pi_1 D = 1_X$ and $\pi_2
D = 1_X$. In the category of sets, this map is given by $D(x) =
(x,x)$ for all $x \in X$. In a category with finite products, we
also have a {\it transposition} morphism given by
 $\tau: X \times X \rightarrow  X \times X$
by $\tau = (\pi_2 \times \pi_1) D_{X \times X}$.
 %%%% An old def kept:
 %to be a  morphism such that
%$\pi_1 \tau = \pi_2$ and
%$\pi_2 \tau = \pi_1$. It exists uniquely from the universal
%property, and is
%viewed as a map that switches inputs.

\begin{definition} \label{selfdist} {\rm
Let $X$ be an object in a category $K$ with finite products.
A map $q: X \times X \rightarrow X$ is a {\it self-distributive
map} if the following diagram commutes:

\[\begin{xy}
    (-2, 20)*+{X \times X \times X}="1";
    (-30,10)*+{X \times X \times X \times X}="2";
    (-30,-10)*+{X \times X \times X \times X}="3";
    (-15,-20)*+{X \times X \times X}="4";
    (15,-20)*+{X \times X}="4a";
    (30,-10)*+{X}="5";
    (30, 10)*+{X \times X}="6";
        {\ar^{q \times 1} "1";"6"};
        {\ar_{1 \times 1 \times \Delta} "1";"2"};
        {\ar_{1 \times \tau \times 1} "2";"3"};
        {\ar_{1 \times 1 \times q} "3";"4"};
        {\ar_{q \times 1} "4";"4a"};
        {\ar_{q} "4a";"5"};
        {\ar^{q} "6";"5"};
\end{xy}
\]
where $\Delta : X \rightarrow X \times X$ is
the
diagonal morphism in $K$ and $\tau : X \times X \rightarrow X
\times X$ is
the transposition.
We also say that a map $q$ satisfies the {\it self-distributive law}.

} \end{definition}

\begin{definition} \label{rshelfobj}
{\rm Let $K$ be a category with finite products. A  {\it shelf in
$K$} is a pair $(X, q)$ such that $X$ is an object in $K$ and $q :
X \times X \rightarrow X$ is
a morphism in $K$ that satisfies the self-distributive law of Definition~\ref{selfdist}.

} \end{definition}

\begin{example}{\rm
A quandle $(X, q)$ is a shelf in the category of sets, with
the cartesian products and the diagonal map $D: X \rightarrow X \times X$
defined by $D(x)=(x, x)$ for all $x \in X$.
Thus the language of shelves and self-distributive maps in categories
unifies all examples discussed in this paper, in particular those constructed
from Lie algebras.
} \end{example}

\begin{remark}\label{categories}{\rm
Throughout this paper, all of the categories considered have
finite products:}

\begin{itemize}
\setlength{\itemsep}{-3pt}

\item $\Set$, the category whose objects are sets and whose
morphisms are functions.

\item $\Vect$, the category whose objects are vector spaces
over a field $k$ and whose morphisms are linear functions.

\item $\Coalg$, the category whose objects are coalgebras
with counit
over a field $k$ and whose morphisms are coalgebra homomorphisms
and compatible with counit.

\item $\CoComCoalg$, the category whose objects are cocommutative
coalgebras
with counit
over a field $k$ and whose morphisms are cocommutative
coalgebra homomorphisms
and compatible with counit.

%\item $\LieAlg$, the category whose objects are Lie algebras over
%a field $k$ and whose morphisms are Lie algebra homomorphisms.
%
%\item $\HopfAlg$, the category whose objects are Hopf algebras
%over a field $k$ and whose morphisms are Hopf algebra
%homomorphisms.

%%{\rm We will remind the reader of the definitions of coalgebras,
%%Lie algebras, and Hopf algebras in the appropriate sections.}

\end{itemize}
\end{remark}

\begin{figure}
\begin{center}
\mbox{ \epsfxsize=4in \epsfbox{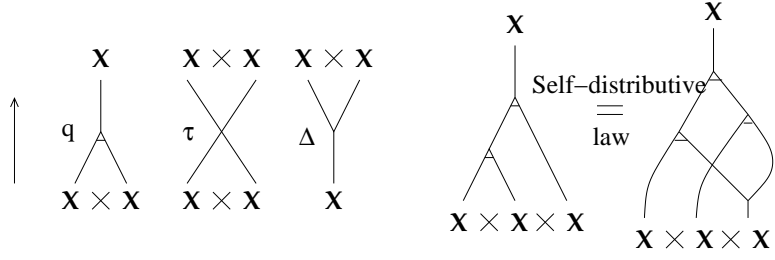} }
\end{center}
\caption{Internal Shelf Axioms}
 \label{distribute}
\end{figure}

It is convenient for calculations to express the maps and axioms
of a shelf in $K$ diagrammatically as we do in the left and right
of Fig.~\ref{distribute}, respectively. The composition of the
maps is read from right to left $(gf)(x)=g(f(x))$ in text and from
bottom to top in the diagrams. In this way, when reading from left
to right one can draw from top to bottom and when reading a
diagram from top to bottom, one can display the maps from left to
right. The argument of a function (or input object from a
category)  is found at the bottom of  the diagram.

\subsection{Coalgebras}\label{coalgsec}

A {\it coalgebra} is a vector space $C$ over a field $k$ together
with a {\it comultiplication} $\Delta: C \rightarrow C \otimes C$
that is bilinear and {\it coassociative}: $(\Delta \otimes 1)
\Delta = (1 \otimes \Delta)\Delta$.  A coalgebra is {\it
cocommutative} if the comultiplication satisfies $\tau \Delta =
\Delta$, where $\tau: C \otimes C \rightarrow  C \otimes C$ is the
transposition $\tau(x \otimes y)=y \otimes x$. A {\it coalgebra
with counit} is a coalgebra with a linear map called the {\it
counit} $\epsilon: C \rightarrow k$ such that $ (\epsilon \otimes
1)\Delta = 1 = (1 \otimes \epsilon)\Delta $ via $k \otimes C \cong
C$. Diagrammatically, this condition says that the following
commutes:
\[ \vcenter{ \xymatrix{
 k \times C \ar[dr]
& C \otimes C \ar[l]_{\epsilon \times 1} \ar[r]^{1 \times
\epsilon}
& C \times k \ar[dl] \\
& C \ar[u]_{\Delta}}}
\]

Note that if $(C, \Delta, \epsilon)$ is a coalgebra with counit,
then so is  the tensor product $C \otimes C.$
%JSC REMOVED THE REST BC IT ADDS TO CONFUSION. THE PROBLEM IS
%SOME ALGEBRAIST CONSIDER DELTA OX DELTA TO INCLUDE THE %TRANSPOSITION.
%, \Delta \otimes \Delta, \epsilon \otimes
%\epsilon)$.

\begin{lemma}
If $C$ is a coalgebra with counit, the comultiplication $\Delta_C
: C \rightarrow C \otimes C$ is the diagonal map in the category
of coalgebras with counits.
\end{lemma}
{\it Proof.} Since $C \otimes C$ is the product in the category of
coalgebras with counits, there is a
diagonal, that is a
unique morphism $\phi : C
\rightarrow C \otimes C$ which makes the following diagram
commute:
\[
  \vcenter{
  \xymatrix@!C{
   &   C
    \ar[ddl]_{1}
    \ar[dd]^{\phi}
    \ar[ddr]^{1}
    \\ \\
  C &
   C \otimes C
   \ar[l]^{\pi_{1}}
   \ar[r]_{\pi_{2}} &
   C}}
\]
where the $\pi_1$ and $\pi_2$ are projection maps
defined by
$$\pi_1 : = \xymatrix{
    A \otimes B \ar[rr]^{1 \otimes \epsilon_B}
    &&
    A \otimes k \ar[rr]^{\sim}
    &&
    A}$$
$$\pi_2 : = \xymatrix{
    A \otimes B \ar[rr]^{\epsilon_A \otimes 1}
    &&
    k \otimes B \ar[rr]^{\sim}
    &&
    B}$$
where $\epsilon_A$ and $\epsilon_B$ are the counit maps for
coalgebras $A$ and $B$.
Since the comultiplication $\Delta_C$ satisfies the same property as
$\phi$ and $\phi$ is unique,
they must coincide.
$\Box$

A linear map $f$ between coalgebras
is said to be
{\it compatible with comultiplication},
or {\it preserves comultiplication},  if it satisfies the condition
$\Delta f = (f \otimes f) \Delta$.
Diagrammatically, the following commutes:
\[
  \xymatrix{
   C
     \ar[rr]^{\Delta_C}
     \ar[dd]_{f}
     && C \otimes C
     \ar[dd]^{f \otimes f} \\ \\
   D
     \ar[rr]^{\Delta_D}
     && D \otimes D }
\]

A linear map $f$ between coalgebras
is said to be {\it compatible with counit}, or {\it preserves counit},  if it
satisfies the condition $\epsilon f=\epsilon$, which,
diagrammatically says the following diagram commutes:
\[
\xymatrix{
     C
     \ar[r]^{\epsilon_C}
      \ar[dr]_{f}
      & k \\
     & D \ar[u]_{\epsilon_D}}
\]

In particular, if $(C, \Delta, \epsilon)$ is
a coalgebra with counit,  a linear map $q: C \otimes C \rightarrow C $ between
coalgebras is compatible with comultiplication
if and only if it satisfies
$\Delta q= (q \otimes  q)(1  \otimes \tau  \otimes  1)(\Delta  \otimes  \Delta)$,
and it is compatible with counit if and only if
it satisfies
$\epsilon q= q (\epsilon \otimes \epsilon)$.

A {\it morphism} $f$ in  the category of coalgebras with counit
is a linear map that preserves comultiplication and counit.
As suggested by the categories listed in Remark ~\ref{categories},
we will focus our main attention on coalgebras with counits. Thus,
we use the word `coalgebra' to refer to a coalgebra with counit
and the phrase `coalgebra morphism' to refer to a linear map that
preserves comultiplication and counit.
On the other hand, we wish to consider examples in which
the self-distributive map is not compatible with the counit (see the sequel).
For categorical hygiene, we are distinguishing a function that satisfies
self-distributivity and is compatible with comultiplication
from {\it a morphism in the category $\Coalg.$ }

\section{Self-Distributive Maps for Coalgebras }
\label{coalsec}

In this section we give concrete and broad examples of
self-distributive maps
for
cocommutative coalgebras. Specifically, we
discuss examples constructed from quandles/racks used as bases,
Lie algebras, and Hopf algebras.

\subsection{Self-Distributive Maps
for Coalgebras
Constructed From Racks}
 \label{rackshelfsec}

In this section we note that quandles and racks  can be used to
construct self-distributive maps in $\CoComCoalg$
  simply by using their elements as basis.

Let $X$ be a rack.
Let $V=kX$ be the vector space over a
field $k$ with the elements of $X$ as basis. Then $V$ is a
cocommutative coalgebra with counit, with
comultiplication $\Delta$ induced by the diagonal map $\Delta(x )=
x \otimes x$, and the counit induced by $\epsilon(x)=1$ for
$x \in X$. This is a standard construction of a coalgebra
with counit from a set.

Set $W=k \oplus kX$.
We denote an element of $W=k \oplus kX$ by
$a + \sum_{x \in X} a_x x$ or more briefly by
$a + \sum_x a_x x$,
and when context is understood by
  $a + \sum a_x x$.
 Extend $\Delta $ and $\epsilon$ on $V=kX$ to $W$ by
linearly extending $\Delta (1)=1 \otimes 1$ and $\epsilon(1)=1$
for $1 \in k$. More explicitly,
$$\Delta (a + \sum a_x x ) = a (1 \otimes 1) + \sum a_x (x \otimes x), $$
and $\epsilon(a + \sum a_x x ) =a+ \sum a_x$. With these
definitions, one can check
 that $(W, \Delta, \epsilon)$ is an object in
$\CoComCoalg$.

Define $q: W \otimes W \rightarrow W$ by
 linearly extending $q(x \otimes y)=x\lt y$,
 $q (1 \otimes x)=1$, $q(x \otimes 1)=0$,
  and   $q(1 \otimes 1)=0$. More explicitly,
$$q( \  (a+ \sum a_x x) \otimes (b+ \sum b_y y )\ )
=  \sum_{y } a b_y + \sum_{x, y } a_x b_y (x \lt y) . $$

\begin{proposition}
\label{quandleQprop} The extended map $q$ given above
is a self-distributive linear map compatible with comultiplication.
\end{proposition}
{\it Proof.\/} We begin by checking that $q$ satisfies
self-distributivity and continue by showing that $q$ is compatible
with comultiplication.
 In the second case,
we check that $\Delta q= (q\otimes q)(1 \otimes \tau \otimes
1)(\Delta \otimes \Delta)$.

Then one computes:
\begin{eqnarray*}
\lefteqn{q(q\otimes 1)
 ( \  (a+ \sum a_x x) \otimes (b+ \sum b_y y ) \otimes (c + \sum c_z z)\ )}\\
 &=& q(\  ( \sum_{y } a b_y + \sum_{x, y } a_x b_y (x \lt y) )
        \otimes (c+\sum c_z z) \ ) \\
 &=& \sum_{y,z } a b_y c_z + \sum_{x, y, z} a_x b_y c_z ((x \lt y) \lt z ) ,\\
 \lefteqn{q(q \otimes q) (1 \otimes \tau \otimes 1) (1 \otimes \Delta)
  ( \ (a+ \sum a_x x) \otimes (b+ \sum b_y y ) \otimes (c + \sum c_z z)   \ ) }\\
&=& q(q \otimes q) (1 \otimes \tau \otimes 1) (\  (a+ \sum a_x x)
\otimes (b+ \sum b_y y ) \otimes
(c (1 \otimes 1) + \sum c_z (z \otimes z) ) \ ) \\
&=& q( q \otimes q)
( \ c ((a+ \sum a_x x) \otimes 1 \otimes  (b+ \sum b_y y ) \otimes 1)\\
& & \quad + \sum_z c_z  (a+ \sum a_x x) \otimes z \otimes  (b+ \sum b_y y ) \otimes z \ )\\
&=& q( 0 + \sum_z c_z  (a + \sum a_x (x \lt z) ) \otimes  (b+ \sum b_y (y \lt z)  \ )\\
&=& \sum_z c_z (\sum_y a b_y + \sum_{x, y} a_xb_y (x \lt z)\lt (y \lt z) ),
 \end{eqnarray*}
 as desired. Compatibility with comultiplication is checked as follows:
 \begin{eqnarray*}
\lefteqn{ \Delta q ( \  ( a + \sum a_x x ) \otimes (b + \sum b_y y ) \ ) }\\
&=& \Delta  (\sum_y a b_y +  \sum_{x,y}  a_x b_y (x \lt y)  )
  \  = \  \sum_y  a b_y (1 \otimes 1) + \sum_{x,y}  a_x b_y (x \lt y)  \otimes (x \lt y) ,\\
\lefteqn{  (q \otimes q ) (1 \otimes \tau \otimes 1) (\Delta
\otimes \Delta )
( \  ( a + \sum a_x x ) \otimes (b + \sum b_y y ) \ ) } \\
&=&
  (q \otimes q) (1 \otimes \tau \otimes 1) (\ ( a (1 \otimes 1) + \sum a_x ( x \otimes x) ) \otimes
  ( b (1 \otimes 1) + \sum b_y (y \otimes y) \ ) \\
  &=&
  (q \otimes q)  (\ ( ab (1 \otimes 1\otimes 1\otimes 1)
  + \sum b a_x (  x\otimes 1\otimes x \otimes 1)  \\
  & & \quad
  + \sum a  b_y (1\otimes y\otimes 1 \otimes y )
 + \sum_{x,y}  a_x b_y ( x \otimes y \otimes x \otimes y) \ ) \\
 &=&
 \sum_y a b_y (1 \otimes 1) + \sum_{x,y}  a_x b_y (x \lt y)  \otimes (x \lt y) . \quad \Box
 \end{eqnarray*}

The pair
$(W, q)$ falls short of being a shelf in $\CoComCoalg$
due to the following:

 \begin{proposition}
The extended map $q$ defined above
is not compatible with the counit, but satisfies $\epsilon q = q
(\epsilon \otimes 1)$.
\end{proposition}
{\it Proof.\/}
The counit $\epsilon$ has as its image $k \subset W$. Thus the image of $\epsilon \otimes 1$ is $W\otimes W$. We compute
the following three quantities:
 \begin{eqnarray*}
 \epsilon  q ( \ ( a + \sum a_x x ) \otimes (b + \sum b_y y ) \ )
 &=& \epsilon  ( \sum a b_y + \sum_{x,y} a_x b_y (x \lt y) )\\
 &=&  a \sum b_y +  \sum_{x,y} a_x b_y ,\\
\epsilon  \otimes \epsilon
( \ ( a + \sum a_x x ) \otimes (b + \sum b_y y ) \ )
&=& (a + \sum a_x)  (b + \sum b_y )  \\
&=&  ab + a \sum b_y + b \sum a_x + \sum a_x b_y , \ {\mbox{\rm and}} \\
 q (\epsilon  \otimes 1   )
( \ ( a + \sum a_x x ) \otimes (b + \sum b_y y ) \ )
&=& q(\  (a + \sum a_x) \otimes (b + \sum b_y y  )\  ) \\
&=& (a + \sum a_x) \sum b_y .
\end{eqnarray*}
The first and third coincide.
 $\Box$

\subsection{Lie Algebras}\label{Liesec}

A {\it Lie algebra} $\g$ is a vector space over a field $k$ of characteristic other than 2, with an
antisymmetric
bilinear form $[ \cdot , \cdot ]: \g \otimes \g \rightarrow \g$
that satisfies the Jacobi identity $[[x,y],
z]+[[y,z],x]+[[z,x],y]=0$ for any $x,y,z \in \g$. Given a Lie
algebra $\g$ over $k$ we can construct a coalgebra $N = k \oplus
\g$.  We will denote elements of $N$ as either $(a,x)$ or $a+x$,
depending on clarity,
where $a \in k$ and $x \in \g$.

In fact, $N$ is a
cocommutative
coalgebra with
comultiplication and counit
given by $\Delta(x)=x \otimes 1 + 1 \otimes x$ for $x \in \g$ and
$\Delta(1)=1 \otimes 1$, $\epsilon(1)=1$, $\epsilon(x)=0$ for $x
\in \g$. In general we compute,
for $a \in k$ and $x \in \g$,
\begin{eqnarray*}
\lefteqn{
\Delta((a,x)) =
\ \Delta(a + x)  =  \Delta(a)+\Delta(x) } \\
&=& a (1 \otimes 1) + x \otimes 1 + 1 \otimes x
\ = \ (a \otimes 1 + x \otimes 1 ) + 1 \otimes x \\
&=& (a+x) \otimes 1 + 1 \otimes x
 \
=    (a, x) \otimes (1, 0) + (1,0) \otimes (0, x).
\end{eqnarray*}

The following
map
 is found in quantum group theory
(see for example, \cite{MajidGreen}, and studied   in \cite{Alissa}
in relation to Lie $2$-algebras).
Define $q: N \otimes N \rightarrow N$ by
linearly extending $q( 1 \otimes (b+y) )=\epsilon(b+y)$, $q((a+x) \otimes 1)=a+x$
and $q(x, y )=[x,y]$ for
$a, b \in k$ and
$x, y \in \g$,
i.e.,
$$
q((a,x) \otimes (b,y))=
q((a+x ) \otimes (b+y) )= q((a + x ) \otimes (b+y) ) = ab + bx + [x,y]
= (ab, bx + [x,y]).$$

Since the solution to the classical YBE follows from the Jacobi
identity, and the YBE is related to self-distributivity (see next
section) via the third Reidemeister move, it makes sense to expect
that there is a relation between the Lie bracket and the
self-distributivity axiom.

\begin{lemma}
The above defined $q$ satisfies the self-distributive law in
Definition~\ref{selfdist}.
 \end{lemma}
{\it Proof.\/} We compute
\begin{eqnarray*}
\lefteqn{q (q \otimes 1 ) ((a, x) \otimes (b,y) \otimes (c, z) )  } \\
&=& q ( (ab + bx +  [x,y]) \otimes (c,z) ) \ = \ abc + bcx + c[x,y] + b[x,z] +  [[x,y],z], \\
\lefteqn{ q (q \otimes q) (1 \otimes \tau \otimes 1 ) (1 \otimes 1
\otimes \Delta)
(  (a+x) \otimes (b+y) \otimes (c+z) ) } \\
&=& q (q \otimes q)  (1 \otimes \tau \otimes 1 ) ( (a+x) \otimes
(b+y) \otimes
\{ (c+z) \otimes 1  +  1 \otimes  z \}  ) \\
&=& q (q \otimes q) ( (a+x) \otimes (c+z)  \otimes (b+y) \otimes 1
+
(a+x) \otimes 1  \otimes (b+y) \otimes z \\
&=& q(  ( ac + cx + [x,z]  ) \otimes  (b+y ) ) +
q( (a+x) \otimes [y,z] ) \\
&=& (abc + bcx + c[x,y] + b[x,z] + [[x,z],y] ) + [x,[y,z]],
\end{eqnarray*}
and the Jacobi identity in $\g$ verifies the condition. $\Box$

\begin{lemma} The map
 $q$ constructed above is a coalgebra
morphism.
\end{lemma}
{\it Proof.\/}
We compute:
$$\Delta q
((a+x) \otimes (b+y))=(ab+b x + [x,y])\otimes 1 + 1 \otimes (bx + [x,y]).$$
On the other hand,
we have
\begin{eqnarray*}
\lefteqn{
(q \otimes q)(1 \otimes \tau \otimes 1 )(\Delta \otimes \Delta) ((a+x) \otimes (b+y)) } \\
& = & (q \otimes q)(1 \otimes \tau \otimes 1 )
((a+x)\otimes 1 + 1 \otimes x)
\otimes
((b+y) \otimes 1 +  1 \otimes y )\\
&=  & q((a+x)\otimes(b+y))\otimes q(1 \otimes 1)
+ q(1 \otimes 1) \otimes q(x \otimes y)+ \epsilon (b+y) \otimes x
 + (a+x)\otimes \epsilon (y) \\
& = & ((a+ b x + [x,y]) \otimes 1)
+ 1 \otimes [x,y]
+ b
\otimes x
+ (a+x)\otimes 0
\\
&=& ((a+ b x + [x,y]) \otimes 1)  +1 \otimes (b x+[x,y])
\end{eqnarray*}
For the counit, we compute:
$$ \epsilon q (\ (a+x)\otimes (b+y) \ ) = \epsilon (ab + bx + [x, y])=ab =
(\epsilon \otimes \epsilon) (\ (a+x)\otimes (b+y) \ ). \quad
\Box$$
Combining these two lemmas, we have:
\begin{proposition}\label{lieQprop}\begin{sloppypar}
The coalgebra $N$ together with map $q$ given above defines a
shelf $(N, q)$ in $\CoComCoalg.$ \end{sloppypar}
%%%%We need cocom later.
\end{proposition}

Groups have quandle structures given by conjugation, and their
subset Lie groups are related to Lie algebras through tangent
spaces and exponential maps. In the above proposition we
constructed shelves in $\CoComCoalg$ from Lie algebras, so we see
this proposition as a step in completing the following square of
relations.
\[
  \xymatrix{
   {\rm Lie\  groups}
     \ar[rr]^{}
     \ar[dd]_{}
     && {\rm Lie\ algebras}
     \ar[dd]^{} \\ \\
   {\rm Quandles}
     \ar[rr]^{}
     && ??? }
\]

\subsection{Hopf Algebras}
\label{Hopfsec}

A  {\it bialgebra} is an algebra $A$ over a field $k$ together
with a linear map called the {\it unit} $\eta: k \rightarrow A$,
satisfying $\eta(a)=a {\bf 1}$ where ${\bf 1} \in A$ is the
multiplicative identity and with an associative multiplication
$\mu: A \otimes A \rightarrow A$ that is also a coalgebra such
that the comultiplication $\Delta$ is an algebra homomorphism. A
{\it Hopf algebra} is a bialgebra $C$ together with a map called
the {\it antipode} $S: C \rightarrow C$ such that $\mu (S \otimes
1) \Delta = \eta \epsilon = \mu (1 \otimes S) \Delta$, where
$\epsilon$ is the counit.

\begin{figure}[htb]
\begin{center}
\mbox{
\epsfxsize=3.5in
\epsfbox{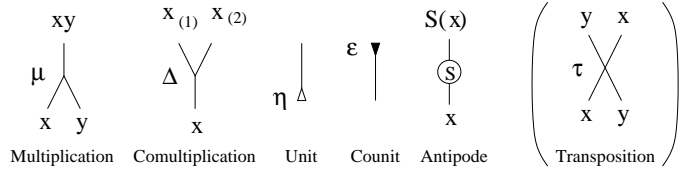}
}
\end{center}
\caption{Operations in Hopf algebras }
\label{Hopfdiag}
\end{figure}

The reader can construct commutative diagrams similar to those
found in Section \ref{coalgsec} for the notions of bialgebra and
Hopf algebra. Our diagrammatic conventions for these maps are
depicted in Fig.~\ref{Hopfdiag}.
Recall that the diagrams are read from bottom to top.
These diagrams have been used (see for example
\cite{Kup1,YetBook}) for proving facts about Hopf algebras
and related invariants.

We review the diagrammatic representation of Hopf algebra axioms. For convenience, assume that the underlying vector space of $A$ is finite dimensional with ordered basis $(e_1,e_2, \ldots, e_n)$. Then the multiplication $\mu$
and comultiplication
$\Delta$
 are determined by the values, $\Lambda_{ij}^\ell, Y^{ij}_\ell \in k$, of the structure constants:
$\mu(e_i \otimes e_j) = \Lambda_{ij}^\ell(e_\ell)$, and $\Delta(e_\ell)=Y^{ij}_\ell e_i \otimes e_j.$ Note that summation conventions are being applied, and so, for example, $ \Lambda_{ij}^\ell(e_\ell)=\sum_{\ell=1}^n  \Lambda_{ij}^\ell(e_\ell)$. Similarly,  the unit can be written as $\eta(1) =\sum_i A^i e_i$. The co-unit can be written as $\epsilon (e_i)= V_i \in k$, so that for a general vector, $\sum_i \alpha^i e_i$, we have $\epsilon(\sum \alpha^i e_i)= (\sum_i A^i) \epsilon(e_i)= \sum_i a^i V_i.$ Finally, the antipode is a linear map so $S(e_i)=s^j_i e_j$ for constants $s^j_i \in k$.

Thus the axioms of a (finite dimensional) Hopf algebra can be formulated in terms of the structure constants. The table below summarizes these formulations. Again summation convention applies, and all super, and subscripted variables are constants in the ground field.
\begin{center}\begin{large}
\begin{tabular}{||l|c||} \hline \hline
associativity &  $\Lambda_{ir}^q\Lambda^r_{j \ell} = \Lambda^q_{p \ell}\Lambda^p_{ij}$ \\  \hline
   coassociativity &
$Y^{ij}_p Y_q^{p \ell} = Y_r^{j \ell}Y_q^{ir}$ \\ \hline
unit  & $\Lambda^\ell_{ij} A^j= A^j \Lambda^\ell_{ji} =\delta^\ell_i$  \\  \hline
co-unit &
$V^i Y^{ij}_\ell = Y^{j i}_\ell V_i = \delta^j_\ell$ \\ \hline
Compatibility & $\Lambda^p_{tv} \Lambda^q_{uw}Y^{tu}_iY^{vw}_j=Y^{pq}_r\Lambda^r_{ij}$ \\ \hline
Antipode & $\Lambda^{i}_{rq}s^{r}_p Y^{pq}_j=\Lambda^{i}_{rq}s^q_pY^{rp}_j=V_iA^j$ \\ \hline \hline
\end{tabular}\end{large}
\end{center}

In the table above, $\delta^{\ell}_i$ denotes a Kronecker delta function. It is a small step, now to translate these
Specifically, the multiplication tensor $\Lambda$ is diagrammatically
represented by the leftmost trivalent vertex read from bottom to top. %%%,
The letter choices $\Lambda,$ $Y$, $A$ and $V$ are meant to suggest the graphical depictions of these operators.
A composition of maps
corresponds to a contraction of the same indices of tensors which,
in turn,
corresponds to
connecting end points of diagrams together
 vertically.
 Figures~\ref{Hopfdiag} and \ref{HopfAxioms} represents such diagrammatic conventions of
 maps that appear in the definition of a Hopf algebra and their axioms.
 %%%%%
%%%%%end add from JSC addend

\begin{figure}[htb]
\begin{center}
\mbox{
\epsfxsize=5in
\epsfbox{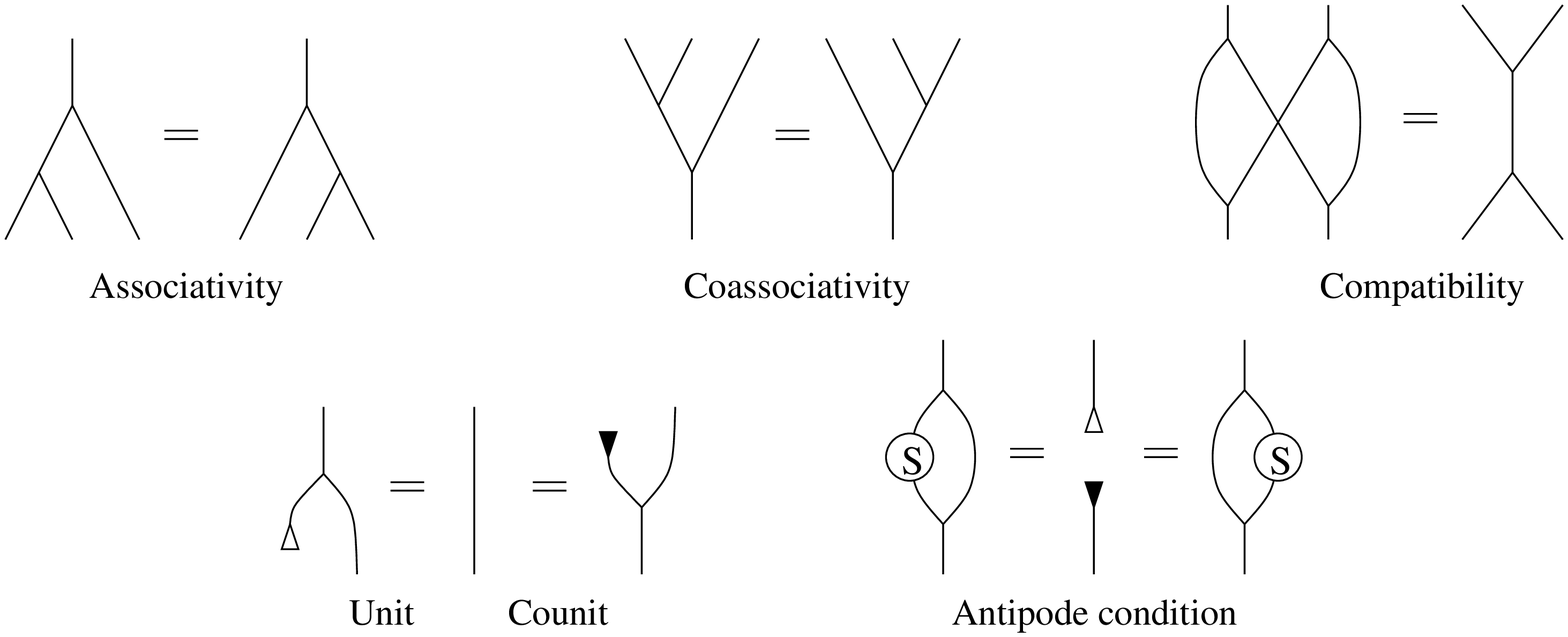}
}
\end{center}
\caption{Axioms of Hopf algebras }
\label{Hopfaxioms}
\end{figure}

\begin{figure}[htb]
\begin{center}
\mbox{
\epsfxsize=1in
\epsfbox{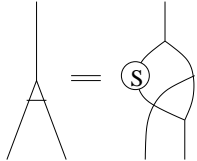}
}
\end{center}
\caption{Self-distributive map in Hopf algebras}
 \label{HopfQ}
\end{figure}

Let $H$ be a Hopf algebra. Define $q: H \otimes H \rightarrow H$
by $q=\mu(1  \otimes \mu)(S \otimes 1 \otimes 1 ) (\tau \otimes
1)(1 \otimes \Delta)$ where $\mu$, $\Delta,$ and  $S$ denote the
multiplication, comultiplication, and  antipode, respectively. If
we adopt the common notation $\Delta(x)=x_{(1)} \otimes x_{(2)}$
and $\mu(x \otimes y)=xy$, then $q$ is written as $q(x \otimes y)=
S(y_{(1)}) x y_{(2)}$.
This
appears as an {\it adjoint map}
  in \cite{Woro,MajidPink},
  and
its diagram is depicted in Fig.~\ref{HopfQ}. Notice the analogy
with the group conjugation as a quandle: in a group ring,
$\Delta(y) = y \otimes y$ and $S(y)=y^{-1}$, so that $q(x \otimes
y)=y^{-1}xy$, and therefore, is of a great interest from point of
view of quandles.

\begin{figure}[htb]
\begin{center}
\mbox{
\epsfxsize=2.5in
\epsfbox{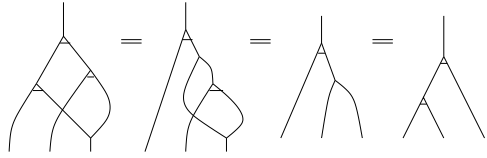}
}
\end{center}
\caption{Proof of self-distributivity in Hopf algebras}
\label{coalgQlem}
\end{figure}

\begin{proposition}\label{hopfQprop}
The above defined
linear
map $q: H \otimes H \rightarrow H$
satisfies the self-distributive law in Definition~\ref{selfdist}.
\end{proposition}
{\it Proof.\/} In Fig.~\ref{coalgQlem}, it is indicated that this
follows from two properties of the adjoint map:
  $q(q  \otimes 1)=q( 1 \otimes \mu)$
  (which is used in the first and the third equalities in the figure),
 and
 $\mu=\mu(1 \otimes q)(\tau \otimes 1)(1 \otimes \Delta)$
 (which is used in the second equality).

\begin{figure}[htb]
\begin{center}
\mbox{
\epsfxsize=4in % 5in
\epsfbox{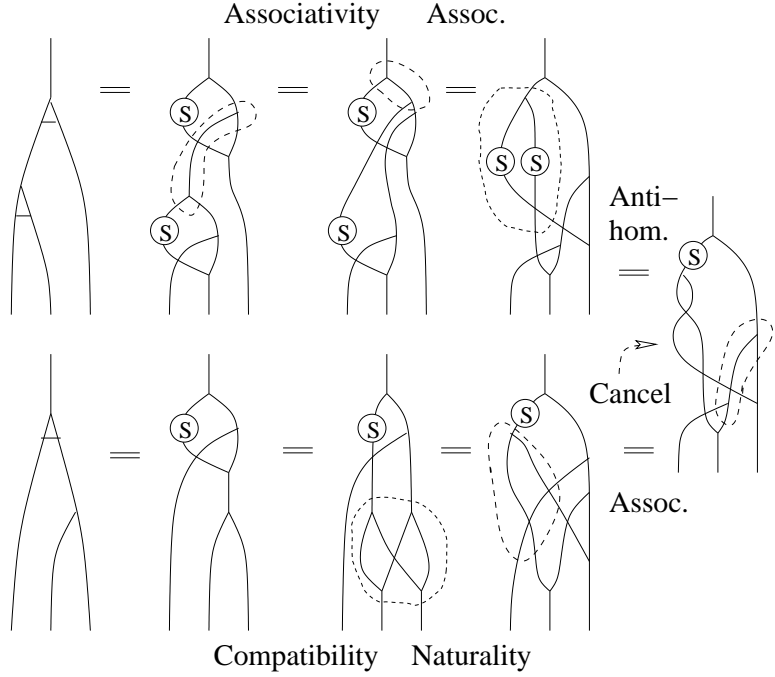}
}
\end{center}
\caption{ $q(q  \otimes 1)=q( 1 \otimes \mu)$ } \label{HopfQlemma}
\end{figure}

\begin{figure}[htb]
\begin{center}
\mbox{
\epsfxsize=2.5in
\epsfbox{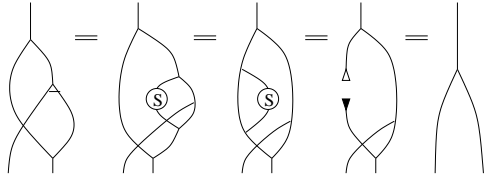}
}
\end{center}
\caption{ $\mu=\mu(1 \otimes q)(\tau \otimes 1)(1 \otimes
\Delta)$} \label{mRlem}
\end{figure}

It is known that these properties are satisfied, and proofs are
found in \cite{Woro,Henn}. Here we include diagrammatic proofs for
reader's convenience in Fig.~\ref{HopfQlemma} and
Fig.~\ref{mRlem}, respectively.
$\Box$

\begin{remark}
\label{RmkH}{\rm \begin{sloppypar}
The definition of $q$ above contains an antipode, which is a
coalgebra {\it anti-homomorphism} and not necessarily a coalgebra
morphism.  Thus, $(H,q)$ is not a shelf in $\Coalg$ in general.
\end{sloppypar}
}
\end{remark}

\subsection{Other Examples}

In this section we observe that there are plenty of examples of
self-distributive linear maps for $2$-dimensional cocommutative
coalgebras and shelves in $\CoComCoalg$.

Let $V$ be the two dimensional vector space over $k$ with basis
$\{ x, y \}$. Define a coalgebra structure on $V$ using the
diagonal map $\Delta(z)=z \otimes z$ for $z \in \{ x, y \}$ and
extending it linearly.

\begin{lemma}
A linear map $q: V \otimes V \rightarrow V$ is self-distributive
and compatible with comultiplication if and only if $q$ is a map
defined on basis elements in the following list:

$$
\begin{array}{ll|l|l|l|l|l|l|l|l|l|l|l|l|l|l|l|l|l|l|l|l}
q(x\otimes x)=&0&0&0&0&0&0&x&x&x&x&x&x&x&x&x&y&y&y&y&y&0\\
q(x\otimes y)=&0&0&0&0&x&y&0&0&0&x&x&x&x&y&y&0&x&y&y&y&0\\
q(y\otimes x)=&0&0&x&x&0&0&0&0&y&x&x&y&y&x&y&0&y&0&x&y&0\\
q(y\otimes y)=&x&y&0&x&y&0&0&y&0&x&y&x&y&y&y&0&y&0&x&y&0
\end{array}
$$
Among these, $(V, q)$ is a shelf in $\CoComCoalg$ if and only if
$q(a,b)\ne0$ for any $a,b \in \{x,y\}.$

\end{lemma}
{\it Proof.\/} Let $q(x \otimes x)=\gamma_1 x + \gamma_2 y$ for
some constants $\gamma_1$, $\gamma_2 \in k$. The compatibility
condition
$$\Delta q (x \otimes x)=(q \otimes q)(1 \otimes \tau \otimes 1)(\Delta \otimes \Delta)(x
\otimes x)$$ implies that $\gamma_1 \gamma_2=0$ and
$\gamma_1^2=\gamma_1$, $\gamma_2^2=\gamma_2$, i.e., $q(x \otimes
x)=0$, $x$, or $y$. The same holds for $x\otimes y$, $y \otimes x$
and $y \otimes y$, so that the value of $q$ for a pair of basis
elements is either a basis element ($x$ or $y$), or $0$.

A case by case analysis (facilitated by {\sc Mathematica} and/or
{\sc Maple}) provides self-distributivity. When $\epsilon{(x)}=\epsilon{(y)}=1$, the only cases for which $\epsilon q= \epsilon \otimes \epsilon$ are those for which $q(a,b)\ne0$ for all four choices of $a,b$.
$\Box$

\bigskip

Another famous example of a cocommutative coalgebra is the
trigonometric coalgebra, $T$,
generated by $a$ and $b$ with
comultiplication given by:
\begin{eqnarray*}
\Delta(a) & =& a \otimes a - b \otimes b \\
\Delta(b) &=& a \otimes b + b \otimes a
\end{eqnarray*}
with counit $\epsilon(a)=1$, $\epsilon(b)=0$, in analogy with
formulas for $\cos(x+y)$ and $\sin(x+y)$ and $\cos(0)=1$,
$\sin(0)=0$.

\begin{lemma} Let $T$ denote the trigonometric coalgebra
over $\C$.
Let
$q:  T \otimes T \rightarrow T$ be a linear map defined by:
\begin{eqnarray*}
q(a \otimes a) \  =\  \alpha_1 a + \beta_1 b, & &
q(a \otimes b) \  =\ \alpha_2 a + \beta_2 b, \\
q(b \otimes a) \  =\  \alpha_3 a + \beta_3 b,  & &
q(b \otimes b) \  =\  \alpha_4 a + \beta_4 b.
\end{eqnarray*}
Then such a linear map  $q$ is self-distributive and
compatible with comultiplication
if and only if the coefficients are found in Table 1, where $i=\sqrt{-1}$.

Among these, $(V, q)$ is a shelf in $\CoComCoalg$ if and only if
$(\alpha_1 ,\alpha_2, \alpha_3, \alpha_4)=(1,0,0,0)$.

\end{lemma}
{\it Proof.\/} This result is a matter of verifying the conditions for self-distributivity and compatibility over all possible choices of inputs. We generated solutions by both {\sc Maple} and {\sc Mathematica}. For the compatibility condition we established a system of $12$ quadratic equations  in eight unknowns. Originally there were $16$ such equations, but $4$ of these are duplicates. In the {\sc Mathematica}  program we used the command ``Solve'' to generate a set of necessary conditions. The self-distributive condition gave a system of cubic equations in the unknowns. We checked these subject to the necessary conditions, and found the $21$ solutions above.

Expressing $\epsilon$ as a $(1 \times 2)$ matrix and $q$ as the $2 \times 4$ matrix
$\left( \begin{array}{cccc} \alpha_1 &\alpha_2& \alpha_3 & \alpha_4 \\ \beta_1 &\beta_2& \beta_3&  \beta_4 \end{array} \right).$
We compute $\epsilon q = (\alpha_1 ,\alpha_2, \alpha_3, \alpha_4)$ and $ \epsilon \otimes \epsilon =(1,0,0,0)$. The result follows.
  $\Box$

\begin{table}[htb]
$$
\begin{array}{||cccc||cccc||}\hline \hline
\alpha_1 & \alpha_2 & \alpha_3 & \alpha_4 & \beta_1 & \beta_2 & \beta_3 & \beta_4 \\ \hline
\hline
 1 & 0 & 0 & 0 &
 0 & 0 & - 1 & 0 \\ \hline
 0 & 0 & 0 & 0 &
 0 & 0 & 0 & 0 \\ \hline
 \frac{1}{2} & - \frac{i}{2} & 0 & 0 &
 0 & 0 & \frac{1}{2} & - \frac{i}{2} \\ \hline
 \frac{1}{2} & \frac{i}{2} & 0 & 0 &
 0 & 0 & \frac{1}{2} & \frac{i}{2} \\ \hline
 1 & 0 & 0 & 0 &
 0 & 0 & 1 & 0 \\ \hline
 \frac{1}{2} & 0 & 0 & - \frac{1}{2} &
 0 & \frac{1}{2} & \frac{1}{2} & 0 \\ \hline
 1 & 0 & 0 & 0 &
 0 & 1 & 0 & 0 \\ \hline
 \frac{1}{4} & - \frac{i}{4} & - \frac{i}{4} & - \frac{1}{4} &
 - \frac{i}{4} & - \frac{1}{4} & - \frac{1}{4} & \frac{i}{4} \\ \hline
 \frac{1}{4} & \frac{i}{4} & - \frac{i}{4} & \frac{1}{4} &
 - \frac{i}{4} & \frac{1}{4} & - \frac{1}{4} & - \frac{i}{4} \\ \hline
 \frac{1}{4} & \frac{i}{4} & \frac{i}{4} & - \frac{1}{4} &
 - \frac{i}{4} & \frac{1}{4} & \frac{1}{4} & \frac{i}{4} \\ \hline
 \frac{1}{4} & \frac{i}{4} & \frac{i}{4} & - \frac{1}{4} &
 \frac{i}{4} & - \frac{1}{4} & - \frac{1}{4} & - \frac{i}{4} \\ \hline
 \frac{1}{4} & - \frac{i}{4} & \frac{i}{4} & \frac{1}{4} &
 \frac{i}{4} & \frac{1}{4} & - \frac{1}{4} & \frac{i}{4} \\ \hline
 \frac{1}{4} & - \frac{i}{4} & - \frac{i}{4} & - \frac{1}{4} &
 \frac{i}{4} & \frac{1}{4} & \frac{1}{4} & - \frac{i}{4} \\ \hline
 1 & 0 & 0 & 0 &
 - \frac{i}{2} & - \frac{1}{2} & \frac{1}{2} & \frac{i}{2} \\ \hline
 \frac{1}{2} & 0 & - \frac{i}{2} & 0 &
 - \frac{i}{2} & 0 & - \frac{1}{2} & 0 \\ \hline
 1 & 0 & 0 & 0 &
 - \frac{i}{2} & \frac{1}{2} & \frac{1}{2} & - \frac{i}{2} \\ \hline
 1 & 0 & 0 & 0 &
 \frac{i}{2} & - \frac{1}{2} & \frac{1}{2} & - \frac{i}{2} \\ \hline
 \frac{1}{2} & 0 & \frac{i}{2} & 0 &
 \frac{i}{2} & 0 & - \frac{1}{2} & 0 \\ \hline
 1 & 0 & 0 & 0 &
 \frac{i}{2} & \frac{1}{2} & \frac{1}{2} & \frac{i}{2} \\ \hline
 1 & 0 & 0 & 0 &
 - i & 0 & 0 & 0 \\ \hline
 1 & 0 & 0 & 0 &
 i & 0 & 0 & 0 \\ \hline \hline
\end{array}$$
\caption{List of self-distributive maps in the trigonometric coalgebra}
\end{table}

\section{Yang--Baxter Equation and
Self-Distributive Maps for Coalgebras }\label{ybesec}

In this section, we discuss relationships between solutions to the Yang-Baxter equations and self-distributive maps. 
\subsection{A Brief Review of YBE}

The Yang--Baxter equation makes sense in any monoidal category.
Originally mathematical physicists concentrated on solutions in
the category of vector spaces
with
the tensor product, obtaining solutions from quantum groups.

Let $V$ be a vector space and $R: V \otimes V \rightarrow V
\otimes V$ an invertible linear map. We say $R$ is a {\it
Yang--Baxter operator} if it satisfies the {\it Yang--Baxter
equation}, (YBE), which says that: $(R \otimes 1)( 1 \otimes R) (R
\otimes 1 ) = ( 1 \otimes R) (R \otimes 1 ) ( 1 \otimes R)$. In
other words, the YBE says that the following diagram commutes:
\[
\begin{xy}
    (-2, 20)*+{V \otimes V \otimes V}="1";
    (-30,10)*+{V \otimes V \otimes V}="2";
    (-30,-10)*+{V \otimes V \otimes V}="3";
    (-2,-20)*+{V \otimes V \otimes V}="4";
    (30,-10)*+{V \otimes V \otimes V}="5";
    (30, 10)*+{V \otimes V \otimes V}="6";
        {\ar^{R \otimes 1} "1";"6"};
        {\ar_{1 \otimes R} "1";"2"};
        {\ar_{R \otimes 1} "2";"3"};
        {\ar_{1 \otimes R} "3";"4"};
        {\ar^{R \otimes 1} "5";"4"};
        {\ar^{1 \otimes R} "6";"5"};
\end{xy}
\]
A solution to the YBE is also called a {\it braiding}.

In general, a braiding operation provides a diagrammatic
description of the process of switching the order of two things.
This idea is formalized in the concept of a braided monoidal
category, where the braiding is an isomorphism
$$R_{X,Y}: X \otimes Y \rightarrow Y \otimes X.$$

If we draw $R : V \otimes V \rightarrow V \otimes V $
by the diagram:

\begin{center}
\mbox{
\epsfxsize=1.3in
\epsfbox{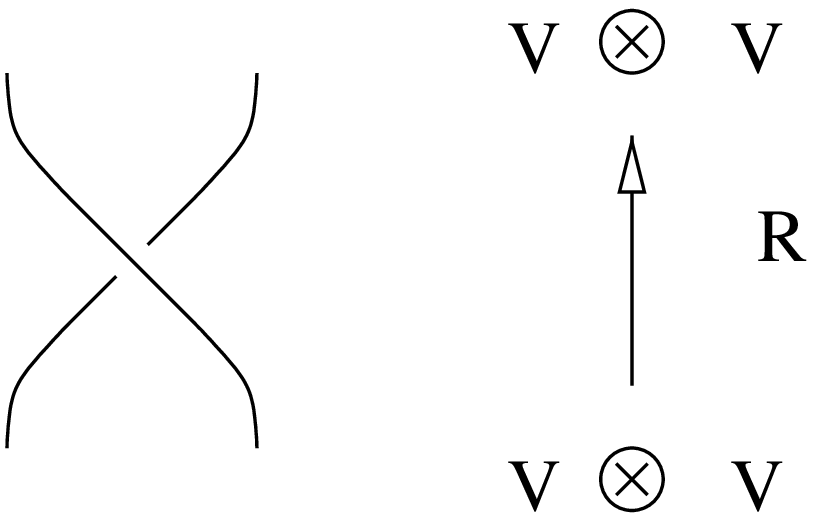}
}
\end{center}

then the Yang--Baxter equation is represented
by:

\begin{center}
\mbox{
\epsfxsize=2in
\epsfbox{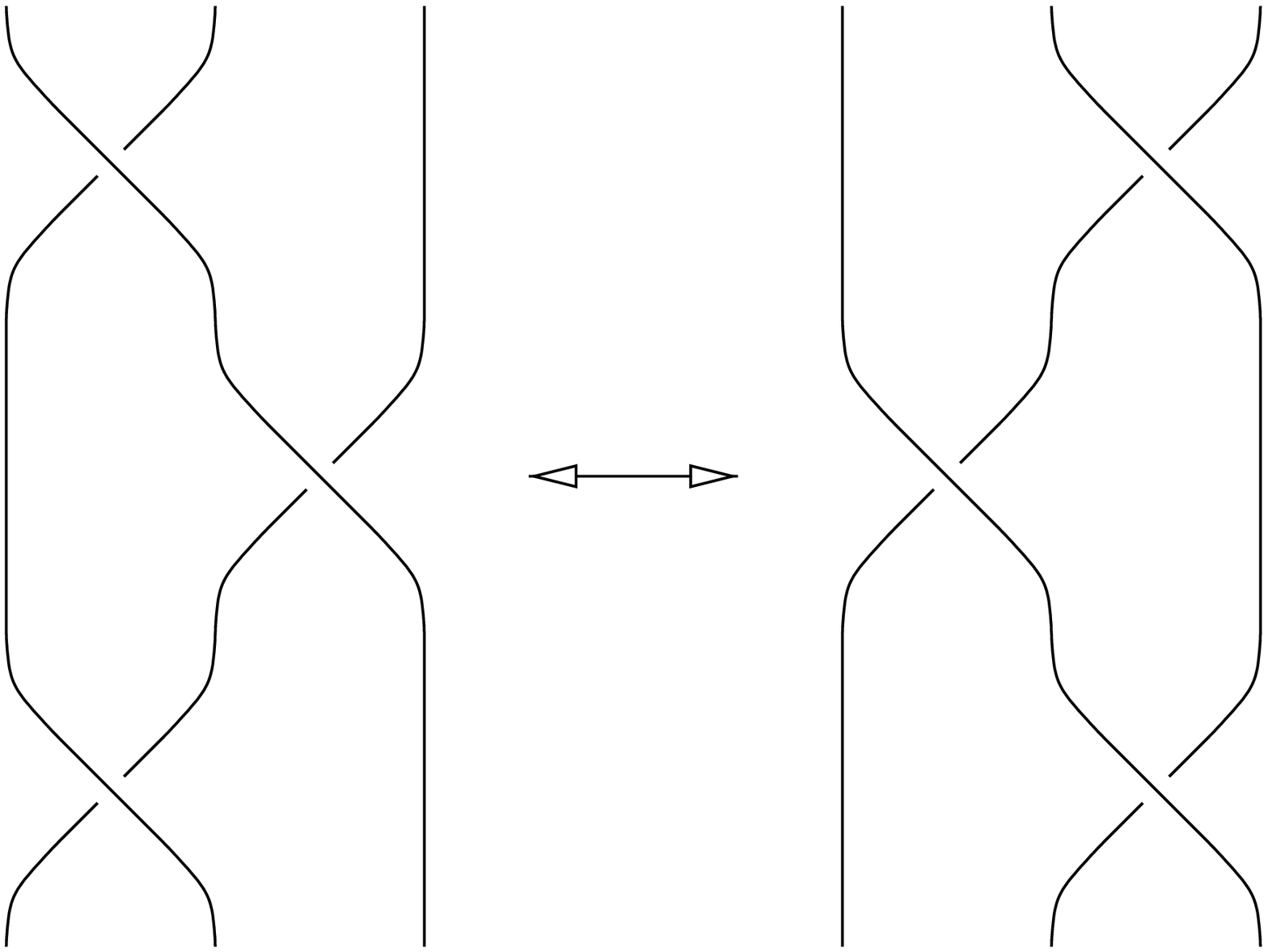}
}
\end{center}

\noindent This diagram represents the third Reidemeister move in
classical knot theory \cite{BZ}, and it gives the most important
relations in Artin's presentation of the braid group
\cite{Birman}. As a result, any  invertible solution of the
Yang--Baxter equation gives an invariant of braids.

\subsection{Shelves in Coalg and Solutions of the YBE}
We now
demonstrate the relationship between self-distributive maps in $\Coalg$ and
solutions to the Yang--Baxter equation.

\begin{figure}[htb]
\begin{center}
\mbox{
\epsfxsize=3.5in
\epsfbox{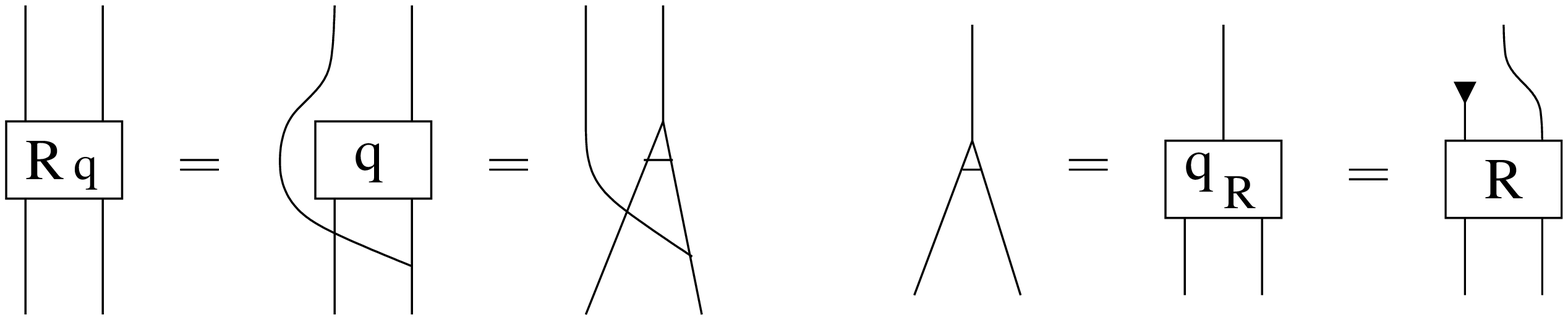}
}
\end{center}
\caption{ Solutions to YBE and shelves in Coalg }
\label{YBEtoQdef}
\end{figure}

\begin{definition} {\rm
Let $X$ be a coalgebra and $q: X \otimes X \rightarrow X$  a
linear map. Then the linear  map $R_q: X \otimes X \rightarrow X
\otimes X $ defined by
$$R_q= (1_X \otimes q)(\tau \otimes 1_X )(1_X \otimes \Delta)$$
is said to be
{\it induced from $q$}.

Conversely, let $R: X \otimes X \rightarrow X \otimes X $ be a
linear map. Then the linear map  $q_R:  X \otimes X \rightarrow X$
defined by $q_R=(\epsilon \otimes 1_X) R $ is said to be {\it
induced from $R$}. } \end{definition}

\begin{figure}[htb]
\begin{center}
\mbox{
\epsfxsize=2.5in
\epsfbox{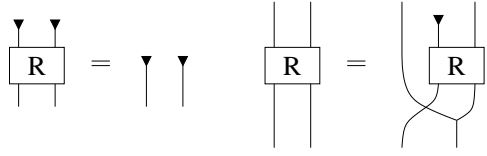}
}
\end{center}
\caption{Hypotheses of Theorem~\ref{risrqr}}
\label{YBEtoQlem}
\end{figure}

Diagrammatically, constructions of one of these maps from the
other are depicted in  Fig.~\ref{YBEtoQdef}. Our goal is to relate
solutions of the YBE and  self-distributive maps in certain categories via these
induced maps.

\begin{figure}[htb]
\begin{center}
\mbox{
\epsfxsize=4.5in
\epsfbox{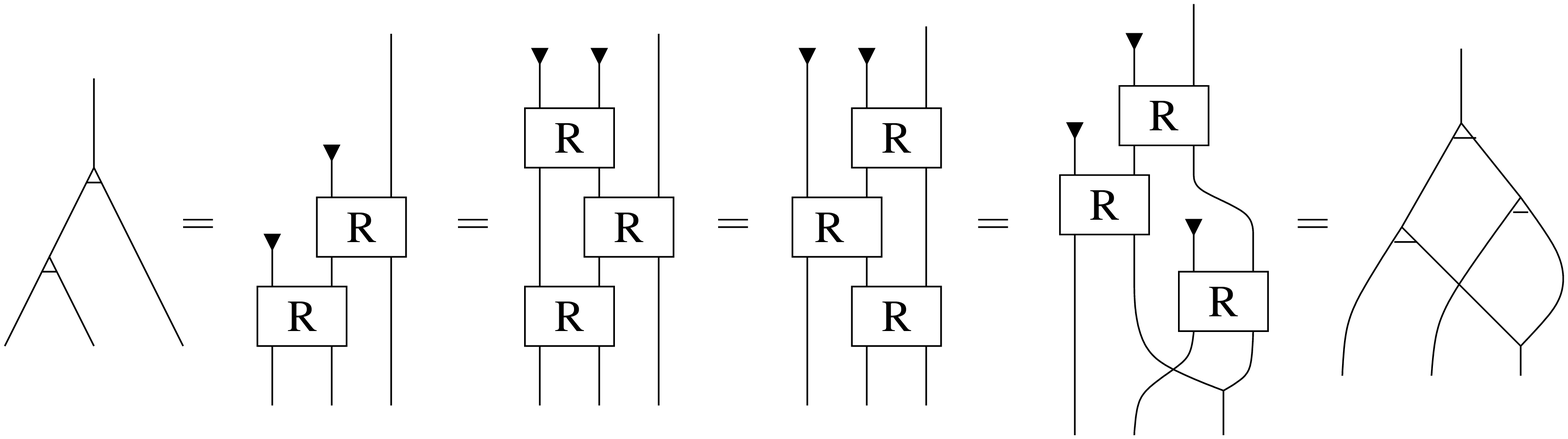}
}
\end{center}
\caption{Proof of Theorem~\ref{risrqr}}
\label{YBEtoQpf}
\end{figure}

\begin{figure}[htb]
\begin{center}
\mbox{
\epsfxsize=1.5in
\epsfbox{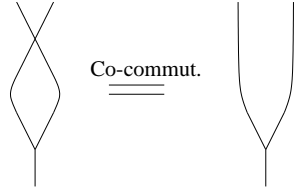} %%%%warning changed
}
\end{center}
\caption{Cocommutativity}
\label{cocom}
\end{figure}

\begin{theorem}\label{risrqr}
Let $R: X \otimes X \rightarrow X \otimes X $ be a solution to the
YBE on
a coalgebra $X$ with counit.
Suppose $R$ satisfies
$(\epsilon \otimes \epsilon) R = (\epsilon \otimes \epsilon)$ and
$R_{q_R}=R$. Then $(X, q_R)$ is a shelf in $\Coalg.$
\end{theorem}
{\it Proof.\/} The conditions in the assumption are presented in
Fig.~\ref{YBEtoQlem}. A proof is presented in Fig.~\ref{YBEtoQpf}.
$\Box$

\begin{figure}[htb]
\begin{center}
\mbox{
\epsfxsize=5.5in
\epsfbox{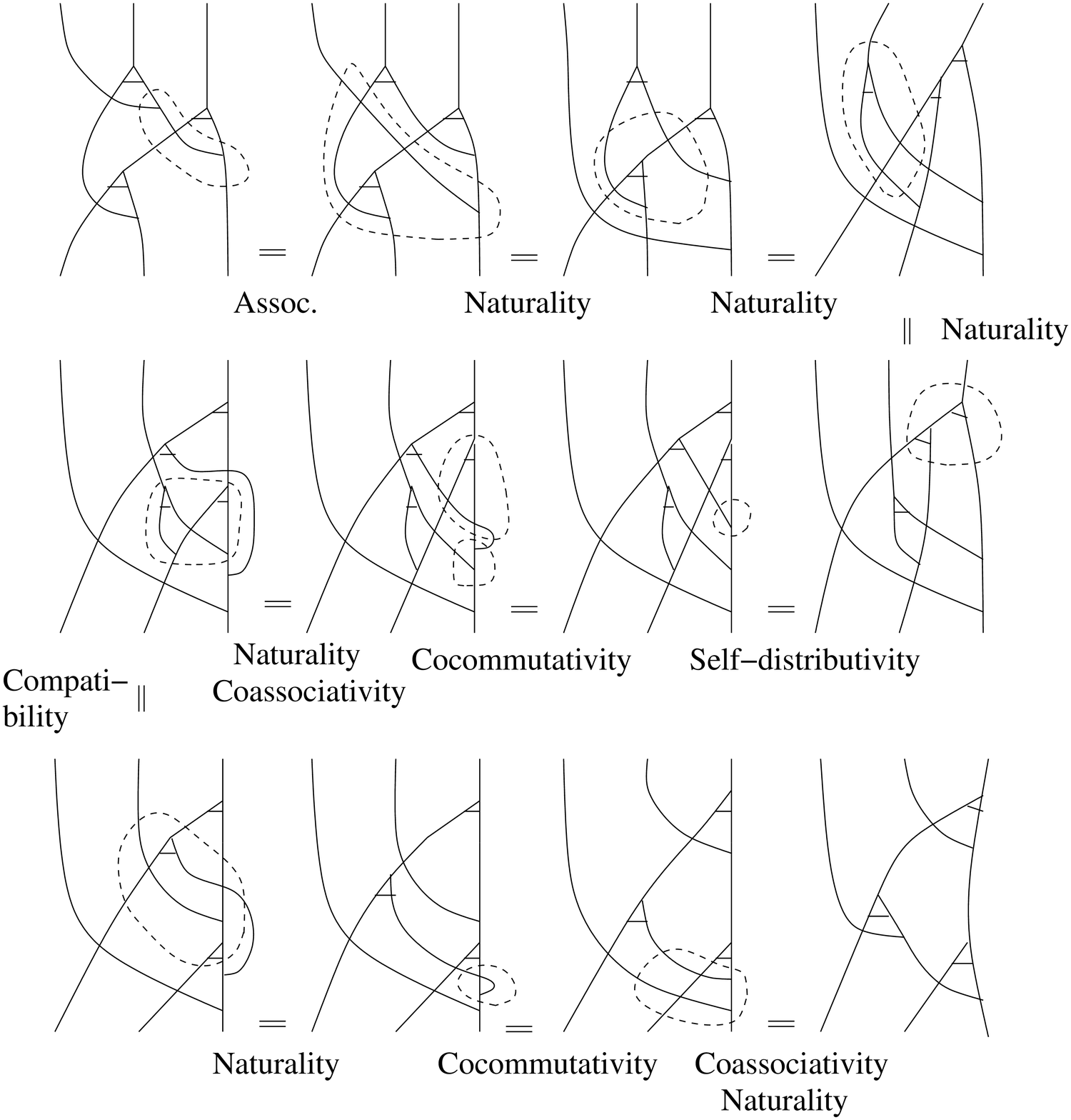}
}
\end{center}
\caption{Proof of Theorem~\ref{QtoYBEthm}}
\label{qtypeIII}
\end{figure}

\begin{theorem} \label{QtoYBEthm}
Let $X$ be an object in $\CoComCoalg$.
Suppose a self-distributive linear map
$q: X \otimes X \rightarrow X$ is compatible with comultiplication.
Then $R_q$ is a solution to the YBE.
\end{theorem}
{\it Proof.\/} The cocommutativity of $\Delta$ is depicted in
Fig.~\ref{cocom}.
A proof, then,
 is depicted in Fig.~\ref{qtypeIII}.
Note here the condition that $q$ is compatible with
comultiplication is that:
%coalgebra morphism is that it preserves the comultiplication:
$\Delta( q(a \otimes b))= q (a_{(1)} \otimes b_{(1)} ) \otimes q(
a_{(2)} \otimes b_{(2)})$ or, equivalently, $\Delta q = q (1
\otimes \tau \otimes 1 ) (\Delta \otimes \Delta)$. This is applied
in Fig.~\ref{qtypeIII} on the bottom row with the equal sign
indicated to follow from compatibility. $\Box$

Propositions~\ref{quandleQprop} and ~\ref{lieQprop}
and
Theorem~\ref{QtoYBEthm} imply the following:

\begin{corollary}
Let $q$ be a map
defined from a quandle/rack as in Proposition~\ref{quandleQprop}
or from a Lie algebra as in Proposition~\ref{lieQprop}. Then the
induced map $R_q$ is a solution to the YBE.
\end{corollary}

In the Lie algebra case, the map is given as follows:
$$R_q ((a,x)\otimes(b,y))=(b,y)\otimes(a, x)+(1,0) \otimes (0, [x,y]).$$
This appears, for example, in \cite{Alissa,MajidGreen}.

\begin{remark} {\rm
Next we focus on the case of  the adjoint map in Hopf algebras.
Remark ~\ref{RmkH} states that
the self-distributive map
$q(x\otimes y)=S(y_{(1)})xy_{(2)}$
is not compatible with comultiplication,
%%%SPELL PER AC AND MS 
and therefore,
Theorem~\ref{QtoYBEthm}
cannot be applied.  However, the induced map $R_q$ does, indeed,
satisfy YBE. This is of course for different reasons, and proved
in \cite{Woro},  which was interpreted
in \cite{Henn} as a restriction of a regular representation of the
universal $R$-matrix of a quantum double. Since it is of a great
interest why the same construction gives rise to solutions to YBE
for different reasons, we include their proofs in diagrams for
reader's convenience, and we specify two conditions from
\cite{Woro} in
our point of view, to construct $R_q$ from $q$,
and make a restatement of his theorem as follows: } \end{remark}

\begin{figure}[htb]
\begin{center}
\mbox{
\epsfxsize=2.5in
\epsfbox{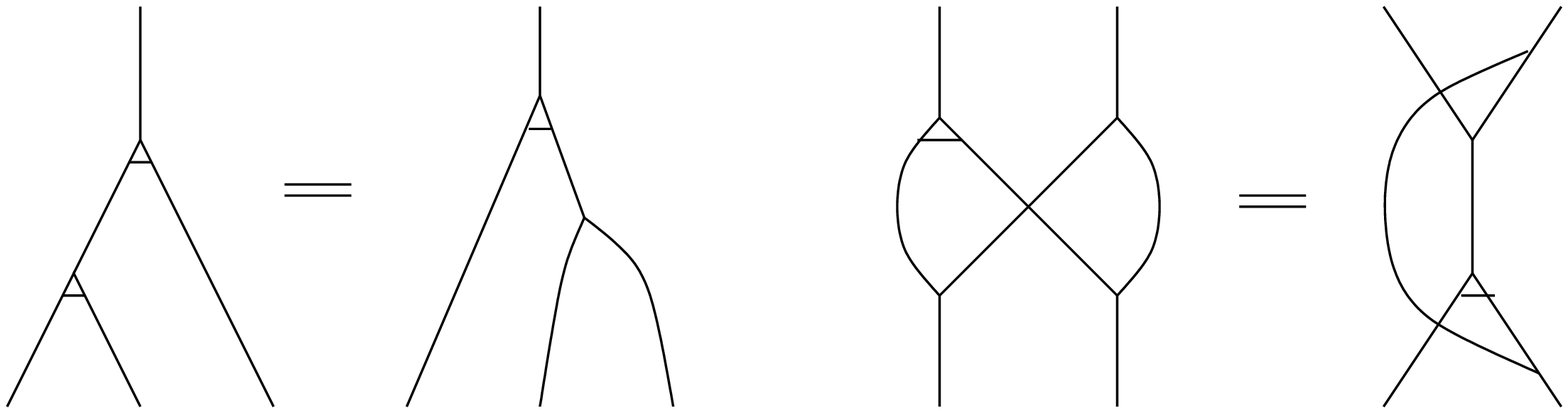}
}
\end{center}
\caption{Conclusion of Prop.~\ref{Hopfeqns}/ Hypothesis of Prop.~\ref{HopfYBEsoln}}
\label{Hopfqcond}
\end{figure}

\begin{figure}[htb]
\begin{center}
\mbox{
\epsfxsize=5.5in
\epsfbox{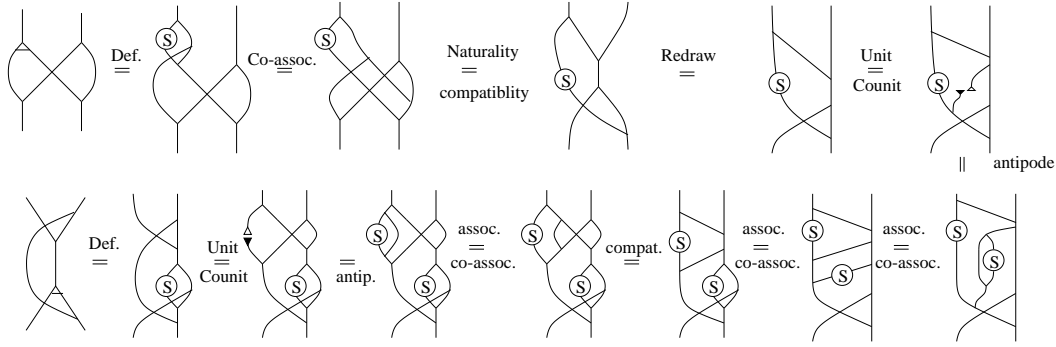}
}
\end{center}
\caption{Proof of Proposition~\ref{Hopfeqns}, second equation}
 \label{hopfYBElem}
\end{figure}

\begin{figure}[htb]
\begin{center}
\mbox{
\epsfxsize=5.5in
\epsfbox{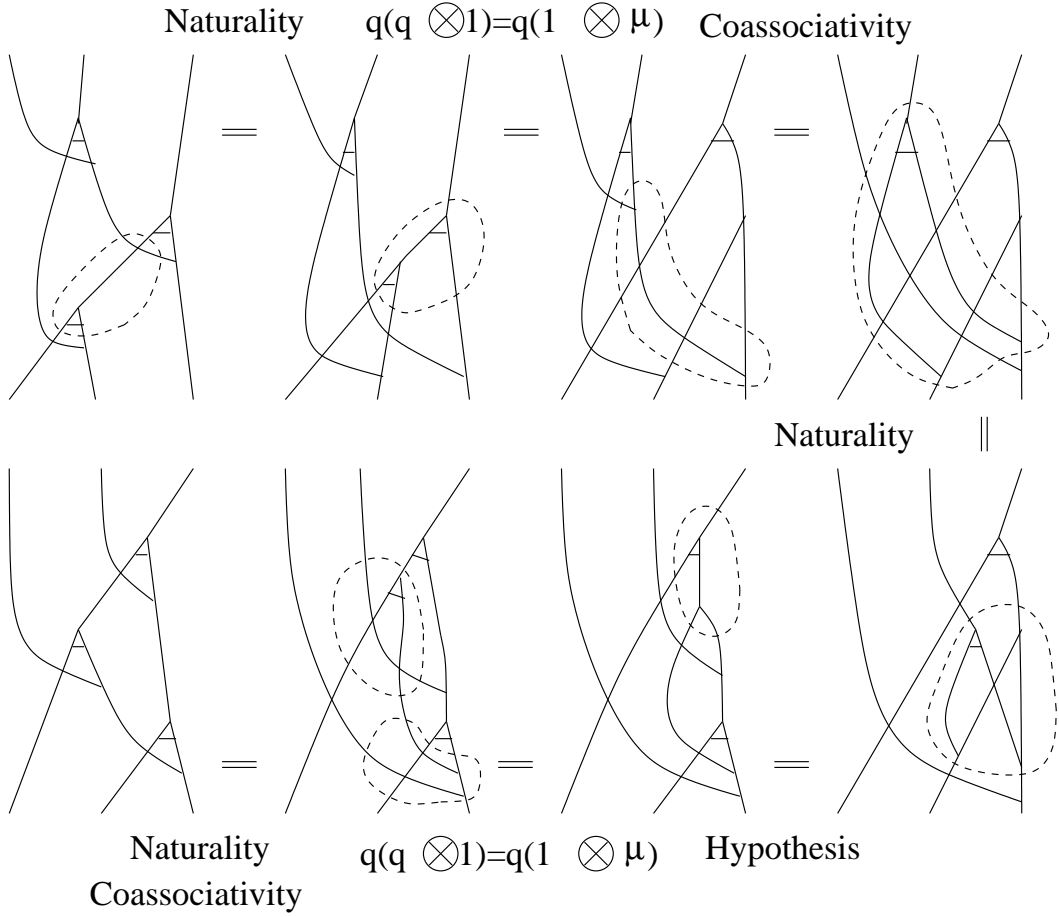}
}
\end{center}
\caption{Proof of Proposition~\ref{HopfYBEsoln}}
\label{HopftypeIII}
\end{figure}

\begin{proposition}
\label{Hopfeqns} In a Hopf algebra, let $q=\mu(1  \otimes \mu)(S
\otimes 1 \otimes 1 ) (\tau \otimes 1)(1 \otimes \Delta)$. Then
$q(q \otimes 1)=q( 1 \otimes \mu)$ and $(q \otimes \mu)(1 \otimes
\tau \otimes 1)(\Delta \otimes \Delta) = (1 \otimes \mu)(\tau
\otimes 1)(1 \otimes \Delta)(1 \otimes q)(\tau \otimes 1)(1
\otimes \Delta)$.
\end{proposition}
{\it Proof.}
The proofs are indicated in Figs.~\ref{HopfQlemma} and ~\ref{hopfYBElem}, respectively.
$\Box$

Recall from Section~\ref{Hopfsec} that in a Hopf algebra, the map
$q=\mu(1  \otimes \mu)(S \otimes 1 \otimes 1 ) (\tau \otimes 1)(1
\otimes \Delta)$ satisfies self-distributivity.

\begin{proposition}
\label{HopfYBEsoln} Suppose  $X$ is a Hopf algebra and $q$
is any
linear
map that
satisfies $q(q  \otimes 1)=q( 1 \otimes \mu)$ and $(q \otimes
\mu)(1 \otimes \tau \otimes 1)(\Delta \otimes \Delta) = (1 \otimes
\mu)(\tau \otimes 1)(1 \otimes \Delta)(1 \otimes q)(\tau \otimes
1)(1 \otimes \Delta)$. Then $R_q$ is a solution to the YBE.
\end{proposition}
{\it Proof.} The required conditions are depicted in
Fig.~\ref{Hopfqcond}. And the proof is given in
Fig.~\ref{HopftypeIII}. $\Box$

In particular, the above
proposition applies when $q(x \otimes y) = S(y_{(1)})x y_{(2)}.$

%
%Note that if we
%use a symbolic shorthand $x*y$ for $q(x \otimes y)$ and $xy$ for
%$m(x \otimes y)$,  then the two conditions are written as
%$(x*y)*z=x*(yz)$ and $y_{(1)}(x*y_{(2)})=xy$, which hold for
%conjugation in groups. A proof, then, is depicted in
%Fig.~\ref{HopftypeIII}.
%
%Among the  two conditions we needed, the first is proved earlier (see
%Fig.~\ref{HopfQlemma}) already, and a proof os the second is
%illustrated in Fig.~\ref{hopfYBElem}.
%$\Box$

% qmultilem.eps and Rladder.eps unified to hopfYBElem
%\begin{figure}[htb]
%\begin{center}
%\mbox{
%\epsfxsize=3.5in
%\epsfbox{hopfYBElem.eps}
%}
%\end{center}
%\caption{Self-distributive map and multiplication}
%\label{hopfYBElem}
%\end{figure}

\section{Graph Diagrams for Bialgebra Hochschild Cohomology}\label{hochsec}

The analogue  of group cohomology for associative algebras is
Hochschild cohomology. Then a natural question is, ``What is an
analogue of quandle cohomology for shelves
in $\Coalg$?"
Since we have developed diagrammatic methods
to study
self-distributivity in $\Coalg,$ we apply these methods to seek such a
cohomology theory, in combination with the interpretations of
cocycles in bialgebra cohomology in terms of deformation theory of
bialgebras. The first step toward this goal is to reestablish
diagrammatic methods for Hochschild cohomology in terms of graph
diagrams. Such approaches are found for homotopy Lie algebras and
operads~\cite{Markl}.
On the other hand, a diagrammatic method
using polyhedra for bialgebra cohomology was given in~\cite{MrSt}.
In this section we follow the
exposition in \cite{MrSt} of cocycles
that appear in bialgebra deformation theory, and establish tree
diagrams that can be used to prove cocycle conditions.

\begin{sloppypar}
First we recall the Hochschild cohomology for bialgebras from
\cite{MrSt}. Let $A=(V, \mu, \Delta)$ be a bialgebra over a field
$k$, where $\mu$, $\Delta$ are multiplication and
comultiplication, respectively, and $d_H: \Hom (V^{\otimes p},
V^{\otimes q} ) \rightarrow \Hom (V^{\otimes (p+1)}, V^{\otimes q}
)$ is the Hochschild differential
$$d_H(f)=\mu( 1  \otimes f ) + \sum_{i=0}^{p-1} (-1)^{i+1}
f(1 ^i \otimes \mu \otimes 1^{n-i-1})
+(-1)^{p+1} \mu (f \otimes 1) $$
where the left and right module structures are
given by multiplication.
Dually $d_C:  \Hom (V^{\otimes p}, V^{\otimes q} )
\rightarrow \Hom (V^{\otimes p}, V^{\otimes (q+1)} )$
denotes the coHochschild differential.
These define the total complex
$(C^*_b(A;A), D)$, where
$C^n_b(A;A)=\oplus_{i=1}^n \Hom(V^{\otimes (n-i+1)} , V^{\otimes i})$.
For example, for a $1$-cochain $f \in \Hom(V , V)$,
$d_H(f)(x\otimes y)=x f(y) - f(xy) + f(x)y$ and
$d_C(f)(x) = x_{(1)} \otimes f(x_{(2)}) - f(x)_{(1)} \otimes f(x)_{(2)}+
 f(x_{(1)} )\otimes x_{(2)}$.
\end{sloppypar}

For the rest of this section, we establish graph
diagrams for Hochschild cohomology and review their aspects in
deformation theory of bialgebras.

\subsection{Graph Diagrams for Hochschild Differentials}

\begin{figure}[htb]
\begin{center}
\mbox{
\epsfxsize=4in
\epsfbox{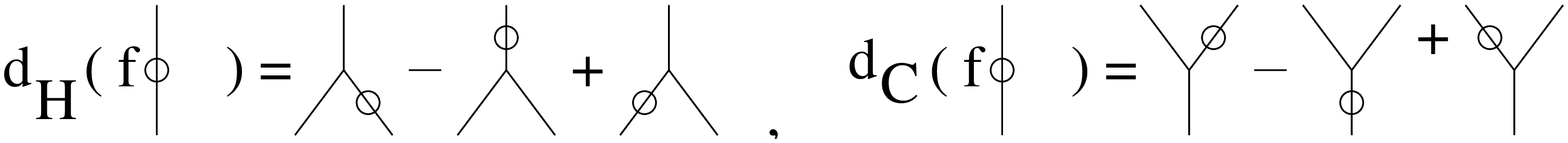}
}
\end{center}
\caption{Hochschild $1$-differentials}
\label{biHoch}
\end{figure}

\begin{figure}[htb]
\begin{center}
\mbox{ \epsfxsize=.5in
\epsfbox{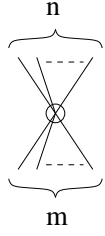}}
\end{center}
\caption{Hochschild $(m+n)$-cochains}
\label{hochcocy}
\end{figure}

A $1$-cochain $f \in \Hom(V, V)$ is represented by a circle on a vertical
segment as shown in Fig.~\ref{biHoch}, where the images of $f$ under
the first differentials $d_H(f)$ and $d_C(f)$, 
%%%deleted (x) as per AC comment
as computed
 above,
are also depicted.
In general, a $(m+n-1)$-cochain in $\Hom(V^{\otimes m}, V^{\otimes n})$
is represented by a diagram
in Fig.~\ref{hochcocy}.

For $(\phi_1, \phi_2)$, where $\phi_1 \in \Hom (V^{\otimes 2}, V
)$ and
  $\phi_2 \in  \Hom ( V, V^{\otimes 2}) $, the differentials are
    \begin{eqnarray}
  d_H (\phi_1) & = & \mu ( 1 \otimes \phi_1) - \phi_1(\mu \otimes 1) +
  \phi_1 ( 1 \otimes \mu) - \mu (\phi_1 \otimes 1),  \label{Hoch1} \\
   d_C(\phi_1) & = &
   (\mu \otimes \phi_1) \tau_2 (\Delta \otimes \Delta)
   - \Delta (\phi_1) + ( \phi_1 \otimes \mu ) \tau_2 (\Delta \otimes \Delta), \label{Hoch21} \\
   d_H (\phi_2) & =&
   (\mu \otimes \mu )\tau_2 (\Delta \otimes \phi_2) - \phi_2 \mu
   +  (\mu \otimes \mu )\tau_2 (\phi_2 \otimes \Delta) ,  \label{Hoch22} \\
   d_C (\phi_2) & = & (1 \otimes \phi_2) \Delta - (\Delta \otimes 1) (\phi_2)
   + (1 \otimes \Delta)(\phi_2) - (\phi_2 \otimes 1) \Delta ,  \label{Hoch3}
   \end{eqnarray}
where $\tau_2$ is the homomorphism induced from the transposition
of the second and the third factors.
The $2$-cocycle conditions are $ d_H (\phi_1)=0$,
$ d_C(\phi_1)= d_H (\phi_2) $, and $ d_C (\phi_2)=0$.
The differential $D$ of the total complex is $D=d_H-d_C$,
$D(\phi_1, \phi_2)=d_H(\phi_1) + [d_H (\phi_2)- d_C(\phi_1)] - d_C (\phi_2)$.

\begin{figure}[htb]
\begin{center}
\mbox{
\epsfxsize=5in
\epsfbox{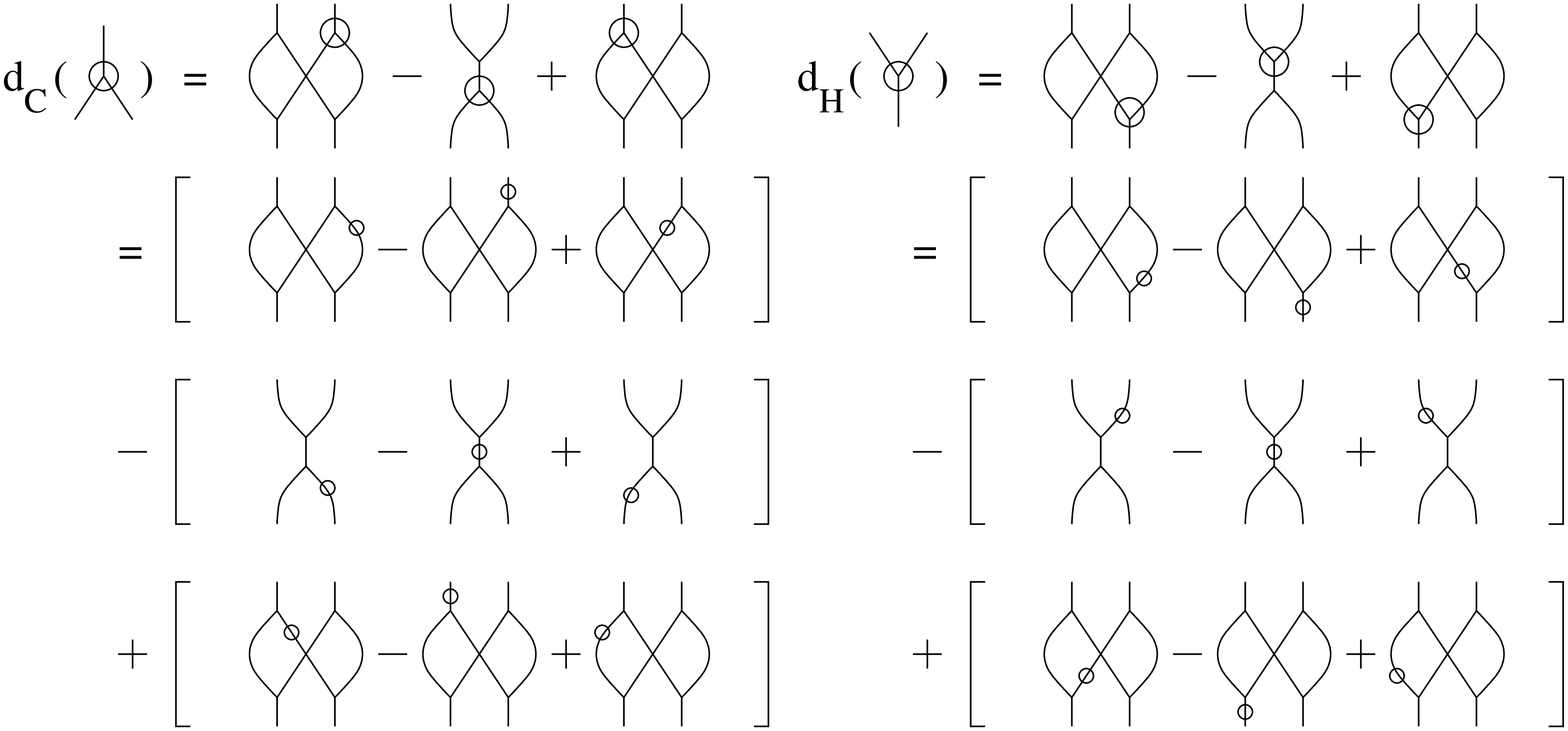}
}
\end{center}
\caption{Hochschild $2$-differentials}
\label{d2hoch12}
\end{figure}

We demonstrate a proof that $(\phi_1, \phi_2)=   (d_H(f), d_C(f))$
satisfies $ d_C(\phi_1)= d_H (\phi_2) $ using graph diagrams.
First, we use encircled vertices as
depicted in Fig.~\ref{hochcocy} to represent an element of $\Hom
(V^{\otimes m}, V^{\otimes n} )$. Then $ d_C(\phi_1)$ and $d_H
(\phi_2) $ are represented on the top line of
Fig.~\ref{d2hoch12}.
Substituting  $(\phi_1, \phi_2)= (d_H(f), d_C(f))$, that are
represented diagrammatically as in Fig.~\ref{biHoch}, we
perform  diagrammatic computations as in the rest of
Fig.~\ref{d2hoch12},
and the equality follows
because multiplication and comultiplication are compatible. In
particular each diagram in the left of the figure for which the
vertex is external to the operations corresponds to a similar
diagram on the right, but the correspondence is given after
considering the compatible structures.

\begin{figure}[htb]
\begin{center}
\mbox{
\epsfxsize=3.2in
\epsfbox{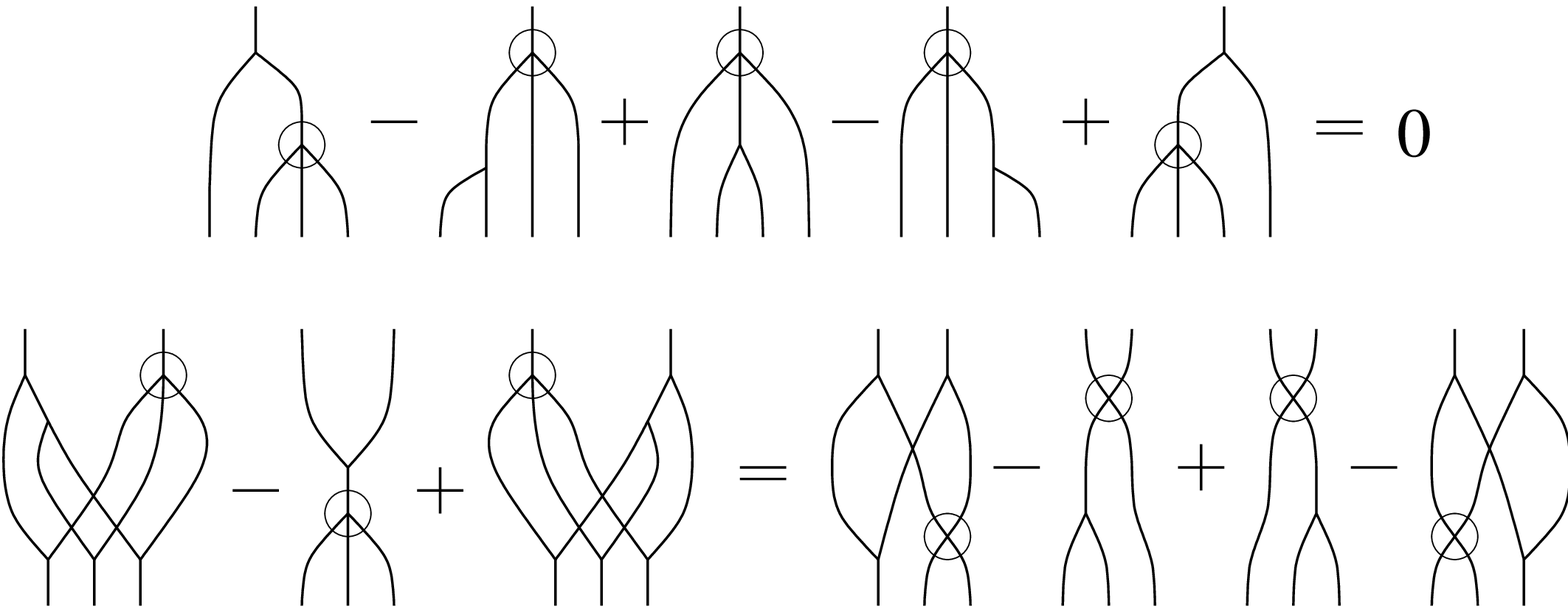}
}
\end{center}
\caption{Hochschild $3$-cocycle conditions}
\label{Hoch31}
\end{figure}

For $3$-cochains $\psi_i \in \Hom (V^{\otimes 3}, V )$, $\psi_2
\in \Hom (V^{\otimes 2}, V^{\otimes 2} )$ and $\psi_3 \in \Hom (V,
V^{\otimes 3} )$, the $3$-cocycle condition is explicitly written
as $d_H(\psi_1)=0$, $d_C(\psi_1)=d_H(\psi_2)$,
$d_C(\psi_2)=d_H(\psi_3)$, and $d_C(\psi_3)=0$, see~\cite{MrSt}
\begin{eqnarray*}
  d_H (\psi_1) & = & \mu ( 1 \otimes \psi_1) - \psi_1(\mu \otimes 1^2) +
  \psi_1(1 \otimes \mu \otimes 1) -
  \psi_1 ( 1^2 \otimes \mu) + \mu (\psi_1 \otimes 1),  \\
   d_C(\psi_1) & = &
   (\mu (1 \otimes \mu) \otimes \psi_1) \tau (\Delta \otimes \Delta \otimes \Delta)
   - \Delta (\psi_1) +
    ( \psi_1 \otimes \mu (\mu \otimes 1) ) \tau (\Delta \otimes \Delta \otimes \Delta), \\
   d_H (\psi_2) & =&
   (\mu \otimes \mu )\tau_2 (\Delta \otimes \psi_2) -
   \psi_2 (\mu \otimes 1) + \psi_2 (1 \otimes \mu)
   -  (\mu \otimes \mu )\tau_2 (\psi_2 \otimes \Delta) ,  \\
   d_C (\psi_2) & = & (\mu \otimes \psi_2)\tau_2 (\Delta\otimes \Delta)
    - (\Delta \otimes 1) (\psi_2)
   + (1 \otimes \Delta)(\psi_2) - (\psi_2 \otimes \mu) \tau_2 (\Delta\otimes \Delta) ,  \\
    d_H (\psi_3) & =&
   (\mu \otimes \mu  \otimes \mu)\tau' (
%%%%\Delta (1 \otimes \Delta)\otimes \psi_3)
(1 \otimes \Delta) \Delta \otimes \psi_3) %%%%%%
   - \psi_3 (\mu)
   +  (\mu \otimes \mu \otimes \mu)\tau' (\psi_3 \otimes
%%%\Delta (\Delta \otimes 1)) ,  \\
(\Delta \otimes 1)\Delta) , \\ %%%%
   d_C (\psi_3) & = & (1 \otimes \psi_3) \Delta - (\Delta \otimes 1^2) (\psi_3)
   + (1 \otimes \Delta \otimes 1)(\psi_3)
   - (1^2 \otimes \Delta)(\psi_3) + (\psi_3 \otimes 1) \Delta ,
   \end{eqnarray*}
where $\tau=\tau_4 \tau_3 \tau_2$
and $\tau'=\tau_5 \tau_2 \tau_3$.
In general, $\tau_i$ indicates the transposition of the  $i\/$th and $(i+1)\/$st factors; the notation is used when type-setting gets complicated.

\begin{figure}[htb]
\begin{center}
\mbox{
\epsfxsize=2.3in
\epsfbox{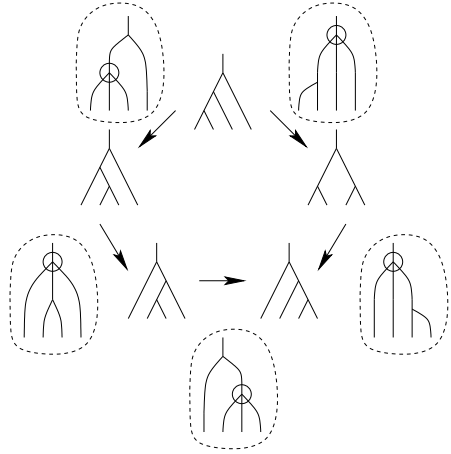}
}
\end{center}
\caption{Hochschild $3$-cocycles as movies, part I}
\label{Hoch32}
\end{figure}

The first two $3$-cocycle conditions, $d_H (\psi_1)=0$ and
$d_C(\psi_1)=d_H(\psi_2)$, are depicted in Fig.~\ref{Hoch31}. Note
that the first is the pentagon identity for associativity. In
particular, $\psi_1$ can be regarded as an obstruction to
associativity. The morphism $\psi_1$ is assigned the difference
between the two diagrams that represent the two expressions
$(ab)c$ and $a(bc)$. Thus $\psi_1$ and its diagram are assigned to
the change of diagrams corresponding to associativity, and can be
seen to form an actual pentagon, as depicted in Fig.~\ref{Hoch32}.

\begin{figure}[htb]
\begin{center}
\mbox{
\epsfxsize=2.7in
\epsfbox{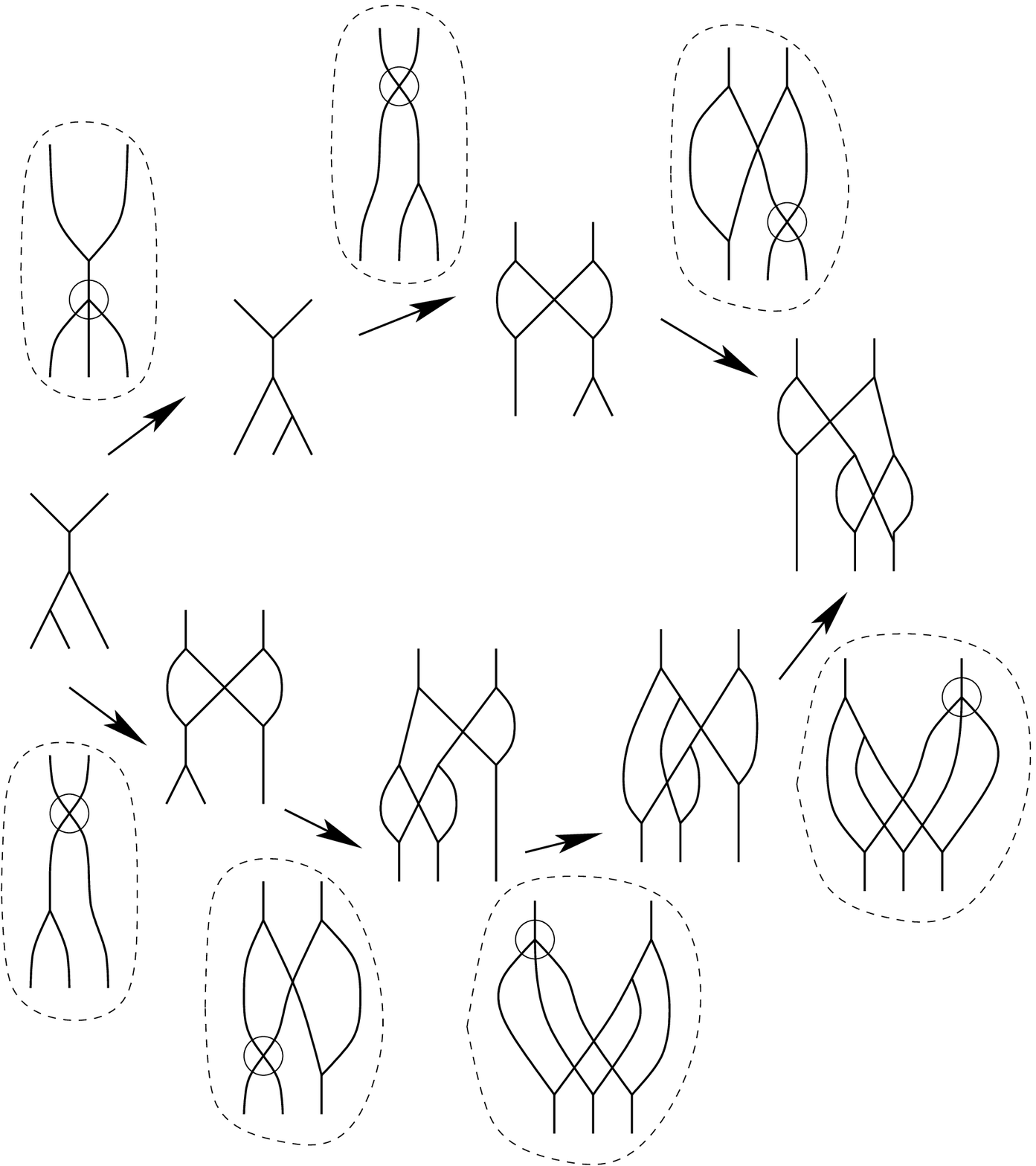}
}
\end{center}
\caption{Hochschild $3$-cocycles as movies, part Il}
\label{Hoch33}
\end{figure}

Similarly, the second condition $d_C(\psi_1)=d_H(\psi_2)$ can be
represented as sequence of applications of the associativity and
compatibility conditions as depicted in Fig.~\ref{Hoch33}.
Furthermore, the relations $d_C(\psi_2)=d_H(\psi_3)$ and
$d_C(\psi_3)=0$ can be obtained by turning the equations in
Fig.~\ref{Hoch31} upside-down. Similarly, the ``movie-moves" in
Figs.~\ref{Hoch32} and \ref{Hoch33} can be turned upside-down.
Thus,  $d_C(\psi_3)=0$ when the pentagon identity for
coassociativity holds, and $d_C(\psi_2)=d_H(\psi_3)$ when
compatibility and coassociativity are compared.

\subsection{Review of Cocycles in Deformation Theory}

Next we follow \cite{MrSt} for deformation of bialgebras.
A deformation of $A=(V, \mu, \Delta) $ is a
$k[[t]]$-bialgebra $A_t=(V_t, \mu_t, \Delta_t)$, where $V_t=V \otimes k[[ t ]]$
and $ A_t/(tA_t) \cong A$. Deformations of $\mu$ and  $\Delta$ are given by
$\mu_t= \mu + t \mu_1 + \cdots + t^n \mu_n + \cdots : V_t \otimes V_t \rightarrow V_t$
and $\Delta_t = \Delta + t \Delta_1
+ \cdots + t^n \Delta_n + \cdots : V_t \rightarrow V_t  \otimes V_t$
where $\mu_i: V \otimes V \rightarrow V $, $\Delta_i : V \rightarrow V  \otimes V$,
$i=1, 2, \cdots$, are sequences of maps.
Suppose $\bar{\mu}=\mu + \cdots + t^n \mu_n$ and
$\bar{\Delta} =\Delta + \cdots + t^n \Delta_n$ satisfy the bialgebra
conditions (associativity, compatibility, and coassociativity) mod $t^{n+1}$,
and suppose
that  there exist $\mu_{n+1}: V \otimes V \rightarrow V$ and
$\Delta_{n+1}: V \rightarrow V \otimes V$ such that
$\bar{\mu}+t^{n+1} \mu_{n+1}$ and $\bar{\Delta}+ t^{n+1} \Delta_{n+1}$
satisfy the bialgebra conditions mod $t^{n+2}$.
Define  $\psi_1 \in \Hom(V^{\otimes 3}, V)$,
$\psi_2\in \Hom(V^{\otimes 2}, V^{\otimes 2})$,
and $\psi_3 \in \Hom(V, V^{\otimes 3})$
by
\begin{eqnarray}
\bar{\mu} (\bar{\mu} \otimes 1) - \bar{\mu}(1 \otimes \bar{\mu})
&=& t^{n+1} \psi_1 \quad {\rm mod}\  t^{n+2} , \label{hoch2d1} \\
\bar{\Delta} \bar{\mu} -
(\bar{\mu} \otimes \bar{\mu})\tau_2 (\bar{\Delta}\otimes \bar{\Delta})
&=&  t^{n+1} \psi_2 \quad {\rm mod}\  t^{n+2} , \label{hoch2d2}\\
(\bar{\Delta} \otimes 1) \bar{\Delta} - (1 \otimes \bar{\Delta})\bar{\Delta}
&=&  t^{n+1} \psi_3 \quad {\rm mod}\  t^{n+2} . \label{hoch2d3}
\end{eqnarray}
For the associativity of $\bar{\mu}+t^{n+1} \mu_{n+1}$ mod
$t^{n+2}$ we obtain:
$$(\bar{\mu}+t^{n+1} \mu_{n+1})((\bar{\mu}+t^{n+1} \mu_{n+1}) \otimes 1)-
(\bar{\mu}+t^{n+1} \mu_{n+1})(1 \otimes (\bar{\mu}+t^{n+1}
\mu_{n+1}))=0 \  {\rm mod}\  t^{n+2} $$ which is equivalent by
degree calculations to:
$$d_H(\mu_{n+1})=\mu (1 \otimes \mu_{n+1})-\mu_{n+1}(\mu \otimes 1)
 + \mu_{n+1}(1 \otimes \mu_{n+1}) - \mu (\mu_{n+1} \otimes 1)=\psi_1. $$
 Similarly, we obtain:
$(\psi_1, \psi_2, \psi_3)=D(\mu_{n+1}, \Delta_{n+1}) $. The
cochains $(\psi_1, \psi_2, \psi_3)$, defined by
deformations~(\ref{hoch2d1},\ref{hoch2d2},\ref{hoch2d3}) then,
satisfy the $3$-cocycle condition $D(\psi_1, \psi_2, \psi_3)=0$.
This concludes the
review  of deformation for the $2$-cocycle conditions cited from
\cite{MrSt}.

\section{Towards a Cohomology Theory for Shelves in Coalg}\label{cohsec}

Let $(X,q)$ be a  coalgebra with a self-distributive
linear
map.
In this section we present
low-dimensional
cocycle conditions for $q$.
We justify our cocycle conditions through the use of analogy
with Hochschild bialgebra cohomology using diagrammatics and the
deformation theories reviewed in the preceding section.
Both
analogies are used interchangeably throughout this section,
both in definitions and computations.

\subsection{Chain Groups}

Following the diagrammatics of the preceding
section,
we define chain groups, for positive integers $n$ and $i =1,
\ldots, n$  by:
\begin{eqnarray*}
C^{n, i}_{\rm sh} (X;X) & =&  \Hom(X^{\otimes (n+1-i)}, X^{\otimes i} ), \\
C^n_{\rm sh} (X;X) & =& \oplus_{i=1}^n C^{n, i}_{sh} (X;X).
\end{eqnarray*}
Specifically, the chain groups in low dimensions of our concern
are:
\begin{eqnarray*}
C^1_{\rm sh}(X;X)&=&  \Hom(X, X) , \\
C^2_{\rm sh}(X;X)&=& \Hom(X^{\otimes 2}, X)\oplus  \Hom(X, X^{\otimes 2}), \\
C^3_{\rm sh}(X;X)&=& \Hom(X^{\otimes 3}, X)\oplus
\Hom(X^{\otimes 2}, X^{\otimes 2})\oplus \Hom(X, X^{\otimes 3}).
\end{eqnarray*}

\begin{figure}[htb]
\begin{center}
\mbox{
\epsfxsize=4in
\epsfbox{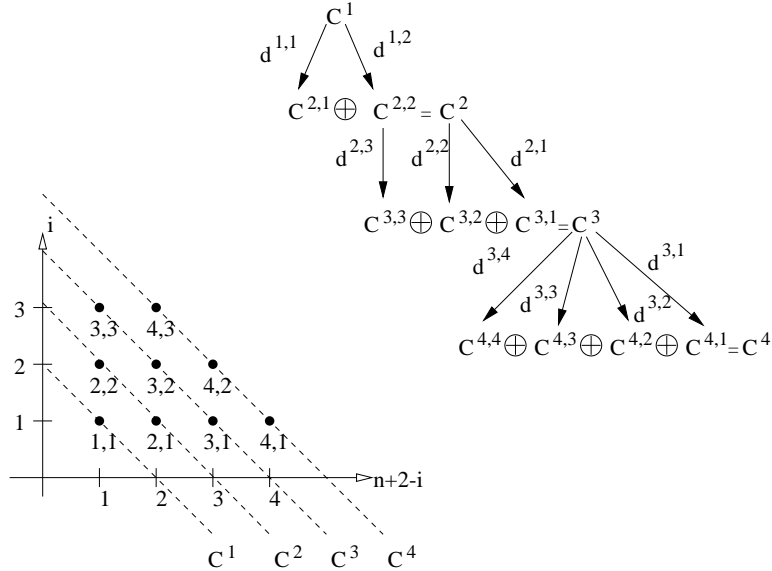}
}
\end{center}
\caption{The lattice of chain groups and differentials  }
\label{complex}
\end{figure}

To help keep track of the chain groups and their indices, we
include the diagram in Fig.~\ref{complex}.
The chain groups $C^{n,i}$ are located at position
$(n+2-i,i)$ in the positive quadrant of the integer lattice. The chain groups $C^j$ are the direct sum of the groups along lines of slope $(-1)$. Differentials in the figure are indicated by arrows that point to the target groups. The differential $d^{2,3}$ has as its source the summand $C^{2,2} \subset C^2$ as indicated.

In the remaining sections we will define differentials that are
homomorphisms between the chain groups:
$$d^{n, i} :  C^n_{\rm sh} (X;X) %\Hom(V^{\otimes (n+1-i)}, V^{\otimes i} )
 \rightarrow C^{n+1, i}_{\rm sh}(X;X) (= \Hom(X^{\otimes (n+2-i)}, X^{\otimes i} ) )$$
and will be  defined individually for $n=1,2,3$ and $i=1, \ldots,
n+1$, and
\begin{eqnarray*}
 D_1 & =& d^{1, 1}-d^{1,2}: C^1_{\rm sh} (X;X) \rightarrow C^2_{\rm sh} (X;X),\\
 D_2&=& d^{2, 1} + d^{2, 2}+d^{2, 3}:  C^2_{\rm sh} (X;X) \rightarrow C^3_{\rm sh} (X;X),\\
 D_3&=& d^{3, 1} +d^{3, 2} +d^{3, 3} +d^{3, 4}:
  C^3_{\rm sh} (X;X) \rightarrow C^3_{\rm sh} (X;X).
\end{eqnarray*}

\subsection{First Differentials}

We take
$$d^{1,2}: \Hom(X;X) (=C^{1, 1}_{\rm sh} (X;X) )  \rightarrow
 \Hom(X, X^{\otimes 2}) (=C^{2,2}_{\rm sh} (X;X) ) $$
 to be the coHochschild differential for the comultiplication
 $d^{1,2} (f) = (1 \otimes f) \Delta - \Delta f + (f \otimes 1)\Delta$.
Again by analogy with the differential for multiplication, we
take:
$$d^{1,1}:  \Hom(X, X)(=C^{1, 1}_{\rm sh} (X;X) )  \rightarrow
\Hom(X^{\otimes 2}, X) (=C^{2,1}_{\rm sh} (X;X) ) $$
to be  $d^{1,1}(f)= q(1 \otimes f) - fq + q(f \otimes 1)$.
Then define $D_2: C^1_{\rm sh} (X;X)  \rightarrow C^2_{\rm sh} (X;X)$ by
$D_1=d^{1,1}- d^{1,2}$.
%%%D_2 mapsto D_1 above per AC

\subsection{Second Differentials}

We derive second differentials by analogy with deformation theory,
and then show that our definitions carry through in diagrammatics.

Recall that the self-distributivity, compatibility, and
coassociativity are written as:
\begin{eqnarray*}
 q (q \otimes 1) &=&  q (q \otimes q)  \tau_2
 ( 1\otimes 1  \otimes \Delta), \\
 \Delta q &=& (q \otimes q) \tau_2 ( \Delta \otimes \Delta ), \\
 (\Delta \otimes 1) \Delta &=& (1 \otimes \Delta ) \Delta .
\end{eqnarray*}
where $\tau_2$ is  the transposition acting on the second and third tensor factors.
As before let $X_t=X \otimes k[[t]]$ and suppose we have partial
deformations $\bar{q} = q + \cdots + t^n q_n$ and $\bar{\Delta}
=\Delta + \cdots + t^n \Delta_n$ satisfying the above three
conditions mod $t^{n+1}$, and suppose there are $q_{n+1}$ and
$\Delta_{n+1}$ such that $\bar{q}+ q_{n+1}$ and $\bar{\Delta} +
\Delta_{n+1}$ satisfy the three conditions mod $t^{n+2}$.

Setting
\begin{eqnarray*}
 \bar{q} (\bar{q} \otimes 1)  -\bar{q} (\bar{q} \otimes \bar{q})\tau_2
 ( 1 \otimes 1  \otimes \bar{\Delta})  &=&
 t^{n+1} \xi_1 \quad {\rm mod}\  t^{n+2} , \\
 \bar{\Delta} \bar{q} -( \bar{q} \otimes  \bar{q}) \tau_2
  ( \bar{ \Delta } \otimes  \bar{\Delta} )&=&
 t^{n+1} \xi_2 \quad {\rm mod}\  t^{n+2} , \\
 ( \bar{\Delta} \otimes 1) \bar{ \Delta} - (1 \otimes  \bar{\Delta} )  \bar{\Delta} &=&
  t^{n+1} \xi_3 \quad {\rm mod}\  t^{n+2} ,
\end{eqnarray*}
we obtain:
\begin{eqnarray*}
\left[ q (q_{n+1} \otimes 1) + q_{n+1} (q \otimes 1) \right]
 -
 \left[
q_{n+1} (q \otimes q)\tau_2( 1 \otimes 1  \otimes \Delta)  +
q ( q_{n+1}\otimes q)\tau_2( 1 \otimes 1  \otimes \Delta)   + \right. & & \\
\left.
q (q \otimes q_{n+1})\tau_2( 1 \otimes 1  \otimes \Delta)  +
q (q \otimes q)\tau_2( 1 \otimes 1  \otimes \Delta_{n+1})
\right]
&=&
  \xi_1  , \\
 \left[ \Delta q_{n+1} + \Delta_{n+1} q \right]
 -
 \left[
 (q_{n+1}  \otimes q) \tau_2 ( \Delta \otimes \Delta ) +
 (q \otimes q_{n+1} ) \tau_2 ( \Delta \otimes \Delta ) + \right. & & \\
\left.
 (q \otimes q) \tau_2 ( \Delta_{n+1}  \otimes \Delta ) +
 (q \otimes q) \tau_2 ( \Delta \otimes \Delta_{n+1}  )
 \right]
  &=&
  \xi_2 , \\
\left[
 (\Delta_{n+1} \otimes 1) \Delta +
   (\Delta \otimes 1) \Delta_{n+1}
   \right]
   -
   \left[
 (1 \otimes \Delta_{n+1} ) \Delta +
 (1 \otimes \Delta ) \Delta_{n+1}
\right]
 &=&
  \xi_3 .
\end{eqnarray*}

\begin{sloppypar}
A natural requirement is  $D_2(q_{n+1}, \Delta_{n+1})=(\xi_1, \xi_2, \xi_3)$, so we define
$D_2: C^2_{\rm sh}(X;X) \rightarrow C^3_{\rm sh}(X;X)$ by
$D_2=d^{2,1} + d^{2,2} + d^{2,3}$, where
\end{sloppypar}

\begin{eqnarray*}
d^{2,1} (\eta_1, \eta_2) & = &
 \left[ q (\eta_1 \otimes 1) + \eta_1 (q \otimes 1) \right]
 -
 \left[
\eta_1 (q \otimes q)\tau_2( 1 \otimes 1  \otimes \Delta)  +
q ( \eta_1 \otimes q)\tau_2( 1 \otimes 1  \otimes \Delta)   \right. \\
 &  & \left. \quad +
 q (q \otimes \eta_1)\tau_2( 1 \otimes 1  \otimes \Delta)  +
q (q \otimes q)\tau_2( 1 \otimes 1  \otimes \eta_2)
\right] \\
d^{2,2} (\eta_1, \eta_2) & = &  \left[ \Delta \eta_1 + \eta_2 q \right]
 -
 \left[
 (\eta_1 \otimes q) \tau_2 ( \Delta \otimes \Delta ) +
 (q \otimes \eta_1 ) \tau_2 ( \Delta \otimes \Delta ) \right. \\
 &  & \left. \quad +
 (q \otimes q) \tau_2 ( \eta_2  \otimes \Delta ) +
 (q \otimes q) \tau_2 ( \Delta \otimes \eta_2  )
 \right]
 \\
 d^{2,3} (\eta_1, \eta_2) & = &
\left[
 (\eta_2 \otimes 1) \Delta +
   (\Delta \otimes 1) \eta_2
   \right]
  -
   \left[
 (1 \otimes \eta_2 ) \Delta +
 (1 \otimes \Delta ) \eta_2
\right]  .
\end{eqnarray*}
In fact, $d^{2,3}= d_C$, the same as the coHochschild $2$-differential for
the comultiplication.

\begin{figure}[htb]
\begin{center}
\mbox{
\epsfxsize=2in
\epsfbox{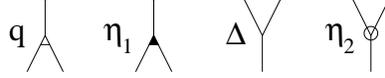}
}
\end{center}
\caption{Diagrams for $2$-cochains}
\label{cocydiag}
\end{figure}

The diagrammatic conventions for $q$, a $2$-cochain $\eta_1 \in
\Hom(X^{\otimes 2}, X)$, and $\Delta$, a $2$-cochain $\eta_2 \in
\Hom(X, X^{\otimes 2})$ are depicted from left to right,
respectively, in Fig.~\ref{cocydiag}.

\begin{figure}[htb]
\begin{center}
\mbox{
\epsfxsize=3in
\epsfbox{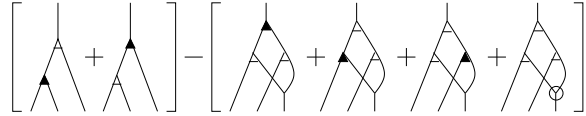}
}
\end{center}
\caption{The first $2$-differential $d^{2,1}$}
\label{diff1}
\end{figure}

\begin{figure}[htb]
\begin{center}
\mbox{
\epsfxsize=3in
\epsfbox{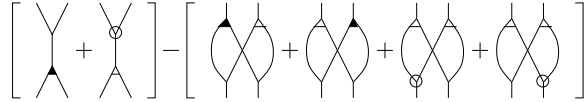}
}
\end{center}
\caption{The second $2$-differential $d^{2,2}$}
\label{diff2}
\end{figure}

The first and second differentials $d^{2,1}(\eta_1, \eta_2)$, $d^{2,2}(\eta_1, \eta_2)$
 are depicted in Fig.~\ref{diff1} and Fig.~\ref{diff2}, respectively.
Here we note that these diagrams agree with those for Hochschild
bialgebra cohomology in the sense that they are obtained by the
following process: (1) Consider the diagrams of the equality in
question (in this case the self-distributivity condition and the compatibility),
(2) Mark exactly one vertex of such a diagram, (3) Take a formal
sum of  such diagrams over all possible markings. In
Fig.~\ref{diff1}, the first two terms correspond to the LHS of
$q(q\otimes 1)=q(q \otimes q)\tau_2(1^2 \otimes \Delta)$, and one
of the two white triangular vertices is marked by a black vertex,
representing the $2$-cochain $\eta_1$, while the remaining white
vertex represents $q$. The negative four terms correspond to the
RHS, and the last term has a circle, representing $\eta_2$ while
unmarked ones in the rest represent $\Delta$.  The same procedure
for the compatibility gives rise to Fig.~\ref{diff2}.

\begin{figure}[htb]
\begin{center}
\mbox{
\epsfxsize=3.1in
\epsfbox{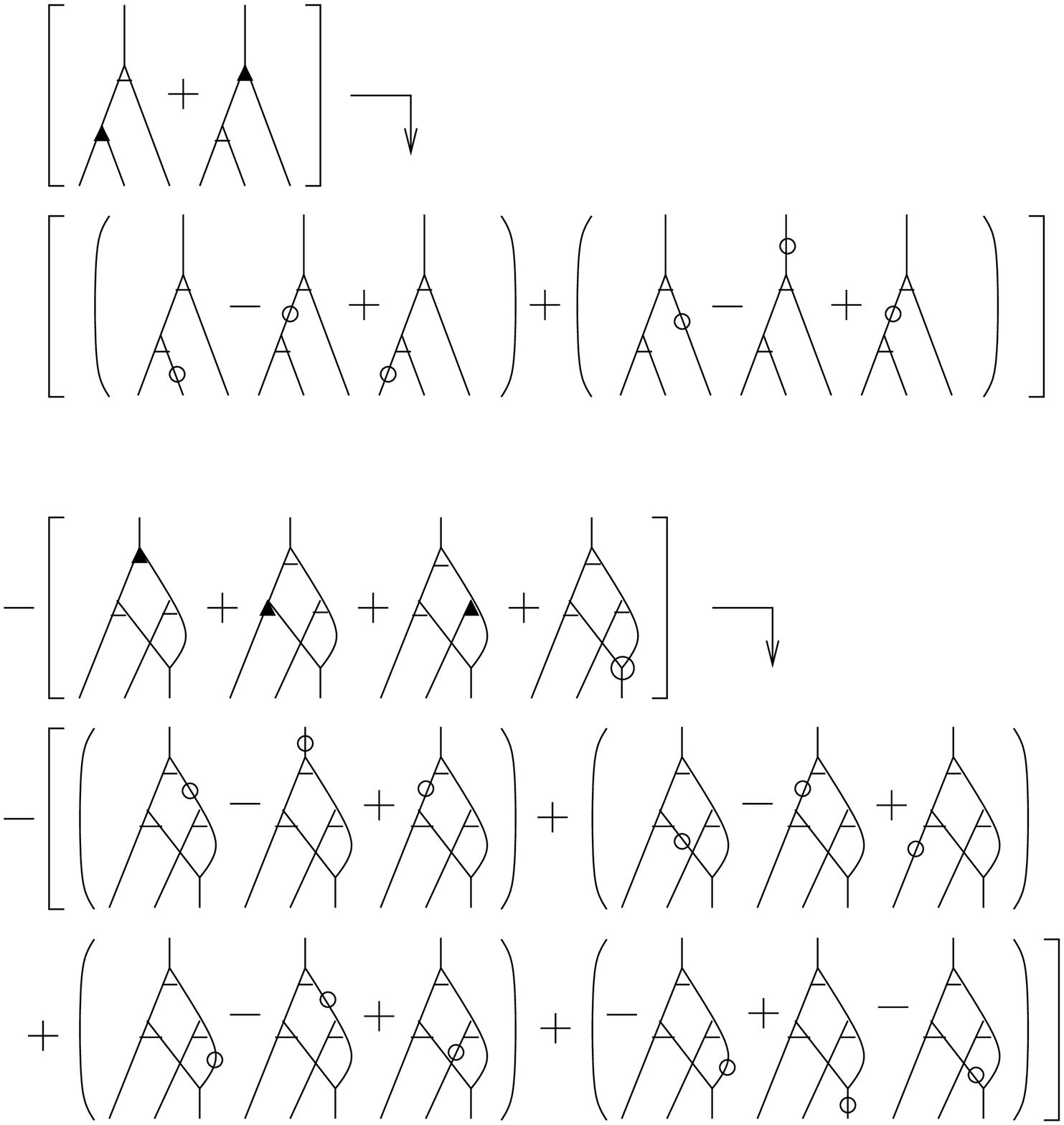}
}
\end{center}
\caption{$d^{2,1}(D_1(f))=0$}
\label{diff21c1}
\end{figure}

\begin{figure}[htb]
\begin{center}
\mbox{
\epsfxsize=3.5in
\epsfbox{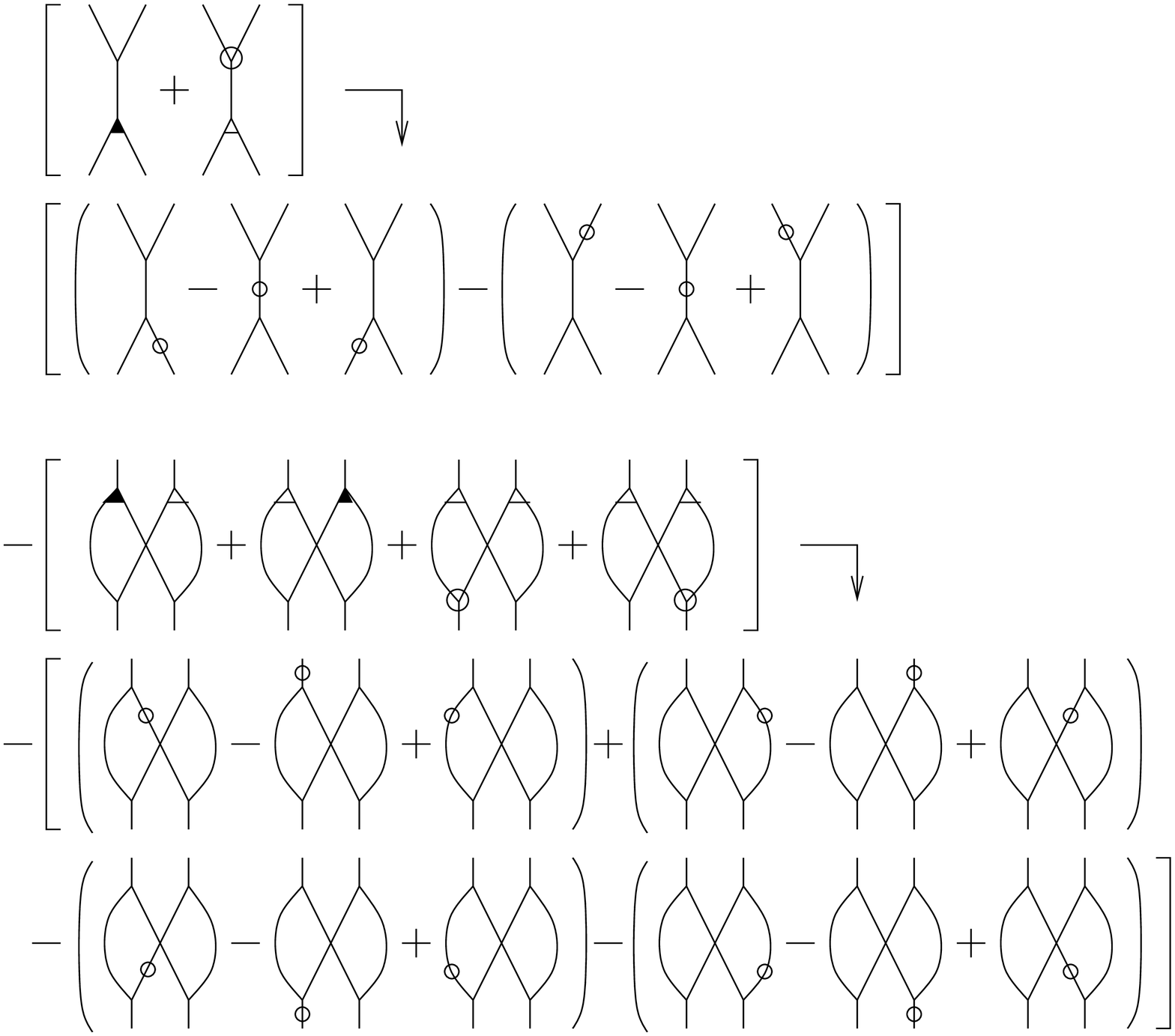}
}
\end{center}
\caption{$d^{2,2}(D_1(f))=0$}
\label{diff22c1}
\end{figure}

\begin{lemma}\label{d2d1lem}
For any $f \in C^{1}_{\rm sh}(X;X)$, we have $D_2 D_1(f)=0$.
 \end{lemma}
{\it Proof.\/} A proof is depicted in Fig.~\ref{diff21c1} and
Fig.~\ref{diff22c1}.
By assumption, $\eta_1=d^{1,1}(f)$ and $\eta_2=d^{1,2}(f)$.
Therefore, as in the case of Hochschild homology,
marked vertices representing $\eta_1$ and
$\eta_2$ are replaced by formal sum of three diagrams representing
$d^{1,1}(f)$ and $d^{1,2}(f)$, see Fig.~\ref{biHoch}.
The situation in which the first two terms are replaced by
three terms each is depicted in the top two lines of
Fig.~\ref{diff21c1}.

 A white circle on an edge represents $f$. The bottom three
lines show replacements for the remaining four negative terms.
Then the terms represented by identical graphs cancel directly. If
a white circle representing $f$
 appears near the boundary, then we use the
self-distributive axiom
 to relate this to another term. For example, the first term on the top left
 cancels with the third term on the bottom row since $f$ is on the second
 tensor factor at the bottom of each.

To facilitate the reader's understanding of the computation we
present the following sequences: $1,-2,3,4,-5,2$ and
$-6,5,-7,-8,7,-3,-9,6,-1,9,-4,8.$ Label the diagrams below the
arrows in Fig.~\ref{diff21c1} in order with these numbers. The
minus sign indicates the sign of the given term on the given side
of the equation, and the number indicates which diagrams cancel
which. A similar labelling can be accomplished in
Fig.~\ref{diff22c1}.
 $\Box$

\begin{figure}[htb]
\begin{center}
\mbox{
\epsfxsize=4in
\epsfbox{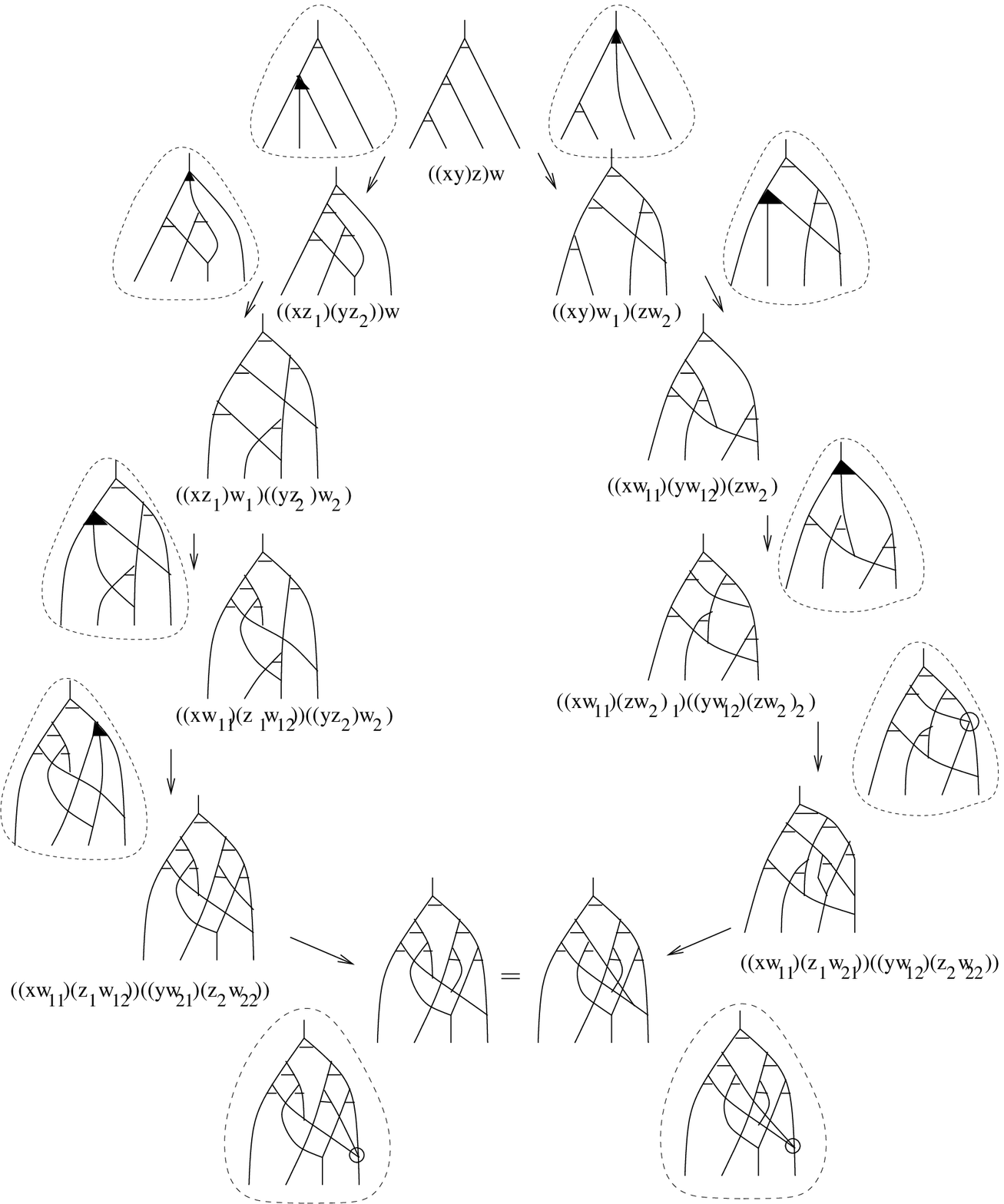}
}
\end{center}
\caption{First $3$-differential, $d^{3,1}$}
\label{D31loop}
\end{figure}

 We also note the following restricted version:

 \begin{lemma}\label{d11lem}
 Let $f \in C^{1}_{\rm sh}(X;X)=\Hom(X;X)$.
 If $d^{1,2} (f)=0 \in C^{2,2}_{\rm sh}(X;X)=\Hom(X;X^{\otimes 2})$,
 then $D_2(d^{1,1}(f), 0)=0$.
 \end{lemma}
 {\it Proof.\/} The conclusion is restated by the following condition:
 $d^{2,i}(d^{1,1}(f), 0)=0$ for $i=1,2$, since $d^{1,1}(f)$ is not
 in the domain of  the differential $d^{2,3}$.
Then one computes $d^{2,i}(\eta_1, 0)$ for $\eta_1=  d^{1,1}(f)$
either directly, or diagrammatically using Figs.~\ref{diff21c1},
and~\ref{diff22c1},
without
trivalent vertices that are encircled.
 $\Box$

\begin{figure}[htb]
\begin{center}
\mbox{
\epsfxsize=3.4in
\epsfbox{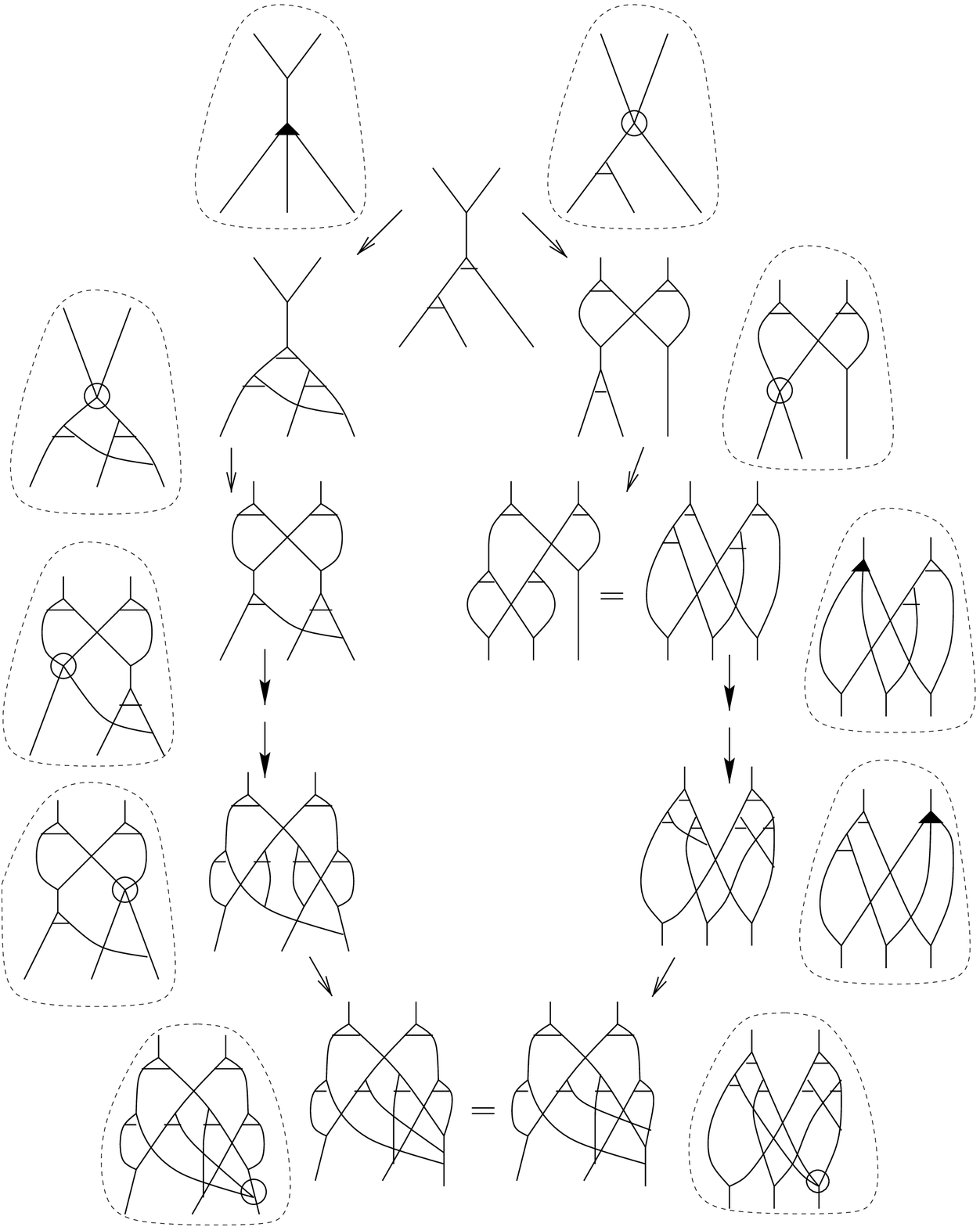}
}
\end{center}
\caption{Second $3$-differential, $d^{3,2}$}
\label{D32loop}
\end{figure}

\subsection{Third Differentials}

Throughout
this section, we consider only self-distributive
linear %% MS added
maps for cocommutative coalgebras with counits.
The maps $q$ need not be compatible with the counit, but there must be such a counit present.
 In this case,
$3$-differentials
$$d^{3, i} :  C^3_{\rm sh} (X;X) (=\oplus_{j=1}^3 \Hom(X^{\otimes (4-j)}, X^{\otimes j} ))
 \rightarrow C^{4, i}_{\rm sh}(X;X) (= \Hom(X^{\otimes (n+2-i)}, X^{\otimes i} ) )$$
are defined below
for $i=1,2,3$, and for $i=4$ it is defined by
 the same map as  the differential for $\Delta$
for co-Hochschild cohomology (the pentagon
identity for the comultiplication). Let $\xi_j \in \Hom(X^{\otimes
(4-j)} , X^{\otimes j})\subset C^3_{\rm sh} (X;X) $, $j=1,2,3$.

\begin{figure}[htb]
\begin{center}
\mbox{
\epsfxsize=2.6in
\epsfbox{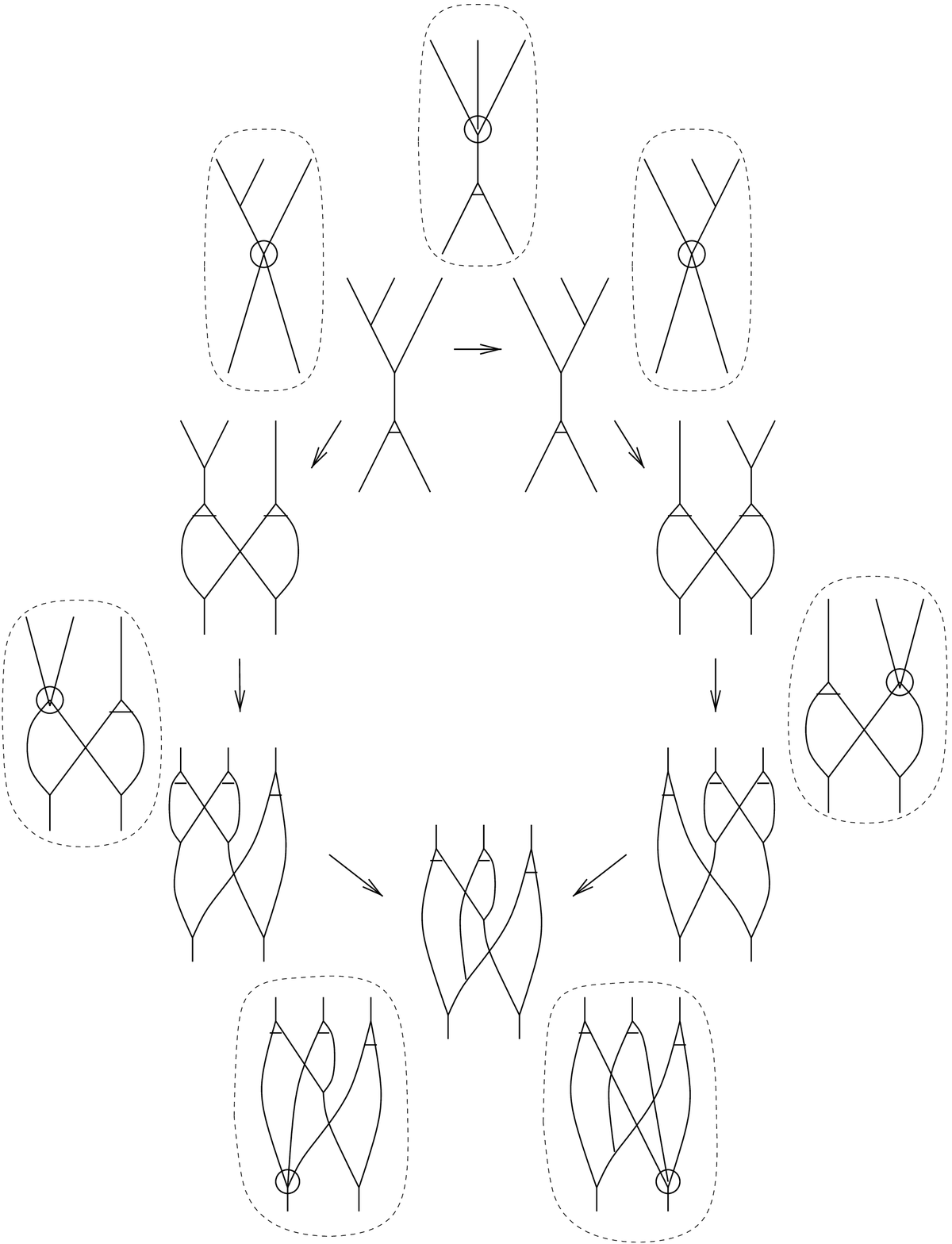}
}
\end{center}
\caption{Third $3$-differential, $d^{3,3}$}
\label{D33loop}
\end{figure}

These differentials are defined by direct analogues with
Hochschild differentials in diagrammatics, and we will justify our
definition in two more
ways: (1)
$2$-cochains vanish under these maps, (2) $3$-cocycles of quandle
and Lie algebra cohomology are realized in these formulas as
discussed in the next section.

First we explain the diagrammatics. Recall that $3$-cocycle
conditions in Hochschild cohomology correspond to two different
sequences of relations applied to graphs that change one graph to
another. At the top of Fig.~\ref{D31loop},
a graph
representing $q(q \otimes 1)(q \otimes 1^2)$ is depicted. There
are two ways to
apply  sequences of self-distributivity
 to this map to get the map
represented by the bottom graphs. A $3$-cochain $\xi_1$
represented by a
black triangular
vertex with three bottom edges and a single top edge corresponds
to applying  the self-distributivity relation to change a graph to another,
and corresponds to where the self-distributivity relation
was applied. The two different sequences are shown at the left and
right of the figure. These sequences give rise to the
LHS and RHS of $d^{3,1}(\xi_1, \xi_2, \xi_3)$.
Similar graphs are obtained as shown in
Figs.~\ref{D32loop} and \ref{D33loop}.

%%% NOTE THE $\xi_3$ TERMS WERE ADDED:
The differentials thus obtained are:
\begin{eqnarray*}
 d^{3,1}(  \xi_1, \xi_2, \xi_3)
& = &
 [ q (\xi_1 \otimes 1) +
 \xi_1 (q \otimes
q \otimes 1)(1 \otimes \tau
\otimes 1^2) (1^2  \otimes \Delta
\otimes 1
  )\\
& &  + \ \
q(\xi_1 \otimes q)
(1^2 \otimes \tau \otimes 1)
 (1^3 \otimes \Delta) (1^2 \otimes q \otimes 1) \\
& & (1 \otimes \tau
 \otimes 1^2)
 (1^2  \otimes \Delta\otimes 1)\\
  & &+ \ \ q (q \otimes \xi_1)(q\otimes q \otimes 1^3)(1 \otimes \tau \otimes 1^4) \\
 & &(1^2 \otimes \Delta \otimes 1^3)(1^2 \otimes \tau \otimes 1^2)(1 \otimes \tau \otimes \tau \otimes 1)(1^2 \otimes \Delta \otimes \Delta)\\
 & &+ \  q(q\otimes q)(q \otimes q \otimes q \otimes q)(1 \otimes \tau \otimes 1^2 \otimes \tau \otimes 1)(1^2 \otimes \Delta \otimes 1^3) \\
& & (1^2 \otimes \tau \otimes 1^2)(1 \otimes \tau \otimes \tau \otimes 1^2)(1^2 \otimes \Delta \otimes \xi_3)] \\
& &- \  [
\xi_1
(q \otimes 1^2)
+q (\xi_1 \otimes
q)
(1^2 \otimes \tau \otimes 1)
(1^3 \otimes \Delta) \\
  & & + \ \
\xi_1 (q \otimes q \otimes q)
(1 \otimes \tau \otimes 1^3)
 (1^2 \otimes \Delta \otimes 1^2)
(1^2 \otimes \tau \otimes 1)
(1^3 \otimes \Delta) \\
 & &
 + \ \ q
(q \otimes q)
(1 \otimes \tau \otimes 1)
(q \otimes q \otimes 1^2)
(1^2 \otimes \Delta \otimes 1^2)
 (1^3 \otimes \xi_2 )
(1^2 \otimes \tau \otimes 1)
(1^3 \otimes \Delta) \\
& & + \  q(q\otimes q)(q \otimes q \otimes q \otimes q)(1 \otimes \tau \otimes \tau \otimes 1^3)(1^4 \otimes \tau \otimes \tau \otimes 1^2) \\
& & (1 \otimes \tau \otimes \Delta \otimes \tau \otimes 1)(1^2 \otimes \Delta \otimes \xi_3) %%%here
]
\\
 d^{3,2}(  \xi_1, \xi_2, \xi_3)
& = &
[ \Delta \xi_1 + \xi_2 (q \otimes q)
(1 \otimes \tau \otimes 1)(1^2 \otimes \Delta)
 \\
 & & + \ \  (q \otimes q)
(1 \otimes \tau \otimes 1)
(\xi_2 \otimes \Delta
)( 1^2 \otimes
q)
(1 \otimes \tau \otimes 1)
 (1^2 \otimes \Delta) \\
 & &+ \   (q \otimes q)(1 \otimes \tau \otimes 1)( \Delta \otimes 1^2)
)(
 q \otimes \xi_2)
(1 \otimes \tau \otimes 1)
(1^2 \otimes \Delta) \\
& & + \ (q \otimes q)(q \otimes \tau \otimes q)(1^2 \otimes q \otimes q \otimes 1^2)(1 \otimes \tau \otimes \tau \otimes \tau \otimes 1) \\
& & (\Delta \otimes 1^2 \otimes \tau \otimes \Delta)(1 \otimes \tau \otimes 1^3)(1 \otimes \Delta \otimes \xi_3)   %%%%added last term
 ] \\
 &-& [
 \xi_2(q \otimes 1)
 + (q \otimes q)
(1 \otimes \tau \otimes 1)
 (\xi_2 \otimes \Delta ) \\
 & & +  \ \ (\xi_1 \otimes q) (1^3 \otimes q \otimes 1)
(1^2 \otimes \tau \otimes 1^2)(1 \otimes \tau \otimes \tau \otimes 1)
 (\Delta \otimes \Delta \otimes \Delta) \\
 & & +  \ \  (q \otimes \xi_1) (q \otimes 1^4)
(1^2 \otimes \tau \otimes 1^2)(1 \otimes \tau \otimes \tau \otimes 1)
 (\Delta \otimes \Delta \otimes \Delta)\\
& & + \ \ (q \otimes q)(q \otimes q \otimes q \otimes q)(1  \otimes \tau \otimes \tau \otimes \tau \otimes 1) (1^2 \otimes \tau \otimes \tau \otimes \Delta ) \\ %%%%added
& & (1 \otimes \tau \otimes \tau \otimes 1^2)(\Delta \otimes \Delta \otimes \xi_3) %%added
 ]
 \\
d^{3,3}(  \xi_1, \xi_2, \xi_3)
& = &
[(\Delta \otimes 1) \xi_2
 + (\xi_2 \otimes q)
(1 \otimes \tau \otimes 1)
(\Delta \otimes \Delta) \\
& & + \ \ (q \otimes q \otimes q)(1 \otimes \tau \otimes 1^3)(1^2 \otimes \Delta \otimes 1^2)(1^2 \otimes \tau \otimes 1^2)(\xi_3 \otimes \Delta) ]\\
&-&
[ \xi_3 q +(1  \otimes \Delta) \xi_2
+ ( q \otimes  \xi_2)
%%%\tau_2
(1 \otimes \tau \otimes 1) %%%%
(\Delta \otimes \Delta) \\
& & + \ \
(q \otimes q \otimes q)(1 \otimes \tau \otimes \tau \otimes 1)(1^2 \otimes \tau \otimes 1^2)(1 \otimes \Delta \otimes 1^3)( \Delta \otimes \xi_3)
]
 \end{eqnarray*}

\begin{figure}
\begin{center}
\mbox{
\epsfxsize=3.6in
\epsfbox{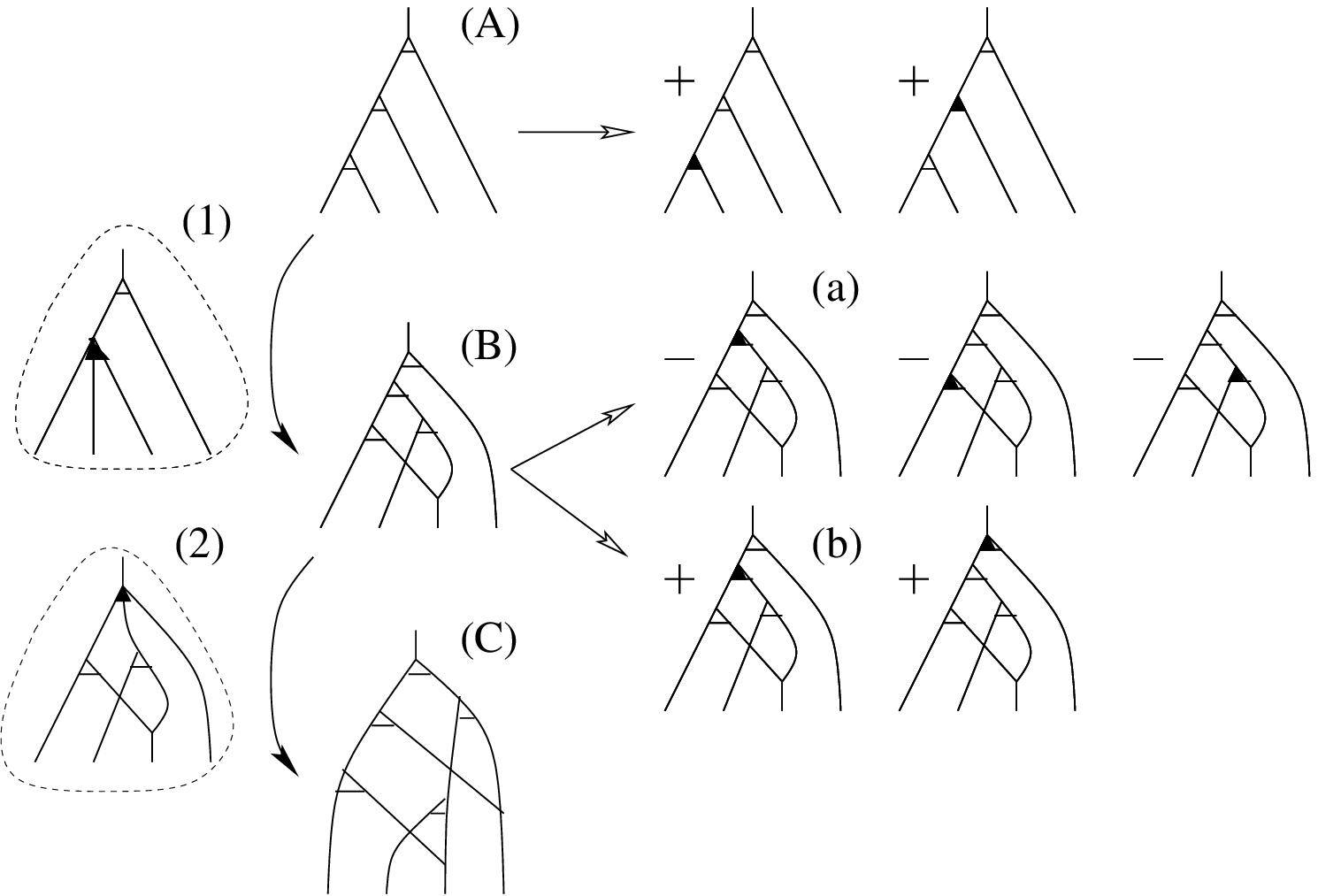}
}
\end{center}
\caption{A strategy for a proof}
\label{d3proof}
\end{figure}

%%%D_3 added here
The third differential is defined as $D^3= d^{3,1}+d^{3,2}+d^{3,3}+d^{3,4}.$
%%%%%
\begin{lemma} \label{d321lem}
Let $(\eta_1, 0) \in \Hom(X^{\otimes 2}, X) \subset C^2_{\rm sh}(X;X)$
$($so that $\eta_2=0)$.
Then $d^{3,1} d^{2,1}(\eta_1, 0)=0$.
\end{lemma}
{\it Proof.\/} This is proved by calculations that seem
complicated without diagrammatics. We sketch our computational
method. For $\xi_1 \in \Hom(X^{\otimes 3}, X) \subset C^3_{\rm
sh}(X;X)$, the first two terms of $d^{3,1}(\xi_1)$ are $q(\xi_1
\otimes 1)$ and $\xi_1 (q \otimes q \otimes 1)(\tau_2)(1^2 \otimes
\Delta \otimes 1)$, that are diagrammatically represented by left
of Fig.~\ref{d3proof}, (1) and (2), respectively. The black
triangular four-valent vertex represents $\xi_1$. On the other
hand, the first term $q(\xi_1 \otimes 1)$ corresponds to the
change of the diagrams represented in (A) and (B). Such a change
of diagrams corresponds to $d^{(2)}_1(\eta_1)$ as depicted in
Fig.~\ref{diff1}. Therefore the first terms of  $d^{3,1} d^{2,1}
(\eta_1, 0)$ are $q(\xi_1 \otimes 1)=q(  d^{2,1} (\eta_1, 0)
\otimes 1)$ consisting of five terms represented by the diagrams
on the right top two rows in Fig.~\ref{d3proof}. The third row
consists of the positive terms of the second term (2), $d^{2,1}
(\eta_1, 0) (q \otimes q \otimes 1)(\tau_2)(1^2 \otimes \Delta
\otimes 1)$. Thus to prove this lemma, we write out all terms and
check that they cancel. For example, the terms on the right of
Fig.~\ref{d3proof} labelled with (a) and (b) cancel.

%% MS bringing these figs to Appendix.
The essential steps in the proofs are found in Figs.~\ref{D31L} and \ref{D31R}
in the appendix. %% MS
 The
first rows of Fig.~\ref{D31L} coincide with those of
Fig.~\ref{d3proof}. The remaining left columns indicate the
different diagrams that are obtained by replacing the four-valent
black vertices by the two sides of
the self-distributive law. %% MS relation.
The right-hand entries are the expansions of the terms in the next differential.
The terms are numbered and those in Fig.~\ref{D31L} and Fig.~\ref{D31R} cancel.

It is somewhat difficult to see the cancellation of the terms
labelled $7,10,11,12,13,$  and $14.$  The terms labelled $15$
coincide by applications of coassociativity and cocommutativity.
The identity between these terms becomes obvious after one works
through the preceding terms. The proofs that the diagrams
represent the same linear maps are provided in the appendix.
$\Box$

Similarly we obtain:
\begin{lemma} \label{d322lem}
$d^{3, 2} d^{2, i} (\eta_1, 0)=0$ for $i=2,3,4$,
and $d^{3,3} d^{2,i}(\eta_1, 0)=0$ for $i=2,3,4$.
\end{lemma}

%% MS add
Diagrams that represent the proofs are included in the appendix. %% MS

\subsection{Cohomology Groups}

Now we use these differentials to define cohomology groups for
%%%%  MS
%shelves in $\Coalg.$
self-distributive linear maps for objects in $\Coalg$.
%% Let $(X, \Delta, q)$ be an object in  $\Coalg$ that satisfies self-distributivity.
Let $(X, \Delta)$ be an object in $\Coalg$, and
$q: X \otimes X \rightarrow X$ be a self-distributive linear map. %% end MS
Then Lemma~\ref{d2d1lem} implies:
\begin{corollary}
$0 \rightarrow C^{1}_{\rm sh}(X;X)
\stackrel{D_1}{\rightarrow} C^{2}_{\rm sh}(X;X) $
is a chain complex.
\end{corollary}
This enables us to define the following cohomology related groups:
\begin{definition}{\rm
The $1$-cocycle and cohomology groups are defined by:
$$H^{1,i}_{\rm sh}(X;X) =
Z^{1,i}_{\rm sh}(X;X) = \{ f \in C^{1,i}_{\rm sh}(X;X) \ | \ d^{1,i}(f)=0 \} $$
for $i=1,2$, and
$$H^{1}_{\rm sh}(X;X) = Z^{1,1}_{\rm sh}(X;X) \oplus Z^{1,2}_{\rm sh}(X;X) .$$
} \end{definition}

Since the $2$-cocycle conditions were formulated directly from a
deformation theory formulation, we have the following:

\begin{proposition}
Let $X_t=X \otimes k[[t]]$ and suppose we have partial
deformations $\bar{q} = q + \cdots + t^n q_n$ and $\bar{\Delta}
=\Delta + \cdots + t^n \Delta_n$ satisfying the above three
conditions mod $t^{n+1}$, so that they define a
self-distributive map
in $\Coalg$ mod $t^{n+1}$. Then there exist $q_{n+1}: X \otimes X
\rightarrow X$ and $\Delta_{n+1}: X \rightarrow X \otimes X $ such
that $\bar{q}+ t^{n+1}q_{n+1}$ and $\bar{\Delta} + t^{n+1}
\Delta_{n+1}$ satisfy the three conditions mod $t^{n+2}$, so that
they define a
self-distributive
linear %% MS
map
%% MS delete: %% in $\Coalg$
mod $t^{n+2}$, if and only if $(q_{n+1},
\Delta_{n+1})$ satisfy the $2$-cocycle condition: $D_2(q_{n+1},
\Delta_{n+1})=0$.
\end{proposition}

For $3$-cocycles, we recall that $(X, \Delta, q)$
consists of an object $(X, \Delta)$ in $\CoComCoalg$,
with a self-distributive linear map $q$.
 Let
$d_1^{n,i}= d^{n,i} | (C^{n, 1}_{\rm sh} (X;X) )$ be the
restriction of $d^{n,i}$  to $C^{n, 1}_{\rm sh} (X;X)=
\Hom(X^{\otimes (n)}, X)$, and $D'_{(1)}=d^{1,1}$,
$D'_{n}=\sum_{i=1}^{n+1} d_1^{n,i}$ for $n=2,3$ and $i=1,2,3$.
Then consider the sequence
$$
{\cal C} : \quad 0 \rightarrow Z^{1,2}_{\rm sh}(X;X)
\stackrel{D'_1}{\rightarrow} C^{2,1}_{\rm sh}(X;X)
\stackrel{D'_2}{\rightarrow} C^{3,1}_{\rm sh}(X;X)
\stackrel{D'_3}{\rightarrow} C^{4,1}_{\rm sh}(X;X) .
$$
Then Lemmas~\ref{d11lem}, \ref{d321lem}, and \ref{d322lem} are
summarized as:
\begin{theorem}
%% MS:
%Let $(X, \Delta, q)$  be an object in
Let $(X, \Delta)$  be an object in
$\CoComCoalg$ and %% $q$
$q: X \otimes X \rightarrow X $ be %% MS
a self-distributive
linear %% MS
map.  Then ${\cal C} $
is a chain complex.
\end{theorem}

This enables us to define:
\begin{definition}{\rm
The % MS delete (never used) %  {\it cocommutative}
$1$-cocycle and cohomology group
are % MS % is
defined as:
$$H'{}^{1,1}_{\rm sh}(X;X) =
Z'{}^{1,1}_{\rm sh}(X;X) = \{ f \in Z^{1,2}_{\rm sh}(X;X) \ | \
d^{1,1}(f)=0 \}, $$ and the
% MS delete (never used) % {\it cocommutative}
$2$- and
$3$-coboundary, cocycle, and cohomology groups are defined as:
\begin{eqnarray*}
B^{j,1}(X;X) &=& {\rm Image}(D'_{j-1}) , \\
Z^{j,1}(X;X) &=& {\rm Ker}(D'_j) , \\
H^{j,1}(X;X) &=& Z^{j,1}(X;X) / B^{j,1}(X;X)
\end{eqnarray*}
for $j=2,3$.
} \end{definition}

%WATCH OUT HERE:
The cocycles in these theories are called {\it shelf cocycles}. The name is a bit of a notational compromise.
They should be called ``cocycles for self-distributive linear maps
for objects %%
in the category of cocommutative coalgebras with counit,''
which would  inevitably get shortened to cocococo-cycles.
There are two points here. First, the analogy
 ``quandle is to rack as rack is to shelf" {\it does not extend}
 to the terminology for shelf-cohomology. More importantly,
we do not require $q$ to be compatible with counit
in defining cohomology theories, yet we call them shelf cocycles for short. %% end

\section{Relations to Other Cohomology Theories}\label{relsec}

In this section we examine relations of these cocycles to those in
other cohomology theories, specifically the original  quandle
cohomology theories~\cite{CJKLS} and Lie algebra cohomology.
%% Mention this as needed:
%Throughout this section we assume that $\Delta $ is not deformed, i.e., $\eta_2=0$.

\subsection{Quandle Cohomology}

In this section we present  procedures
that produce  shelf $2$- and
$3$-cocycles from quandle $2$- and $3$-cocycles, respectively, and
show that non-triviality is inherited by these processes.

First we briefly review the definition of  quandle $2$- and
$3$-cocycles. A quandle {\it $2$-cocycle} is a linear function
$\phi$ defined on the free abelian group generated by pairs of
elements $(x,y)$ taken from a quandle $X$ such that
$$\phi(x,y) - \phi(x,z) + \phi(x\lt y, z) - \phi(x\lt z, y \lt z)=0, \quad \forall x,y,z \in X$$
and $\phi(x,x)=0$ for all $x \in X$. The function $\phi$ takes
values in some fixed abelian group $A$. Similarly a {\it
$3$-cocycle} is a function $\theta$ with the properties that
$$\theta(x,y,z) + \theta(x \lt z, y \lt z, w) + \theta ( x,z, w) = \theta (x\lt y,z,w) +\theta(x,y,w) + \theta( x \lt w, y \lt w, z \lt w),$$
and
$$\theta(x,x,y)=\theta(x,y,y)=0$$ for all $x,y,z,w\in X.$
Quandle cohomology groups $H^n_{\rm Q}(X;A)$ were defined based on these conditions,
see \cite{CJKLS,CKamS} for details.

These cocycles were used to develop invariants of classical knots
and knotted surfaces. We summarize the construction as follows.
Given a quandle homomorphism from the fundamental quandle  of a
codimension $2$ embedding to the finite quandle $X$, and given a
cocycle ($\phi$ or $\theta$), we evaluate the cocycle at the
incoming quandle elements near each $0$-dimensional multiple point
(crossing and triple point, respectively), in the projection of
the knot or knotted surface. These values are added together in
the abelian group $A$, and the collection of the results are
formally collected together as a multiset over all homomorphisms.
The cocycle invariants are fairly powerful in
 determining properties of knots and knotted surfaces.
Generalizations have been discovered
\cite{AG,CES,CESbiq}.

Recall that $W= k \oplus kX$ ($V=kX$) is the direct sum of the field $k$ and the vector space whose basis is comprised of the elements in $X$,
and the self-distributive map $q$ defined on $V$ was extended to $W$.

\begin{theorem}
\label{q2prop}
For a quandle $2$-cocycle $\phi$
with the coefficient group $A=k$,
define $\hat{\phi}: W \otimes W \rightarrow W$ by
linearly extending $\hat{\phi} ( x \otimes y)=\phi(x,y)$,
$\hat{\phi} (1 \otimes x)=1$, and $\hat{\phi}(x \otimes 1)=\hat{\phi}(1 \otimes 1)=0$
for $x, y \in X$.
Then $\hat{\phi}$ 
%%%is a shelf  $2$-cocycle: 
satisfies %%%%%
$d^{2,1} (\hat{\phi}, 0) = 0 $.
\end{theorem}
{\it Proof.\/}
We write expressions such as $(a+\sum_x a_x x)$ in the more
compact form $(a + Ax)$. Then

\begin{eqnarray*}
\lefteqn{ \ (a+\sum_x a_x x)\otimes (b+\sum_y a_y y) \otimes (c+\sum_z a_z z)\  }  \\
&= &
 abc (1 \otimes 1 \otimes 1) + Abc (x \otimes 1 \otimes 1) + aBc (1 \otimes y \otimes 1) + ABc (x \otimes y \otimes 1)   \\
&+& abC (1 \otimes 1 \otimes z) + AbC (x \otimes 1 \otimes z) +
aBC (1 \otimes y \otimes z) + ABC (x \otimes y \otimes z)
\end{eqnarray*}

In  order to compute $d^{2,1}(\hat{\phi},0)$ on expressions such
as the one above, we must compute
 it on the eight tensor
products $(1 \otimes 1 \otimes 1)$ through $(x \otimes y \otimes
z)$. These calculations are summarized in the table below
(Juxtaposition or commas are used in place of $\otimes$ for
typesetting purposes.):

\vspace*{1cm}

\noindent
\begin{tabular}{|c||c|c|c|c|c|c|c|c||} \hline \hline
\negthinspace \negthinspace \negthinspace \mbox{
\epsfxsize=.5cm
\epsfbox{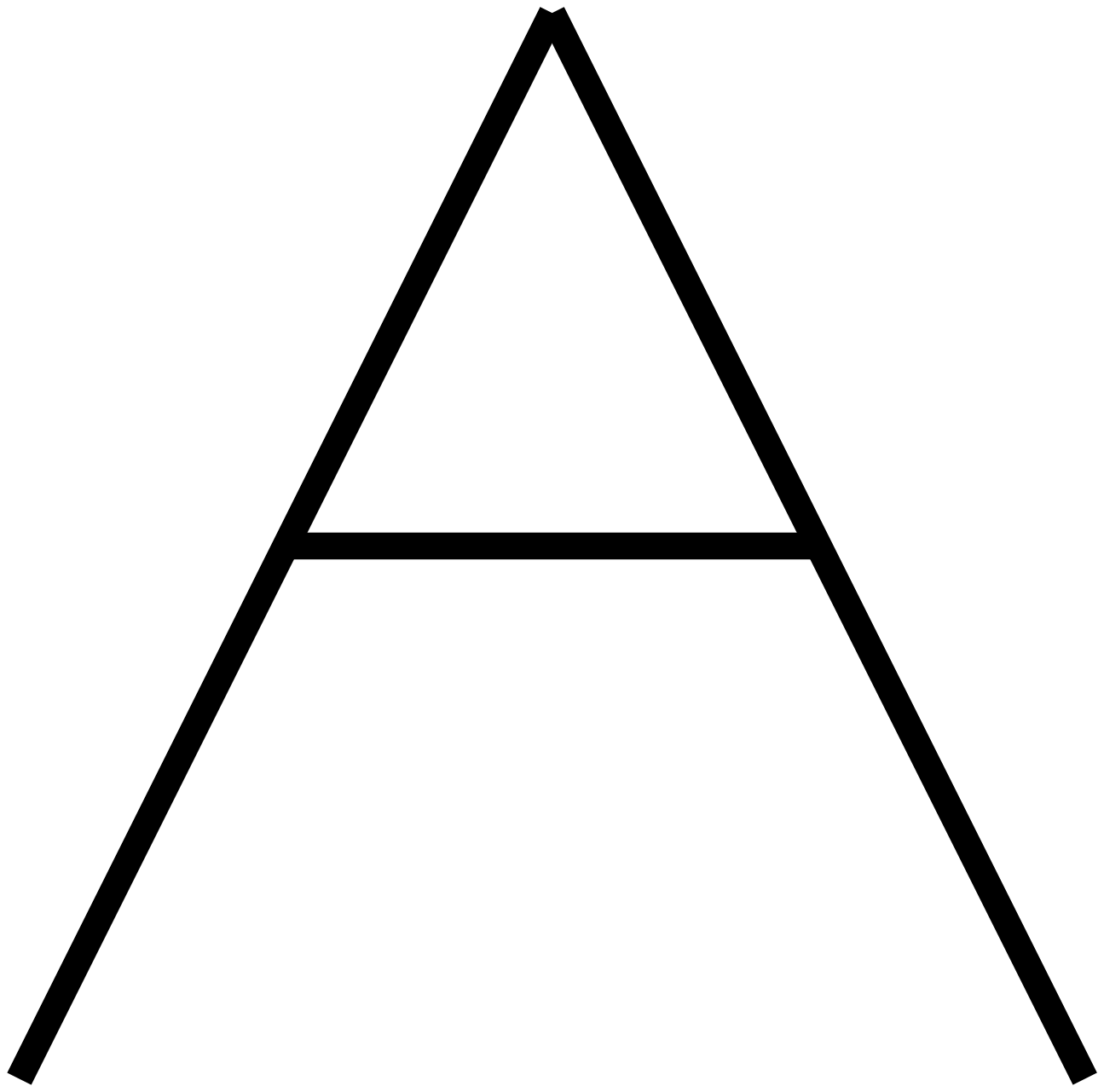}
}  \negthinspace \negthinspace \negthinspace \negthinspace & $0$ & $0$& $0$ & $0$ & $0$  &$0$  & $1$  & $ \phi(x,y)$ \\  \hline
\negthinspace \negthinspace \negthinspace \mbox{
\epsfxsize=.7cm
\epsfbox{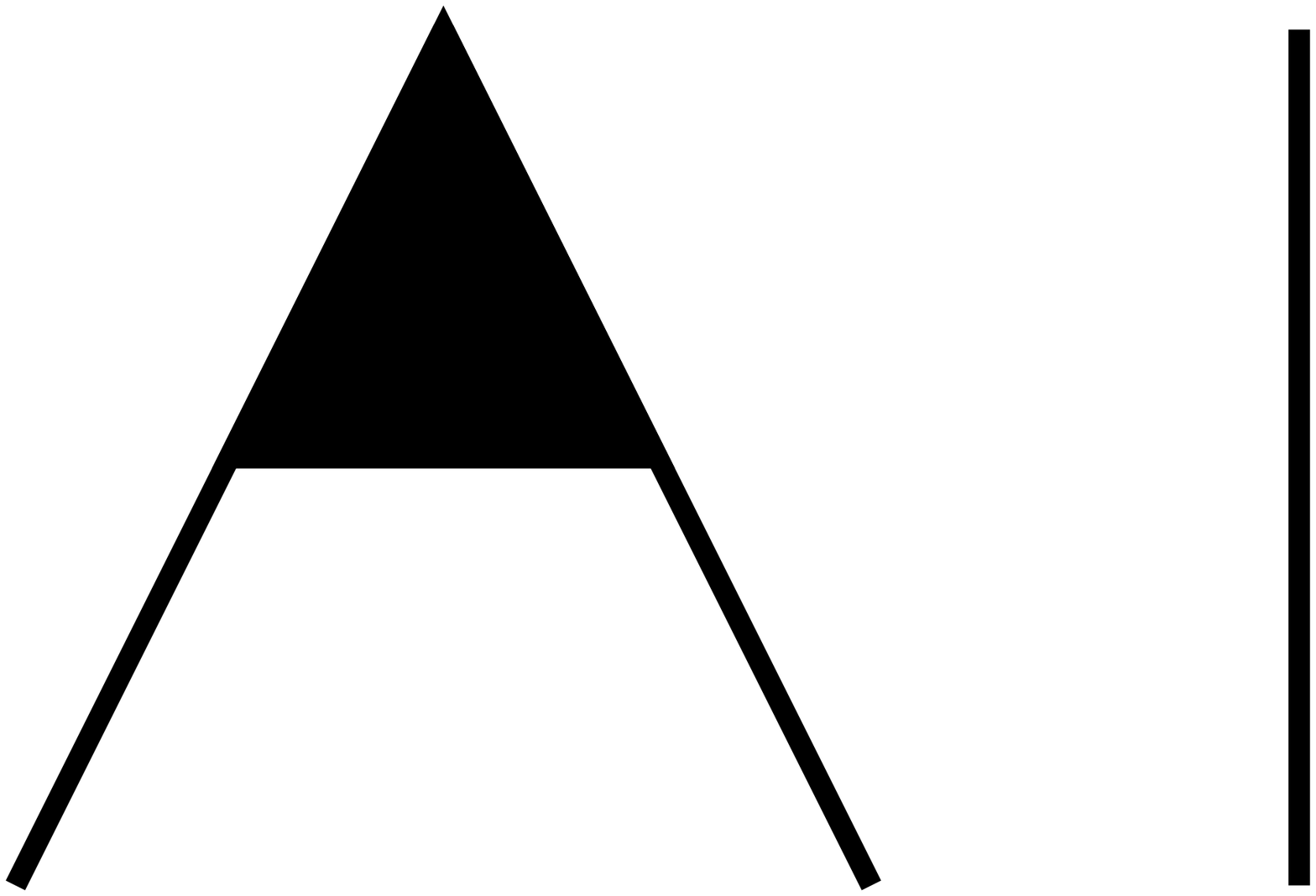}
} \negthinspace
\negthinspace \negthinspace \negthinspace & $0 \otimes 1$ & $0 \otimes 1$ & $1 \otimes 1$ & $\phi(x,y)\otimes 1$ & $0 \otimes z$ & $0\otimes z$ &$1 \otimes z$ & $\phi(x,y) \otimes z$ \\ \hline \hline
\negthinspace \negthinspace \negthinspace \mbox{
\epsfxsize=.5cm
\epsfbox{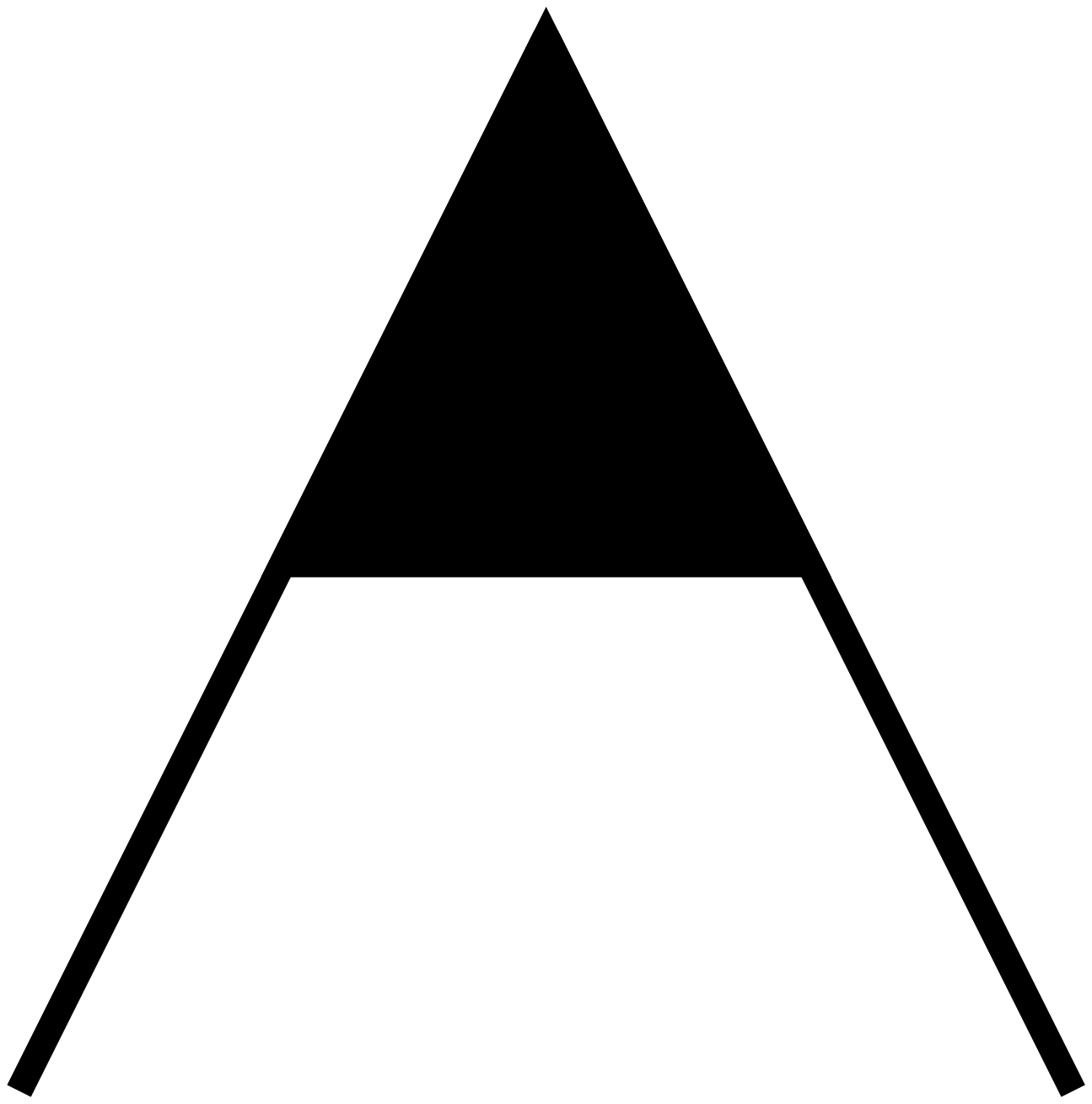}
}\negthinspace \negthinspace \negthinspace \negthinspace
    & $0$ & $0$& $0$ & $0$ & $0$  &$0$  & $1$  & $ \phi(x\lt y,z)$ \\  \hline
\negthinspace \negthinspace \negthinspace \mbox{
\epsfxsize=.7cm
\epsfbox{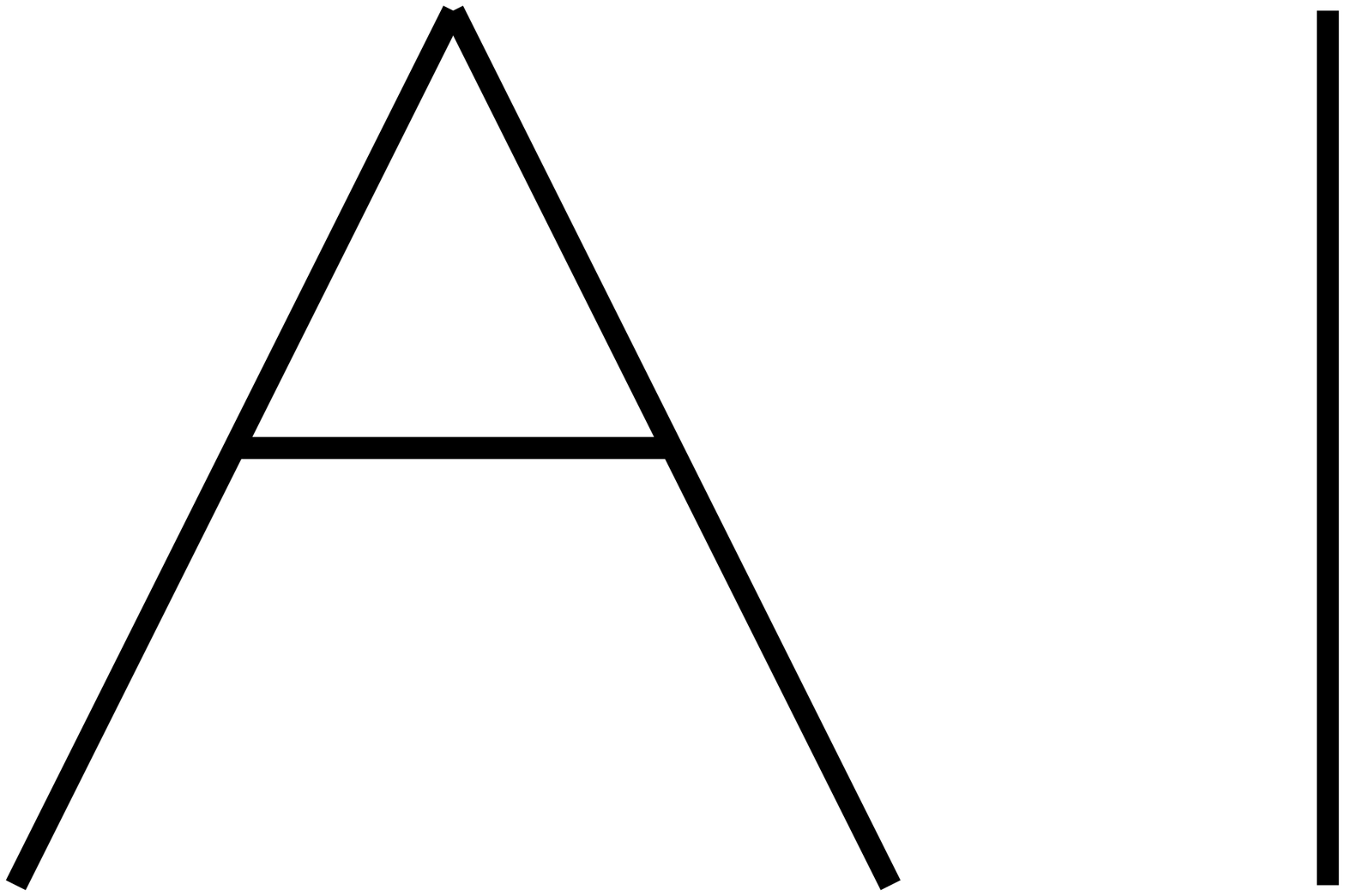}
} \negthinspace \negthinspace \negthinspace \negthinspace
 & $0 \otimes 1$ & $0 \otimes 1$ & $1 \otimes 1$ & $(x\lt y) \otimes 1$ & $0 \otimes z$ & $0\otimes z$ &$1 \otimes z$ & $(x \lt y) \otimes z$ \\ \hline \hline
\negthinspace \negthinspace \negthinspace \mbox{
\epsfxsize=.5cm
\epsfbox{q.eps}}
\negthinspace \negthinspace  \negthinspace \negthinspace &$0$ & $0$ & $0$ &$0$ & $0$ & $0$ & $1$ & $\phi(x,z)$ \\ \hline
\negthinspace \negthinspace \negthinspace \mbox{
\epsfxsize=.7cm
\epsfbox{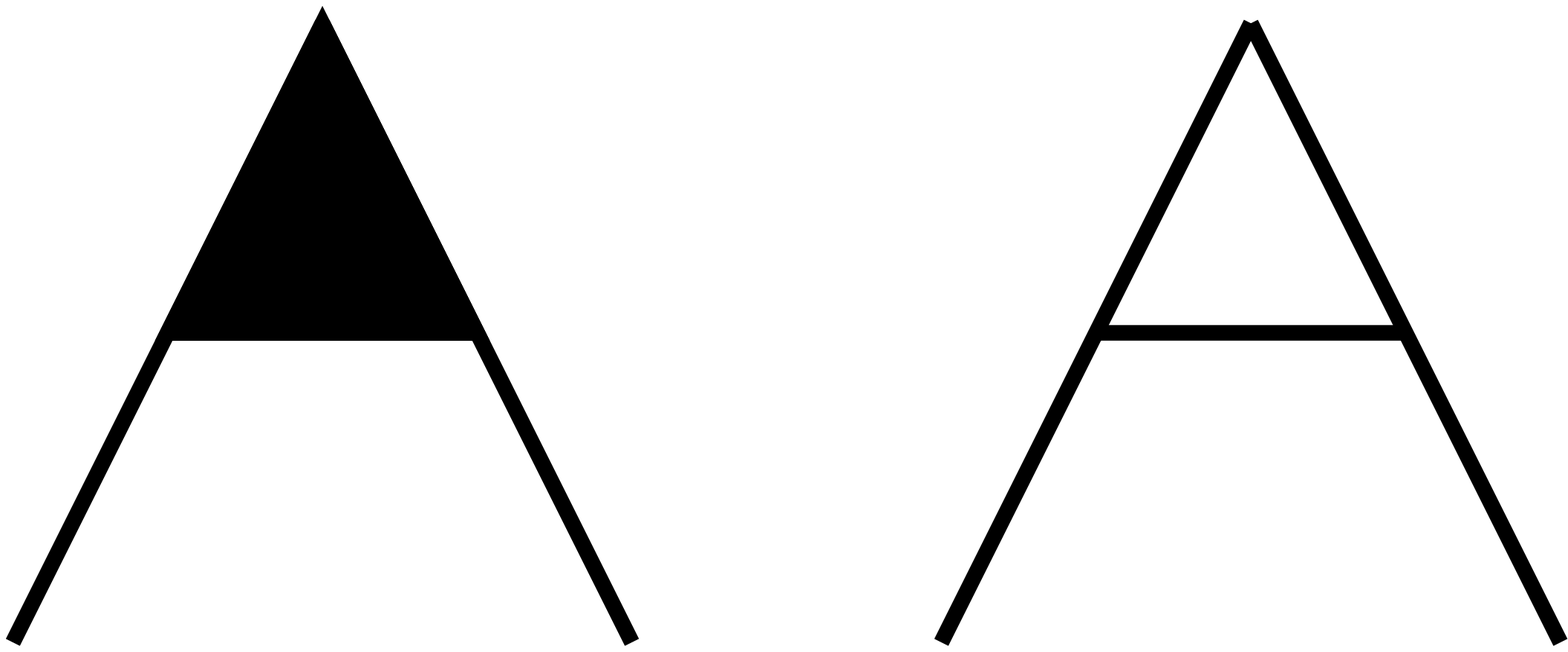}}\negthinspace \negthinspace \negthinspace \negthinspace
& $0 , 0$ & $0 , 0$ & $0 , 0 $ & $0 , 0$ & $1 ,1$ & $\phi(x,z),1$ &$1 , y\lt z$ & $\phi(x,z), y\lt z$ \\ \hline \hline
\negthinspace \negthinspace \negthinspace \mbox{
\epsfxsize=.5cm
\epsfbox{phi.eps}}
\negthinspace \negthinspace  \negthinspace \negthinspace &$0$ & $0$ & $0$ &$0$ & $0$ & $0$ & $1$ & $\phi(x\lt z,y \lt z)$ \\ \hline
\negthinspace \negthinspace \negthinspace \mbox{
\epsfxsize=.7cm
\epsfbox{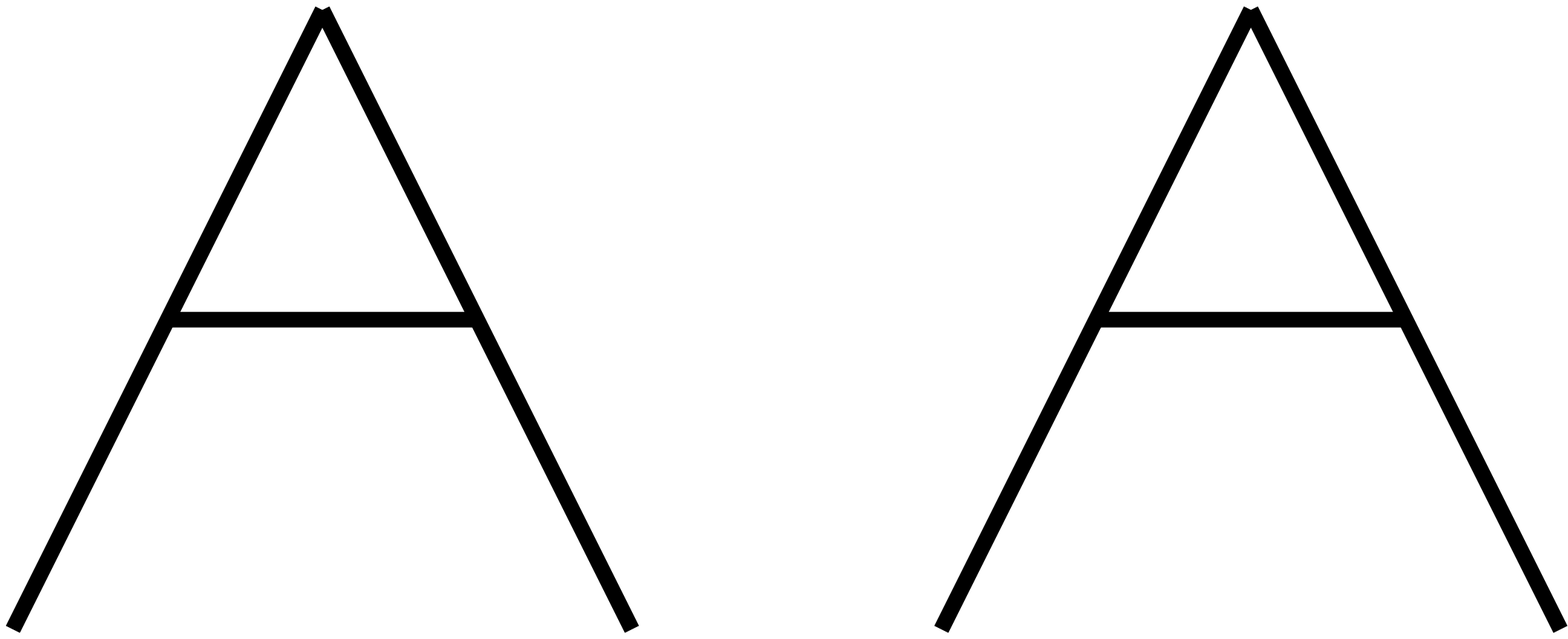}}\negthinspace \negthinspace \negthinspace \negthinspace
& $0 , 0$ & $0 , 0$ & $0 , 0 $ & $0 , 0$ & $1 ,1$ & $x \lt z,1$ &$1 , y\lt z$ & $x\lt z, y\lt z$ \\ \hline \hline
\negthinspace \negthinspace \negthinspace \mbox{
\epsfxsize=.5cm
\epsfbox{q.eps}}
\negthinspace \negthinspace  \negthinspace \negthinspace &$0$ & $0$ & $0$ &$0$ & $0$ & $0$ & $0$ & $0$ \\ \hline
\negthinspace \negthinspace \negthinspace \mbox{
\epsfxsize=.7cm
\epsfbox{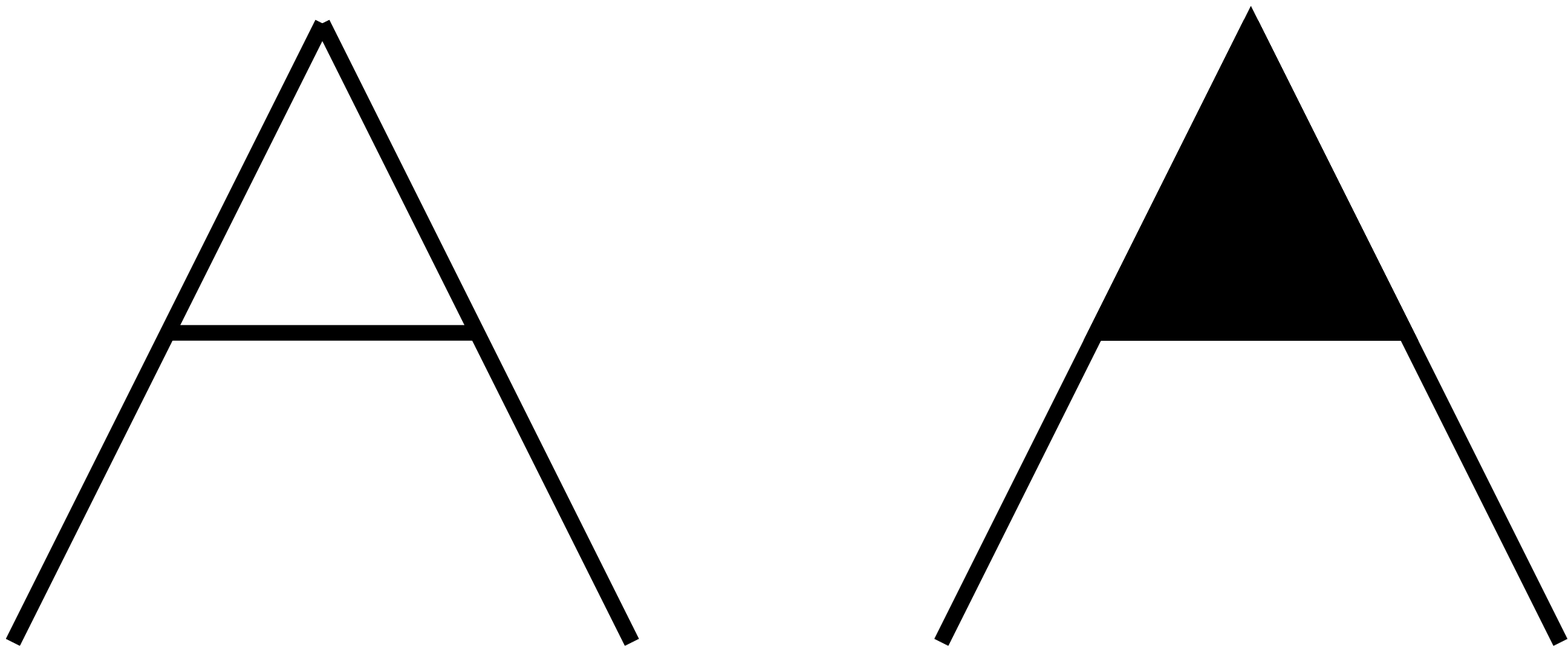}}\negthinspace \negthinspace \negthinspace \negthinspace
& $0 , 0$ & $0 , 0$ & $0 , 0 $ & $0 , 0$ & $1 ,1$ & $x\lt z,1$ &$1 , \phi(y,z)$ & $x \lt z , \phi(y,z)$ \\ \hline \hline
\negthinspace \negthinspace \negthinspace \mbox{
\epsfxsize=.7cm
\epsfbox{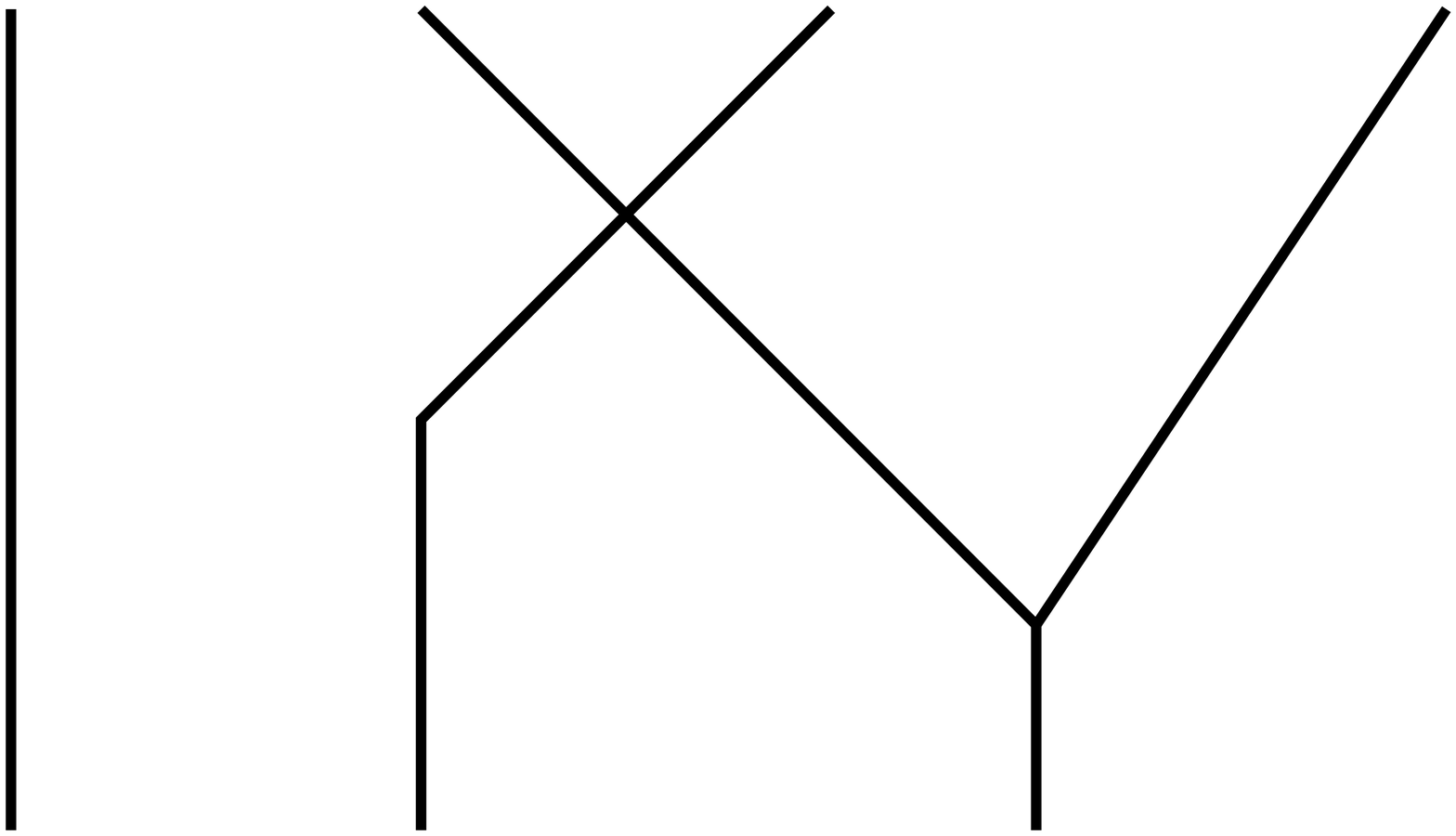}
} \negthinspace \negthinspace  \negthinspace \negthinspace & $1 11 1$ &
 $x 11 1$ &
 $1 1y 1$ &
 $x  1y  1$ &
 $1  z1  z$ &
 $x z1z$ &
 $1 zyz$ &
 $x zy z$ \\ \hline \hline
 &$ 1 \negthinspace  \otimes  \negthinspace  1 \negthinspace  \otimes  \negthinspace  1$ & $x     \negthinspace  \otimes  \negthinspace 1 \negthinspace  \otimes   \negthinspace 1$ & $1  \negthinspace   \otimes    \negthinspace y  \negthinspace \otimes    1$ & $x  \negthinspace   \otimes  \negthinspace   y \otimes  \negthinspace    1$ & $1   \negthinspace \otimes 1     \negthinspace \otimes    \negthinspace  z$ & $x   \negthinspace  \otimes    \negthinspace 1 \otimes   \negthinspace   z$ &$1  \negthinspace  \otimes \negthinspace    y    \negthinspace \otimes  \negthinspace   z$ & $x  \negthinspace   \otimes  \negthinspace   y     \negthinspace \otimes  \negthinspace   z$ \\ \hline \hline \hline
\end{tabular}

\vspace*{1cm}

Thus the calculation becomes:

\begin{eqnarray*}
\lefteqn{ d^{2,1} (\hat{\phi}, 0)\left( \ (a+\sum_x a_x x)\otimes (b+\sum_y a_y y) \otimes (c+\sum_z a_z z)\  \right) }  \\
&= & d^{2,1} (\hat{\phi}, 0) \left(
 abc (1 \otimes 1 \otimes 1) + Abc (x \otimes 1 \otimes 1) + aBc (1 \otimes y \otimes 1) + ABc (x \otimes y \otimes 1)    \right. \\
&&   \left. \quad  \quad+  \ (abC (1 \otimes 1 \otimes z) + AbC (x \otimes 1 \otimes z) + aBC (1 \otimes y \otimes z) + ABC (x \otimes y \otimes z) \right) \\
& =&
 \left[  \ 2 aBC +ABC  \ \phi( x, y)
+  ABC  \ \phi(x\lt y, z) \  \right] \\
& &  -
 \left[ \
2 aBC + ABC  \phi(x\lt z, y \lt z)
 + ABC   \phi(x, z) \ \right] \quad \Box
\end{eqnarray*}

\begin{remark} {\rm
On the other hand, without the factor $k$ in $W$, the original
$2$-cocycles do not give rise to shelf cocycles. Consider $V$ to
have as its basis the trivial quandle $X$ and let $q: V \otimes V
\rightarrow V$ be induced from $\lt$ so that $q(x \otimes y)=x$
for  all $x,y \in X.$ If $\eta_2=0$ and $\eta_1$ is any linear
function, then $d^{2,1}(\eta_1, 0)(x \otimes y \otimes z) = -x \ne
0 \in V$. But in quandle cohomology {\em any function} is a
cocycle.
 } \end{remark}

\begin{theorem}
\label{nontriv2prop} For the cocycles in %%%Proposition
Theorem~\ref{q2prop},
the following holds: If $\phi$ is not a coboundary, then
$\hat{\phi}$  is not a coboundary. In particular, if $H^2_{\rm
Q}(X; k)\neq 0$, then 
%%%%%$H^2,_{\rm sh}(W;W)\neq 0$.
$H^{2,1}_{\rm sh}(W;W)\neq 0$. %%%%%
\end{theorem}
{\it Proof.\/} A function $\phi$ is a coboundary if and only if
there is a $1$-cochain such that $\delta g=\phi$, which is written
as $\phi(x,y)=g(x) - g(x\lt y) $ for any $x, y \in X$ (see
\cite{CJKLS}).

Suppose $\hat{\phi}$ is a coboundary, then there is a $1$-cochain
$f$ such that $D_1 (f)=\hat{\phi}$. A $1$-cochain $f$, in this
case, is a map $f: W \rightarrow W(=k \oplus kX)$, which is
written as
$f(a+\sum a_x x)=f_0(a+ \sum a_x x) + f_1(a+ \sum a_x x)$, where
$a\in k$, $x\in X$, $f_0(a+\sum a_x x) \in k$, and $f_1(a+\sum a_x
x) \in kX$. The condition $D_1(f)=\hat{\phi}$, then, is written
as:
\begin{eqnarray*}
\lefteqn{ \hat{\phi}( \ (a+\sum a_x x) \otimes (b+\sum b_y y) \  )
% = a\sum b_y  + \hat{\phi }(\sum a_x x, \sum a_y y) %% <-- Keep
 = a \sum b_y + \sum_{x, y} a_x b_y\  \phi ( x,  y) }\\
 &=& D_1 (f)(\  (a+\sum a_x x) \otimes (b+\sum a_y y )\  )\\
&=&\{ q(1 \otimes f)-fq +q(f \otimes 1) \} (\  (a+\sum a_x x) \otimes (b+\sum a_y  y)\  )
% Keep:
%\\&=& q(\  (a+\sum a_x x) \otimes ( f_0(b+\sum a_y y) + f_1(b+\sum a_y y) )\  )\\
%& & - f(a + \sum a_x b_y (x\lt y) )
%+ q( \ ( f_0(a+\sum a_x x) + f_1(a+\sum a_x x) ) \otimes (b+\sum b_y y)\  ) .
\end{eqnarray*}
In particular, for $(a + \sum a_x x, b+ \sum b_y y)=(x,y)$, we
obtain:
$$
\phi(x,y)=(x\lt f_1(y) ) - ( f_0(x\lt y) + f_1(x\lt y) ) + ( f_0(x) +  f_1(x) \lt  y),
$$
and by comparing the $k$ and $kX$ factors, this  reduces to
$\phi(x,y)= f_0(x) - f_0(x\lt y)$ and $ f_1(x\lt y)=x\lt f_1(y)+
f_1(x) \lt  y$. In particular, the first equation implies that
$\phi$ is a coboundary and causes a contradiction. $\Box$

Next we consider $3$-cocycles.

\begin{theorem}
\label{q3prop}
For a quandle $3$-cocycle $\theta$
with the coefficient group $A=k$,
define $\hat{\theta}: W \otimes W \rightarrow W$ by
linearly extending $\hat{\theta} ( x \otimes y \otimes z )=\theta(x,y, z)$,
$\hat{\theta} (1 \otimes y \otimes z )=1$, and \\
$\hat{\theta}(x \otimes y\otimes 1 )=\hat{\theta}(x \otimes 1\otimes z)=
\hat{\theta}(x \otimes 1\otimes 1)=\hat{\theta}(1 \otimes y\otimes 1)=
\hat{\theta}(1 \otimes 1\otimes z )=\hat{\theta}(1 \otimes 1\otimes 1)=0$\\
for $x, y, z\in X$. Then $\hat{\theta}$ is a shelf
$3$-cocycle: $d^{3,1} (\hat{\theta}, 0, 0) = 0 $.
\end{theorem}
{\it Proof.\/}
In a manner similar to the proof of Theorem~\ref{q2prop}, we begin by expanding:
\begin{eqnarray*}
\lefteqn{
(a+\sum a_x x ) \otimes (b + \sum b_y y )
 \otimes (c + \sum c_z z )\otimes (d + \sum d_w w) } \\
& = & (a+A x ) \otimes (b + B y )
 \otimes (c + C z )\otimes (d + D w)
 \\
&= & abcd (1 \otimes 1 \otimes 1 \otimes 1) +  Abcd (x \otimes 1 \otimes 1 \otimes 1)  \\ & +&   aBcd (1 \otimes y \otimes 1 \otimes 1)+  ABcd (x \otimes y \otimes 1 \otimes 1)
 \\
&+&  abCd (1 \otimes 1 \otimes z \otimes 1) +  AbCd (x \otimes 1 \otimes z \otimes 1)  \\ & +&   \ aBCd (1 \otimes y \otimes z \otimes 1) +  ABCd (x \otimes y \otimes z \otimes 1)  \\ &+ &
  abcD (1 \otimes 1 \otimes 1 \otimes w) +  AbcD (x \otimes 1 \otimes 1 \otimes w)   \\ &+& aBcD (1 \otimes y \otimes 1 \otimes w) +  ABcD (x \otimes y \otimes 1 \otimes w)  \\ &+&
abCD (1 \otimes 1 \otimes z \otimes w) +  AbCD (x \otimes 1 \otimes z \otimes w)  \\ &+&   aBCD (1 \otimes y \otimes z \otimes w) +  ABCD (x \otimes y \otimes z \otimes w)  \end{eqnarray*}

In a table similar to the one above, the values of the various operators
$q (\hat{\theta}\otimes 1)$ and so forth can be evaluated on each of the sixteen tensors
$(1 \otimes 1 \otimes 1 \otimes 1 )$ through $( x \otimes  y\otimes z\otimes w ).$ Most of these evaluations give $0$ (a result we leave to the reader). The exceptions are the values on
$(1 \otimes y \otimes z \otimes w)$ and $( x \otimes  y\otimes z\otimes w ).$
We remind the reader that $\xi_2=0$ and $\xi_3=0,$ so those terms do not appear below.

\begin{figure}[htb]
\begin{center}
\mbox{
\epsfxsize=2in
\epsfbox{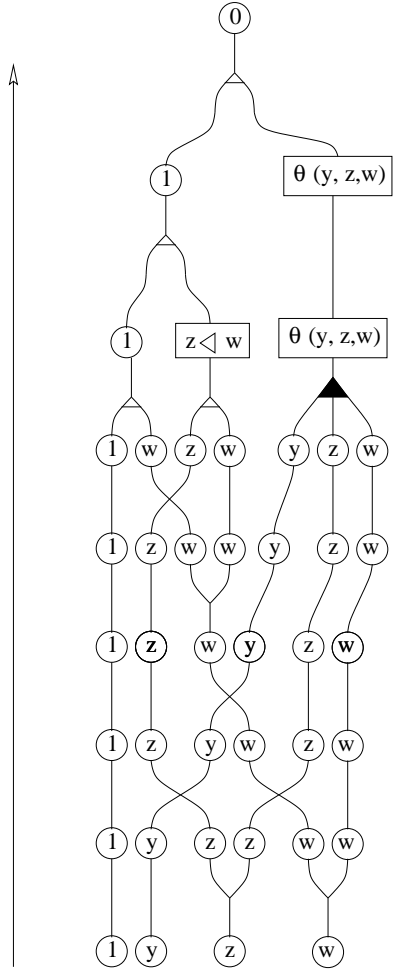}
}
\end{center}
\caption{A sample computation with a $3$-cocycle}
\label{christmastree}
\end{figure}

We compute:
$$ q(\hat{\theta}\otimes 1) (1 \otimes  y \otimes z \otimes w)$$
$$=\hat{\theta} (q\otimes q)(1 \otimes \tau \otimes 1^2)(1^2 \otimes \Delta 1)(1 \otimes  y \otimes z \otimes w)$$
$$=q(\hat{\theta} \otimes q)(1^2\otimes \tau \otimes 1)(1^2 \otimes q \otimes \Delta)(1\otimes \tau \otimes 1^2)(1^2 \otimes \Delta \otimes 1) (1 \otimes  y \otimes z \otimes w)$$
$$=(\hat{\theta})(q\otimes 1^2)(1 \otimes  y \otimes z \otimes w)$$
$$=q(\hat{\theta}\otimes q)(1^2 \otimes \tau \otimes 1)(1^3 \otimes \Delta)      (1 \otimes  y \otimes z \otimes w)$$
$$= \hat{\theta}(q\otimes q\otimes 1)(1 \otimes \tau \otimes 1 \otimes q)(1^2 \otimes \Delta\otimes 1^2)(1^2\otimes \tau \otimes 1)(1^3\otimes \Delta)      (1 \otimes  y \otimes z \otimes w)$$
$$=1,$$
and
$$q(q \otimes \hat{\theta} (q\otimes q\otimes 1^3)(1 \otimes \tau \otimes 1^4)(1^2 \otimes \Delta \otimes 1^3)(1^2 \otimes \tau \otimes 1^2)(1\otimes \tau \otimes \tau \otimes 1)(1^2 \otimes \Delta \otimes \Delta ) (1 \otimes  y \otimes z \otimes w)$$
$$= q (q \otimes q)(1 \otimes \tau \otimes1)(q \otimes q \otimes \xi_2)(1 \otimes \tau 1^3)(1^2 \otimes \Delta \otimes 1^2)(1^2 \otimes \tau \otimes 1)(1^3 \otimes \Delta)(1 \otimes  y \otimes z \otimes w)$$
$$=0.$$
The last equality follows trivially since $\xi_2=0$.  A scheme for making these computations is illustrated in Fig~\ref{christmastree}. Meanwhile,
%%%hats moved to theta twice and extra parens rm
%%%also an \otimes added at second eqn.
$$ q(\hat{\theta}\otimes 1) (x \otimes  y \otimes z \otimes w)= \theta(x,y,z)$$
$$\hat{\theta} (q\otimes q)(1 \otimes \tau \otimes 1^2)(1^2 \otimes \Delta \otimes1)(x \otimes  y \otimes z \otimes w)=\theta(x\lt z, y \lt z, w)$$
$$q(\hat{\theta} \otimes q)(1^2\otimes \tau \otimes 1)(1^2 \otimes q \otimes \Delta)(1\otimes \tau \otimes 1^2)(1^2 \otimes \Delta \otimes 1) (x \otimes  y \otimes z \otimes w)=\theta(x,z,w)$$
$$(\hat{\theta})(q\otimes 1^2)(x \otimes  y \otimes z \otimes w)=\theta(x\lt y,z,w)$$
$$q(\hat{\theta}\otimes q)(1^2 \otimes \tau \otimes 1)(1^3 \otimes \Delta)      (x \otimes  y \otimes z \otimes w)=\theta(x,y,w)$$
$$ \hat{\theta}(q\otimes q\otimes 1)(1 \otimes \tau \otimes 1 \otimes q)(1^2 \otimes \Delta\otimes 1^2)(1^2\otimes \tau \otimes 1)(1^3\otimes \Delta)      (x \otimes  y \otimes z \otimes w)=\theta(x\lt w, y \lt w,z\lt w),$$
and
$$q(q \otimes \hat{\theta} (q\otimes q\otimes 1^3)(1 \otimes \tau \otimes 1^4)(1^2 \otimes \Delta \otimes 1^3)(1^2 \otimes \tau \otimes 1^2)(1\otimes \tau \otimes \tau \otimes 1)(1^2 \otimes \Delta \otimes \Delta ) (x \otimes  y \otimes z \otimes w)$$
$$= q (q \otimes q)(1 \otimes \tau \otimes1)(q \otimes q \otimes \xi_2)(1 \otimes \tau 1^3)(1^2 \otimes \Delta \otimes 1^2)(1^2 \otimes \tau \otimes 1)(1^3 \otimes \Delta)(x \otimes  y \otimes z \otimes w)$$
$$=0.$$

The result follows. $\Box$

\begin{theorem}
\label{quan3}
For the cocycles in
Theorem~\ref{q3prop}, the following holds:
If $\theta$ is not a coboundary, then $\hat{\theta}$  is not a
coboundary. In particular, if $H^3_{\rm Q}(X; k)\neq 0$, then
$H^{3,1}_{\rm sh}(W;W)\neq 0$. %% $H^3_{\rm sh}(W;W)\neq 0$. %%MS0221
\end{theorem}
{\it Proof.\/} The proof is similar to that of
Theorem~\ref{nontriv2prop}. The cochain $\theta$ is a
coboundary if and only if there is a $2$-cochain $\phi$ such  that
$\delta \phi=\theta$,
 which is written as $\theta(x,y,z)=\phi(x,y)+\phi(x\lt
y, z) - \phi(x,z) - \phi(x\lt z, y \lt z) $
 for any $x, y, z \in X$ (see \cite{CJKLS}).

Suppose $\hat{\theta}$ is a coboundary. Then there is a
$2$-cochain $f$ such that $D_2 (f)=\hat{\theta}$. A $2$-cochain
$f$, in this case, is a map $f: W \otimes W \rightarrow W(=k
\oplus kX)$, that is written as:
$$f(\ ( a+\sum a_x x) \otimes (b + \sum b_y y) \ )
=f_0(\ ( a+ \sum a_x x )\otimes (b + \sum b_y y) \ )
+ f_1(\ ( a+ \sum a_x x)\otimes (b + \sum b_y y) \ ) , $$
where $a\in k$, $x, y\in X$, $f_0(a+\sum a_x x) \in k$, and $f_1(a+\sum a_x x) \in kX$.
We take specific values and compute
$\hat{\theta}  = D_2 (f)  $ evaluated at
$x \otimes y \otimes z$.
We have $\hat{\theta} ( x \otimes y \otimes z ) =\theta(x,y,z) $,  and
\begin{eqnarray*}
 D_2 (f) ( x \otimes y \otimes z )&=&
 [\   f_0(x \otimes y) + f_1(x \otimes y) \lt z
 + f_0(( x \lt y )\otimes z) + f_1( ( x \lt y )\otimes z) \ ] \\
 &-& [ \ f_0(  ( x \lt z)\otimes (y \lt z) ) +  f_1(  ( x \lt z)\otimes (y \lt z) )
+
f_0(x \otimes z) \\
& & + f_1(x \otimes z) \lt (y \lt z) + (x \lt z)\lt f_1(y \otimes z) \ ] ,
\end{eqnarray*}
and comparing the elements on $k$, we obtain
$$\theta (x,y,z) =   f_0(x \otimes y)
 + f_0(( x \lt y )\otimes z)  -  f_0(  ( x \lt z)\otimes (y \lt z) )
 -  f_0(x \otimes z), $$
 so that by defining $\phi(x, y)=f_0(x \otimes y)$ for any $x,y \in X$,
 we obtain a contradiction $\theta=\delta \phi$.
 $\Box$

\subsection{Lie Algebra Cohomology}

Let $q:N \otimes N \rightarrow N$ be the map defined in
Lemma~\ref{Liesec},
where $N=k \oplus \g$ for a Lie algebra $\g$ over a ground field $k$.
Let $\psi: \g \times \g \rightarrow \g$ be a Lie algebra $2$-cocycle,
with adjoint action. Then $\psi $ is bilinear and satisfies
$$\psi(y , x)=-\psi(x , y) , $$
$$  [ \psi(x,y), z] + [ \psi(y,z), x] + [ \psi(z,x), y ] +
\psi( [x,y], z)+ \psi( [ y, z], x ) + \psi( [ z, x], y) = 0 .
$$
It defines a linear map $\psi: \g \otimes \g \rightarrow \g$. The following
result says that a Lie algebra $2$-cocycle gives rise to a shelf  $2$-cocycle,
 when the comultiplication is fixed and undeformed
($\eta_2=0$).
\begin{theorem}
\label{Lie2ad}
Let $\psi: \g \times \g \rightarrow \g$ be a Lie algebra
$2$-cocycle with adjoint action. Define $\hat{\psi} : N \otimes
N \rightarrow N $ by $\hat{\psi} ((a+x) \otimes (b+y)) = \psi(x
\otimes y)$ for $a,b,c \in k$, $x,y,z \in \g$. Then $\hat{\psi}$
is a shelf $2$-cocycle: $d^{2,1}
(\hat{\psi}, 0)= d^{2,2} (\hat{\psi}, 0) = 0.$ %%%%
%%%%\tilde -->\hat twice above 
%%%Also d^{2,2}
\end{theorem}
{\it Proof.\/} One computes:
\begin{eqnarray*}
\lefteqn{ d^{2,1} (\hat{\psi}, 0) (\ (a + x)\otimes (b+y) \otimes (c+z) \ ) }\\
%%%(c+z) added
&=&
\{  \ q(\ \psi(x, y) \otimes (c+z) \ ) + \hat{\psi}(\ (ab+bx+[x,y]) \otimes (c+z) \ ) \}  \\
&-& \{  \ \{ \ \hat{\psi}(q \otimes q) + q (\hat{\psi} \otimes q) + q( q \otimes \hat{\psi})\
 \} \tau_2
 (\ (a+x)\otimes (b+y) \otimes ( \ (c+z) \otimes 1 + 1 \otimes z\  ) \ ) \ \}  \\
 &=& \{  \ ( c \psi(x,y) + [\psi(x,y), z ]  \ ) + ( \ b \psi(x,z) + \psi([x,y], z) \ ) \ \}  \\
 &-& \{ \ (\ c \psi(x,y) + \psi( [x,z], y) + \psi (  x, [y,z])  \ ) \\
& &  + ( \ b \psi(x,z) + [ \psi (x,z) , y ] \ ) + ( \ [ x, \psi(y,z) ]\ ) \ \} \quad =0. %period here 
\end{eqnarray*}
%and the result follows. 
The other equality $d^{2,2} (\hat{\psi}, 0) = 0$ is checked similarly. %% MS0221
$\Box$

Next we consider  Lie algebra  $2$-cocycles
$\psi : \g \times \g \rightarrow k$ with the trivial
representation on the ground field $k$.
In this case the $2$-cocycle condition is being skew-symmetric
and satisfying the Jacobi identity:
\begin{eqnarray*}
\psi (y , x)&=& -\psi (x , y) ,\\
\psi ( [x,y], z)+ \psi ( [ y, z], x ) + \psi ( [ z, x], y) &=&  0 .
\end{eqnarray*}
Let $\g '= k \gamma + \g$ where $\gamma \in \g$ and $[\gamma,z]=0$ for all $z\in \g$.  Then $\g'$ is a
Lie algebra with Lie bracket given by $[ a\gamma + x , b\gamma + y
]' = [x, y] \in \g \subset \g '$.
For a given $2$-cocycle $\psi: \g \times \g \rightarrow k$, define
$\psi': \g ' \times \g ' \rightarrow \g '$ by $\psi' (a\gamma + x
, b\gamma + y) = \psi(x,y) \gamma \in \g '$. Then we claim that
$\psi'$ satisfies the $2$-cocycle condition with adjoint action.
We compute:
$$[ \psi' (a \gamma+x, b \gamma+y ), c \gamma+z ]'=
 [\psi (x, y) \gamma, c \gamma+z ]' = [\psi (x, y) \gamma, z]=0. $$
Therefore the first three terms involving the adjoint action, in
fact, vanish by construction. The last three terms reduce to the
$2$-cocycle condition of $\psi$, since
$$\psi' ( [a  \gamma + x, b  \gamma + y]', c  \gamma +z)
= \psi ' ([x,y], c \gamma + z)=\psi( [x, y], z). $$ Hence this
reduces to the previous case. We summarize this situation as:

\begin{theorem}
\label{Lie2} A Lie algebra $2$-cocycle valued in the ground field
with trivial representation gives rise to a shelf $2$-cocycle.
\end{theorem}

Next we investigate relations for $3$-cocycles. A Lie algebra
$3$-cocycle with adjoint  action is a totally skew-symmetric
trilinear map $\zeta : \g \times \g  \times \g \rightarrow \g $
for a Lie algebra $\g$ that satisfies
\begin{eqnarray*}
[ \zeta (x,y,z), w] - [\zeta (x,y,w), z]+ [ \zeta (x,z,w), y] - [\zeta (y,z,w), x] & & \\
- \  \zeta ( [x,y],z,w) + \zeta( [x,z],y,w) - \zeta( [x,w], y,z) & & \\
+ \  \zeta( [y,z],x,w) - \zeta ([ y,w], x, z) + \zeta ([z,w], x,y)
& = & 0 .
\end{eqnarray*}
This defines a linear map $\zeta: \g \otimes \g  \otimes \g \rightarrow \g $.
Recall that we defined $N=k \oplus \g $.

\begin{theorem}
\label{Lie3}
Let $\zeta: \g \times \g\times \g  \rightarrow \g$ be a Lie
algebra $3$-cocycle with adjoint action. Define $\hat{\zeta} : N
\otimes N \otimes N \rightarrow N $ by $\hat{\zeta} ( (a+x)
\otimes (b+y) \otimes (c+z) ) = \zeta(x \otimes y \otimes z)$.
Then $\hat{\zeta}$ 
%% is a shelf $3$-cocycle: 
satisfies %% MS0221
$d^{3,1} (\hat{\zeta}, 0, 0) = 0 $.
\end{theorem}
{\it Proof.\/}
There are  four positive ($L_1, L_2, L_3, L_4$)
and three negative ($R_1, R_2, R_3$) terms in $d^{3,1} (\hat{\zeta}, 0, 0)$
(the last negative term vanishes because $\xi_2=0$ in
$(\xi_1, \xi_2, \xi_3)= (\hat{\zeta}, 0, 0)$).
We evaluate each term for a general element
$$(a+x)\otimes (b+y) \otimes (c+z) \otimes (d+w) $$
as before.
The first term $L_1$ is
$$ q( \ \zeta(x, y, z) \otimes (d + w) \ ) = d \ \zeta(x,y,z) + [ \zeta(x,y,z), w] .$$
The second term $L_2$ is
\begin{eqnarray*}
\lefteqn{ \hat{\zeta} (q \otimes q \otimes 1) ( \
(a+x) \otimes \{ (c+z) \otimes (b+y) \otimes 1 + 1 \otimes (b+y) \otimes z \} \otimes w
\  ) } \\
& = &  \hat{\zeta} ( \
(ac + cx + [x,z] ) \otimes (b+y) \otimes w + (a+x) \otimes [y,z] \otimes w  \ ) \\
&=& c \ \zeta(x,y,w) + \zeta([x,z],y,w) + \zeta(x, [y,z],w)
\end{eqnarray*}
By similar calculations the remaining terms give
\begin{eqnarray*}
L_3 & : &  b\ \zeta (x,z,w) + [\zeta( x,z,w), y]  \\
L_4 & : &  [x, \zeta(y,z,w)] \\
R_1 & : &  b \ \zeta (x,z,w) + \zeta( [x,y],z,w) \\
R_2 & : & c\ \zeta (x,y,w) +[ \zeta (x,y,w), z] \\
R_3 & : &  d\  \zeta(x,y,z) + \zeta([x,w],y,z)+\zeta (x, [y,w],z) + \zeta( x,y,[z,w]) \\
\end{eqnarray*}
and the result $(L_1+L_2+L_3+L_4)-(R_1+R_2+R_3)=0$ follows.
$\Box$

\begin{theorem}\label{nontrivLie2prop}
For the cocycles in
Theorem~\ref{Lie2},
 the following
holds: If $\psi$ is not a coboundary, then $\hat{\psi}$  is not a
coboundary. In particular, if the second cohomology group of the
Lie algebra cohomology with adjoint action is non-trivial
$(H^2_{\rm Lie}(\g; \g)\neq 0)$, then $H^2_{\rm sh}(N;N)\neq 0$.
\end{theorem}
{\it Proof.\/} The proof is similar to that of
Theorem~\ref{nontriv2prop}.
If $\psi$ is  a coboundary, then
there is a $1$-cochain $g$ such  that $\delta g=\psi$, which is
written as $\psi(x,y)=[x, g(y)]+[g(x), y] - g([x,y])$ for any $x,
y \in \g$.

Suppose $\hat{\phi}$ is a coboundary, then there is a $1$-cochain
$f$ such that $D_1 (f)=\hat{\psi}$. A $1$-cochain $f$, in this
case, is a linear map $f: N \rightarrow N(=k \oplus \g)$, that is
written as
$f(a+ x)=f_0(a+  x) + f_1(a+
x)$, where $a\in k$, $x\in \g $, $f_0(a+ x) \in k$, and $f_1(a+ x)
\in \g$. The condition $D_1(f)=\hat{\phi}$, then, is written as
\begin{eqnarray*}
\lefteqn{ \hat{\psi}( \ (a+ x) \otimes (b+y) \  )
 =   \psi ( x,  y) }\\
 &=& D_1 (f)(\  (a+ x) \otimes (b+ y )\  )\\
&=&\{ q(1 \otimes f)-fq +q(f \otimes 1) \} (\  (a+ x) \otimes (b+ y)\  )
\end{eqnarray*}
In particular, for $(a +  x, b+  y)=(x,y)$, we obtain
\begin{eqnarray*}
\psi(x,y) &=& q(x \otimes  f_1(y) )  - ( f_0( q(x \otimes  y))
+ f_1( q (x\otimes  y) ) ) + ( f_0(x) +  q( f_1(x) \otimes   y) ) \\
&=&
[x,  f_1(y)] - f_0([x,y]) - f_1([x,y]) + f_0(x) + [f_1(x), y].
\end{eqnarray*}
Comparing the elements in $k$ and $\g$ in the image,
we obtain
\begin{eqnarray*}
0 &=& - f_0([x,y]) + f_0(x) , \\
\psi(x,y) &=&[x,  f_1(y)] - f_1([x,y]) + [f_1(x), y],
\end{eqnarray*}
and the second implies that $\psi$ is a coboundary.
$\Box$

Let $W_p$ be the Witt algebra, a Lie algebra over the field
${\mathbb F}_p$ with $p$ elements for a prime $p>3$. Specifically,
$W_p$ has basis $e_a$, $a\in {\mathbb F}_p$ and has bracket
defined by $[e_a, e_b]=(b-a)e_{a+b}$. Then it is
known~\cite{Block} (we thank J. Feldvoss for informing us) that
the Lie algebra cohomology with trivial action
 $H^2_{\rm Lie}(W_p; {\mathbb F}_p)$ is one-dimensional and generated by
 the Virasoro cocycle $c(e_a, e_{-a})=a(a^2-1)$ (otherwise zero).
Let  $W_p'=k\gamma \oplus W_p$, $N(W_p')=k \oplus W_p'$ be the
%% MS %% self-distributive object in $\Coalg$
object in $\CoComCoalg$ with a self-distributive linear map $q$ %% MS
constructed in Section~\ref{rackshelfsec}.  Then we have:
\begin{corollary}
$H^2_{\rm sh}(N(W_p'); N(W_p'))\neq 0$.
\end{corollary}

\section{A Compendium of Questions}

What are more precise relationships among the Lie bracket,
self-distributivity, solutions to the Yang-Baxter equations,  Hopf
algebras, and quantum groups? Can the cocycles constructed herein
be used to construct invariants of knots and knotted surfaces? Can
the coboundary maps be expressed skein theoretically? Is there a
spectral sequence that is associated to
a % the %% MS
filtration of the
chain groups? If so, what are the differentials? What does it
compute? Are there non-trivial cocycles among any of the
trigonometric shelves? The proofs of the main theorems come from
grinding through computation. Are there more conceptual  proofs?
How can the theory be extended to higher dimensions, such as to
higher dimensional Lie algebras, or Lie $2$-algebras?  How, if at
all, do the Zamolodchikov tetrahedron
equation and the Jacobiator
identity of a Lie $2$-algebra, relate to shelf cohomology?
%%What
%%%is the missing corner of the following square and how does this
%%%diagram relate to the work done here?
%%%\[
%%%  \xymatrix{
%%%   {\rm Lie groups}
%%%     \ar[rr]^{}
%%%     \ar[dd]_{}
%%%     && {\rm Lie algebras}
%%%     \ar[dd]^{} \\ \\
%%%   {\rm Quandles}
%%%     \ar[rr]^{}
%%%     && ??? }
%%%\]
Can it be shown to be a cohomology theory in the case when $\xi_2$
and $\xi_3$ are non-zero? Is there a spin-foam interpretation of
the $3$-cocycle conditions?
% Does this give topological applications? %% MS % asked already for knots above %
% Do you have any questions? %% MS %% I suppose we delete this in the final version.

\newpage
\appendix
\section{Proving $D_3 D_2=0$}

The next illustrations represent  the proofs
of Lemmas~\ref{d321lem} and \ref{d322lem}.

\bigskip

\begin{figure}[htb]
\begin{center}
\mbox{
\epsfxsize=2.5in
\epsfbox{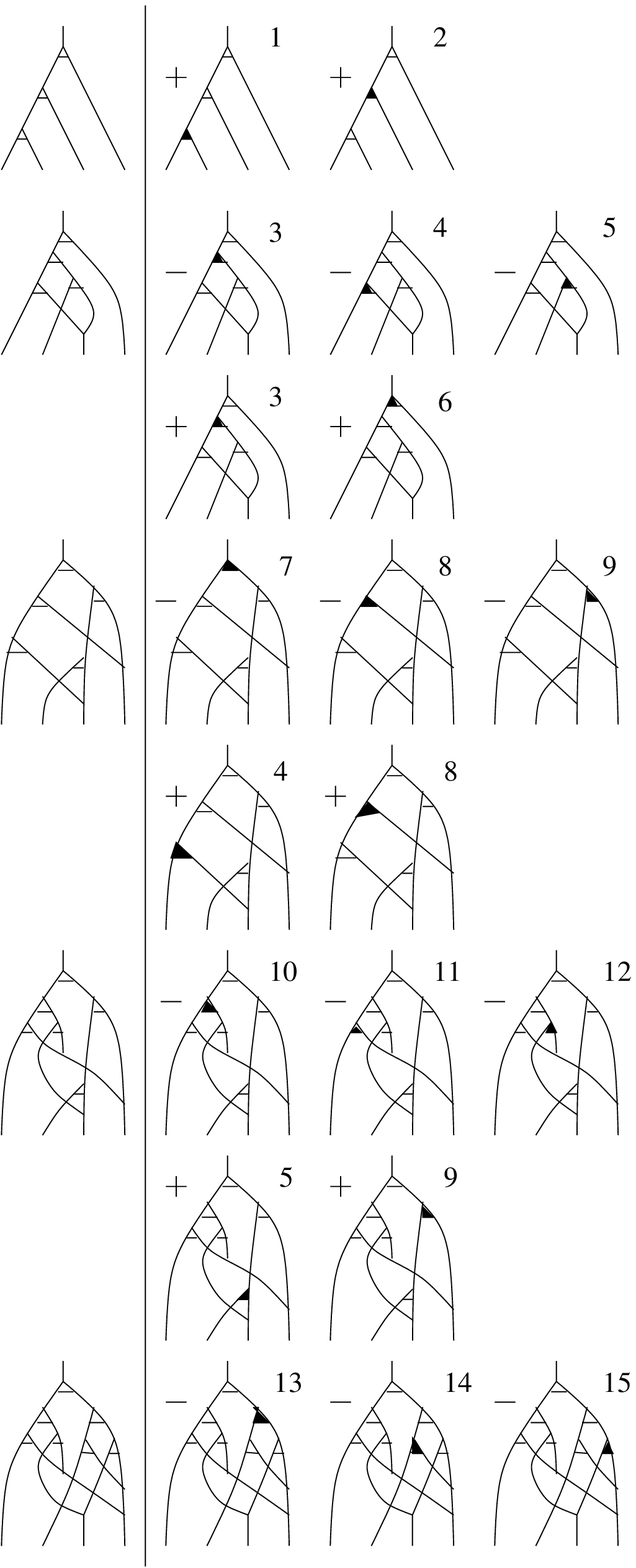}
}
\end{center}
\caption{ $d^{3,1} d^{2,1}(\eta_1, 0)$, LHS}
\label{D31L}
\end{figure}

\begin{figure}[htb]
\begin{center}
\mbox{
\epsfxsize=2.5in
\epsfbox{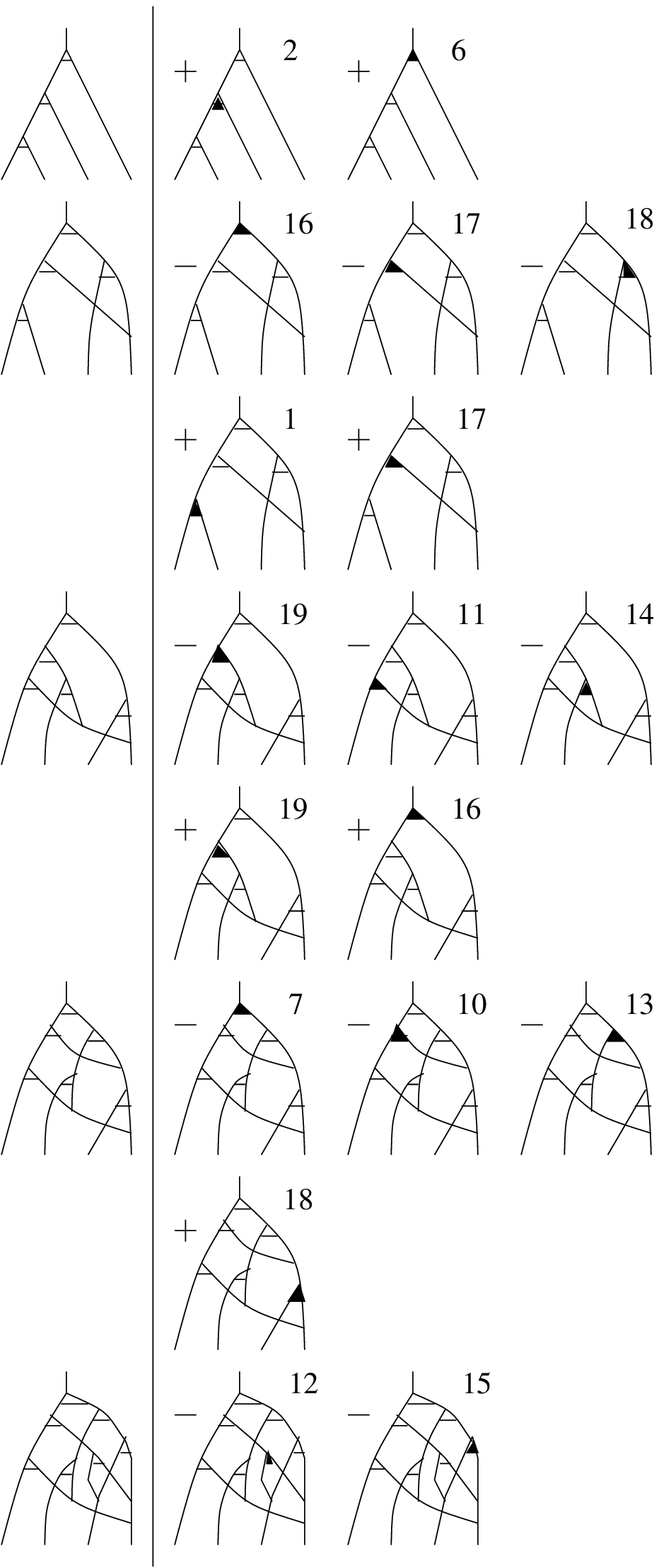}
}
\end{center}
\caption{$d^{3,1} d^{2,1}(\eta_1, 0)$, RHS}
\label{D31R}
\end{figure}

\begin{figure}[htb]
\begin{center}
\mbox{
\epsfxsize=3.2in
\epsfbox{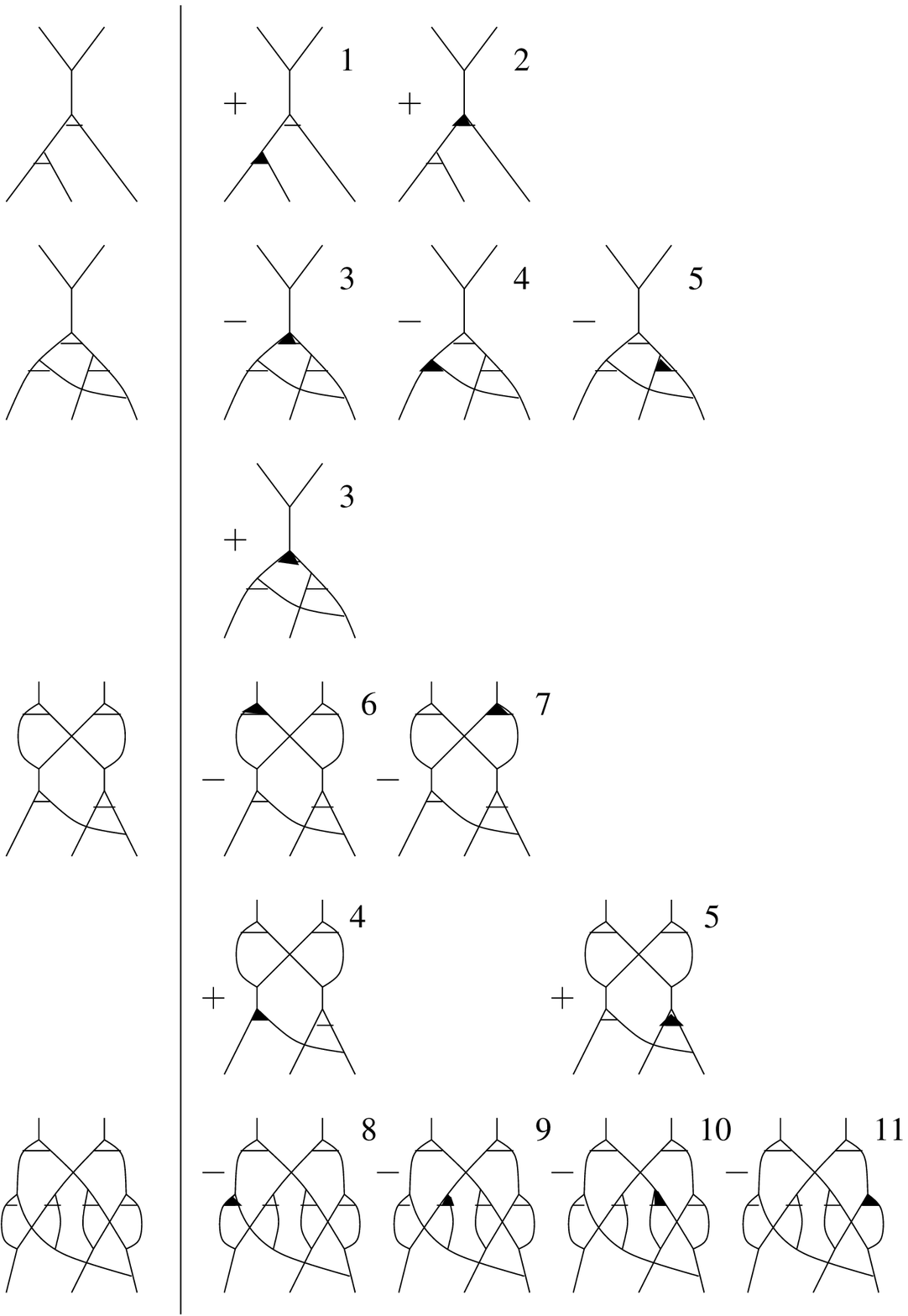}
}
\end{center}
\caption{$d^{3, 2} d^{2, 1} (\eta_1, 0)$, LHS}
\label{D32L}
\end{figure}

\begin{figure}[htb]
\begin{center}
\mbox{
\epsfxsize=3.5in
\epsfbox{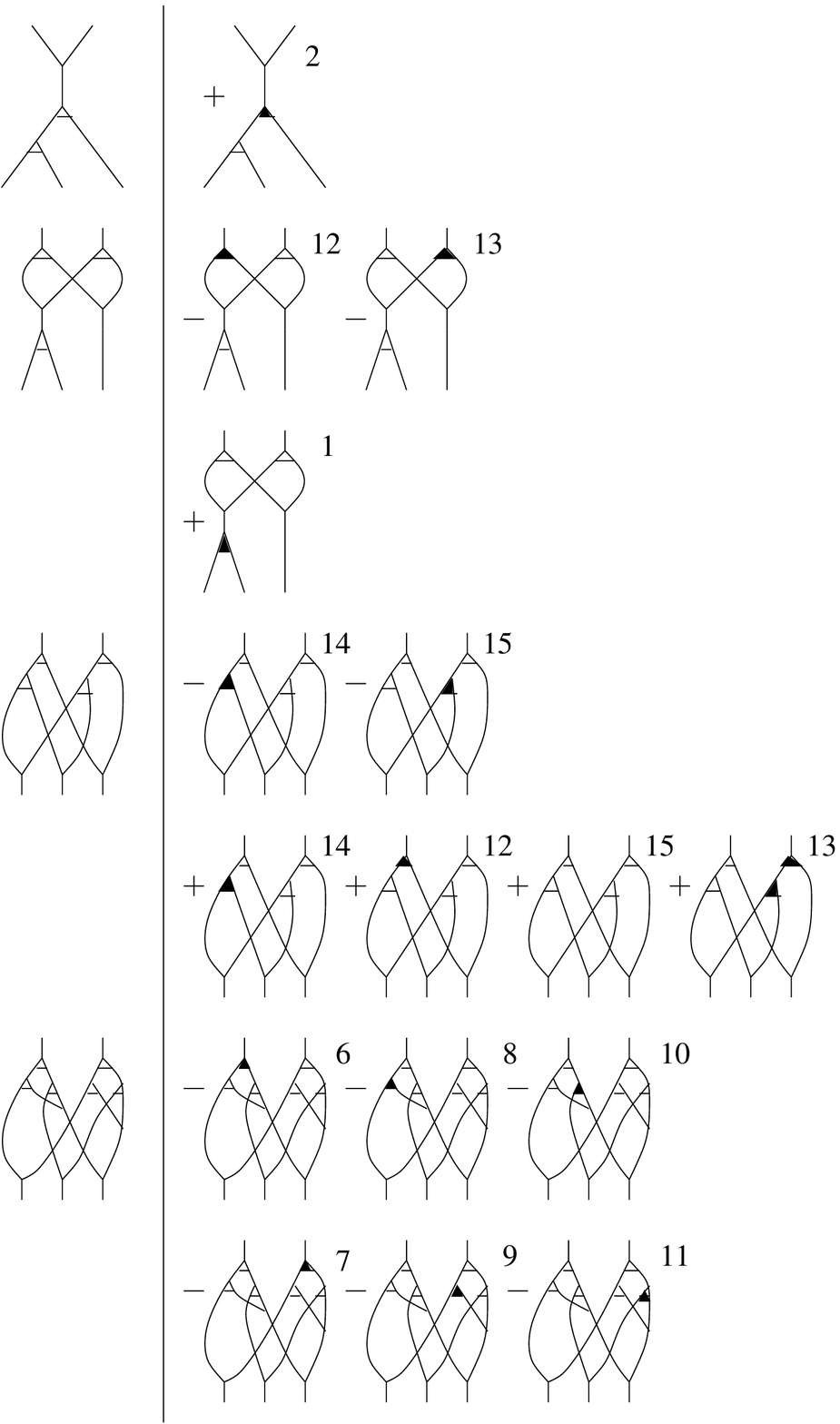}
}
\end{center}
\caption{$d^{3, 2} d^{2, 1} (\eta_1, 0)$, RHS}
\label{D32R}
\end{figure}

\begin{figure}[htb]
\begin{center}
\mbox{
\epsfxsize=4in
\epsfbox{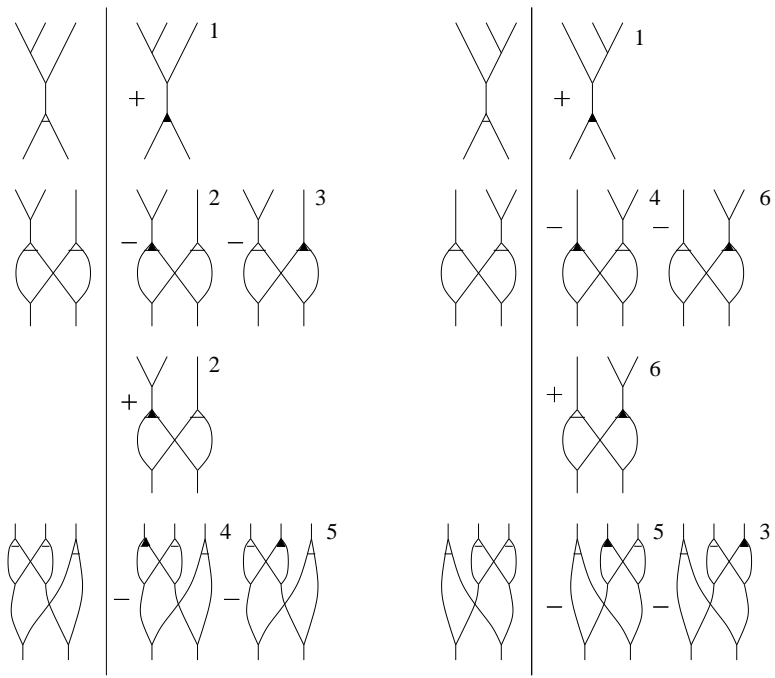}
}
\end{center}
\caption{$d^{3, 3} d^{2, 1} (\eta_1, 0)=0$}
\label{D33L}
\end{figure}

\clearpage

\section{Proving identities between terms in Fig.~\ref{D31L} and \ref{D31R}}

The next illustrations give the outlines of the proofs that the
terms labelled $7$, $10$, $11$, $12$, $13$, and $14$  represent
the same functions in Figs.~\ref{D31L} and \ref{D31R}.

\vspace{1in}

\begin{figure}[htb]
\begin{center}
\mbox{
\epsfxsize=3in
\epsfbox{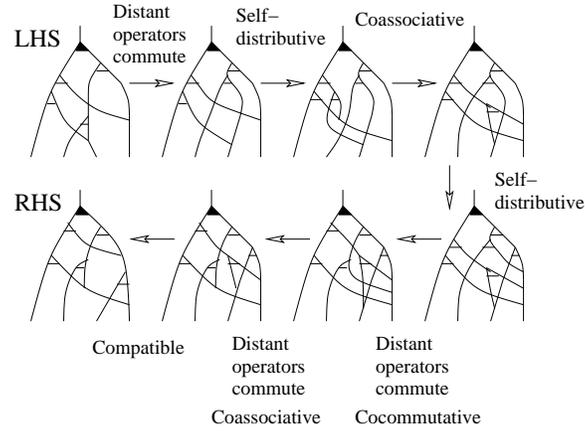}
}
\end{center}
\caption{The term $7$}
\label{trickproof7}
\end{figure}

\begin{figure}[htb]
\begin{center}
\mbox{
\epsfxsize=3in
\epsfbox{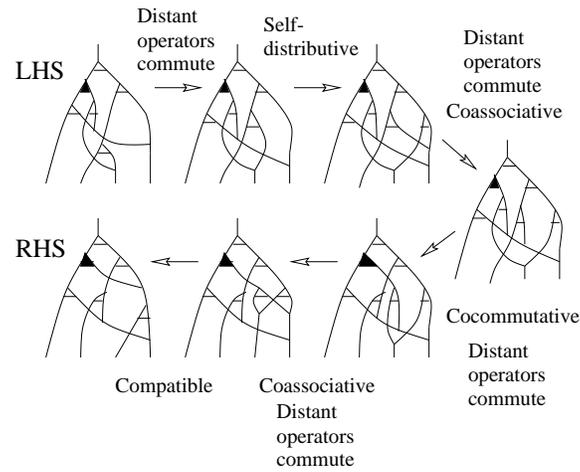}
}
\end{center}
\caption{The term $10$}
\label{trickproof10}
\end{figure}

\begin{figure}[htb]
\begin{center}
\mbox{
\epsfxsize=4in
\epsfbox{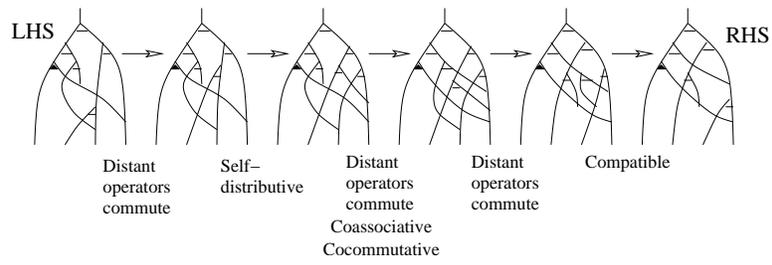}
}
\end{center}
\caption{The term $11$}
\label{trickproof11}
\end{figure}

\begin{figure}[htb]
\begin{center}
\mbox{
\epsfxsize=4in
\epsfbox{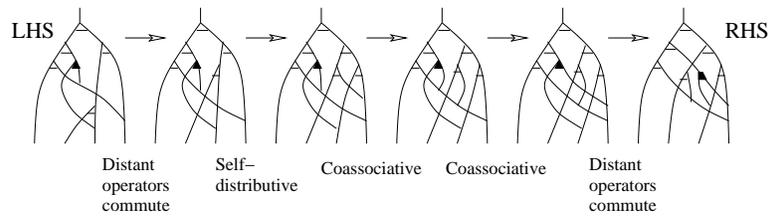}
}
\end{center}
\caption{The term $12$}
\label{trickproof12}
\end{figure}

\begin{figure}[htb]
\begin{center}
\mbox{
\epsfxsize=3.5in
\epsfbox{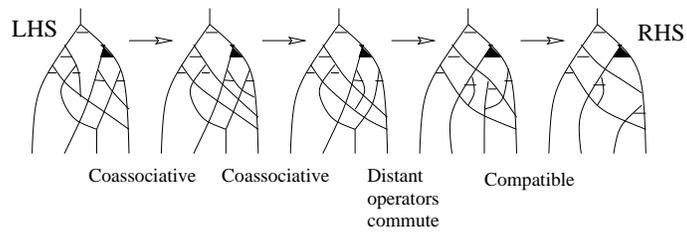}
}
\end{center}
\caption{The term $13$}
\label{trickproof13}
\end{figure}

\begin{figure}[htb]
\begin{center}
\mbox{
\epsfxsize=3.5in
\epsfbox{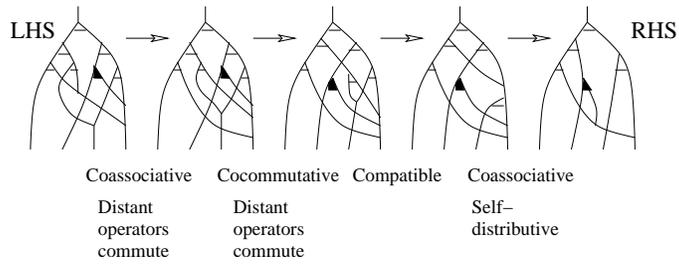}
}
\end{center}
\caption{The term $14$}
\label{trickproof14}
\end{figure}

\end{document}